\documentclass{aims}
\usepackage{amsmath}
\usepackage{paralist}

\usepackage{graphicx,color}
\usepackage[font=footnotesize]{caption}
\usepackage{lscape,float}
\usepackage{authblk}
\usepackage{lineno}
\usepackage{mathtools,amssymb,array,amsthm,amsfonts}
\usepackage[parfill]{parskip}
\usepackage[T1]{fontenc}
\usepackage{bigfoot}
\usepackage[numbered,framed]{matlab-prettifier}
\usepackage[export]{adjustbox}
\usepackage[parfill]{parskip}
\usepackage[font=footnotesize]{subcaption}
\usepackage{lmodern}
\usepackage[utf8]{inputenc}
\usepackage{algorithm,algpseudocode}
\usepackage{dsfont}
\usepackage{bbm}
\usepackage{scalerel,relsize}
\usepackage{booktabs,mathrsfs}
\usepackage{mathptmx}
\usepackage[table]{} 
\usepackage{hhline}
\usepackage{colortbl}
\usepackage{tabularx,ragged2e} 

\usepackage[pdftex, colorlinks]{hyperref}
\hypersetup{colorlinks=true,  	
    linkcolor=red,        	
    citecolor=black,      	
    filecolor=black,      	
    urlcolor=blue         	
    }

\def\currentvolume{X}
 \def\currentissue{X}
  \def\currentyear{200X}
   \def\currentmonth{XX}
    \def\ppages{X--XX}
     \def\DOI{10.3934/xx.xx.xx.xx}

\let\oldReturn\Return
\renewcommand{\Return}{\State\oldReturn}
\newtheorem*{thm}{Theorem}
\newcommand{\vars}{\texttt}
\newcommand{\R}{\mathbb{R}}
\DeclareMathAlphabet{\mathcal}{OMS}{cmsy}{m}{n}
\DeclareMathOperator{\vol}{vol}

\definecolor{S1}{rgb}{0, 0.4470, 0.7410} 
\definecolor{S2}{rgb}{0.8500, 0.3250, 0.0980}
\definecolor{S3}{rgb}{0.9290, 0.6940, 0.1250}
\definecolor{S4}{rgb}{0.4940, 0.1840, 0.5560}
\definecolor{t1}{RGB}{76.7550, 189.9750, 237.9150} 
\definecolor{t2}{RGB}{255, 255, 0}
\definecolor{t3}{RGB}{255, 0 , 0}
\definecolor{t4}{RGB}{118.8300  171.8700   47.9400}
\definecolor{t5}{RGB}{255, 140, 0}
\definecolor{t6}{RGB}{72, 21, 170}
\definecolor{tn1}{RGB}{0, 200, 50} 
\definecolor{tn2}{RGB}{20, 0 , 250}

\graphicspath{{./FigsTo2Vs/}}

\title[A tale of two vortices] 
      {A tale of two vortices: how numerical ergodic
theory and transfer operators reveal fundamental changes to coherent
structures in non-autonomous dynamical systems$^\dagger$}

\author[Chantelle Blachut and Cecilia Gonz\'alez-Tokman]{}


\email{uqcblach@uq.edu.au, cecilia.gt@uq.edu.au}

\thanks{This work has been partially supported by an Australian Research Council Discovery Early Career Researcher Award (DE160100147) and by an Australian Government Research Training Program Stipend Scholarship (CB) }

\thanks{$^\dagger$ To appear in the Journal of Computational Dynamics. $^*$ Corresponding author: uqcblach@uq.edu.au}
\begin{document}
\maketitle

\centerline{\scshape Chantelle Blachut$^*$ and Cecilia Gonz\'alez-Tokman}

{\footnotesize
\centerline{School of Mathematics and Physics,}
\centerline{The University of Queensland,}
\centerline{St Lucia, QLD 4072, Australia}
} 



\vspace{-0.25cm}
\begin{abstract}
Coherent structures are spatially varying regions which disperse minimally over time and organise motion in non-autonomous systems.
This work develops and implements algorithms providing multilayered descriptions of time-dependent systems which are not only useful for locating coherent structures, but also for detecting time windows within which these structures undergo fundamental structural changes, such as merging and splitting events.
These algorithms rely on singular value decompositions associated to Ulam type discretisations of transfer operators induced by dynamical systems, and build on recent developments in multiplicative ergodic theory.
Furthermore, they allow us to investigate various connections between the evolution of relevant singular value decompositions and dynamical features of the system. 
The approach is tested on models of periodically and quasi-periodically driven systems, as well as on a geophysical dataset corresponding to the splitting of the Southern Polar Vortex. 
\end{abstract}
\section{Introduction}
Coherent structures are spatially varying regions that disperse minimally over time and organise motion in non-autonomous systems. In the form of oceanic eddies and atmospheric vortices, they play important roles in biogeophysical phenomena and influence the weather of our planet. Understanding and characterising the dynamical behaviour of such structures, as well as maximising the information about them that can be extracted from data and models of the underlying flows, is important for understanding how transport and mixing properties develop as the dynamical system evolves. In this paper, we consider fundamental structural changes, such as merging and splitting events, to show that the dynamical behaviour of coherent structures can be characterised using transfer operator technology and results from ergodic theory.

\begin{figure}[H]
\centering
\begin{minipage}[b]{\textwidth}
\centering
\begin{minipage}[b]{0.49\textwidth}
\centering
\includegraphics[width=\columnwidth]{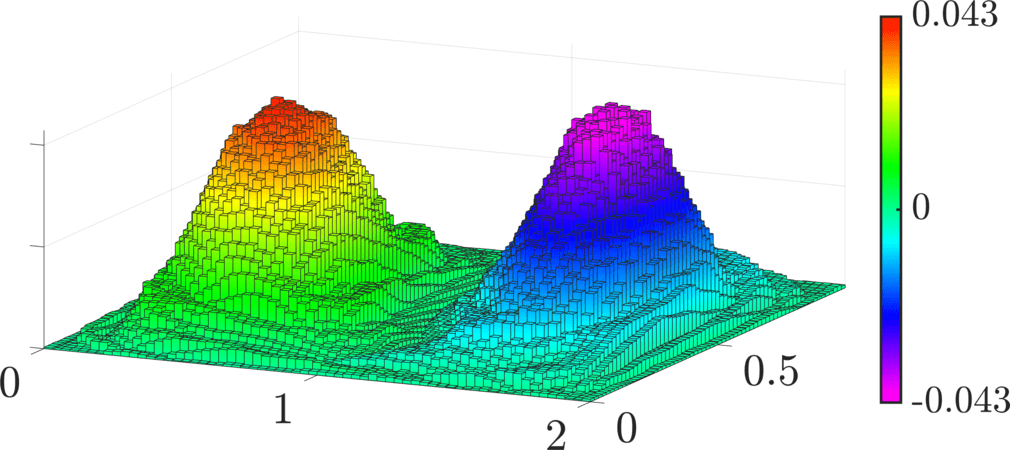}
\subcaption{Periodically driven double gyre at time $t=0$.}\label{fig:AI_Vs_CS_Pics_a}
\end{minipage}
\hfil
\begin{minipage}[b]{0.49\textwidth}
\centering
\includegraphics[width=\columnwidth]{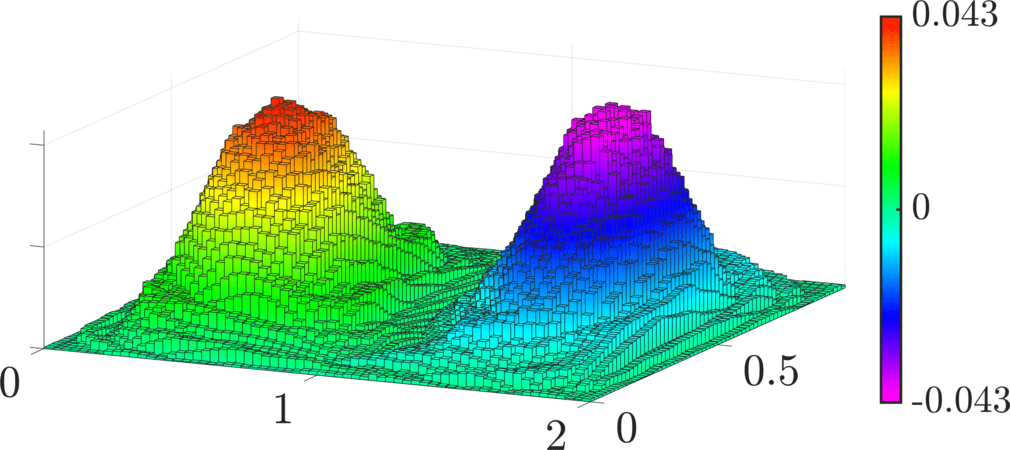}
\subcaption{Periodically driven double gyre at time $t=1$.}\label{fig:AI_Vs_CS_Pics_b}
\end{minipage}
\end{minipage}
\begin{minipage}[b]{\textwidth}
\centering
\begin{minipage}[b]{0.49\textwidth}
\centering
\includegraphics[width=\columnwidth]{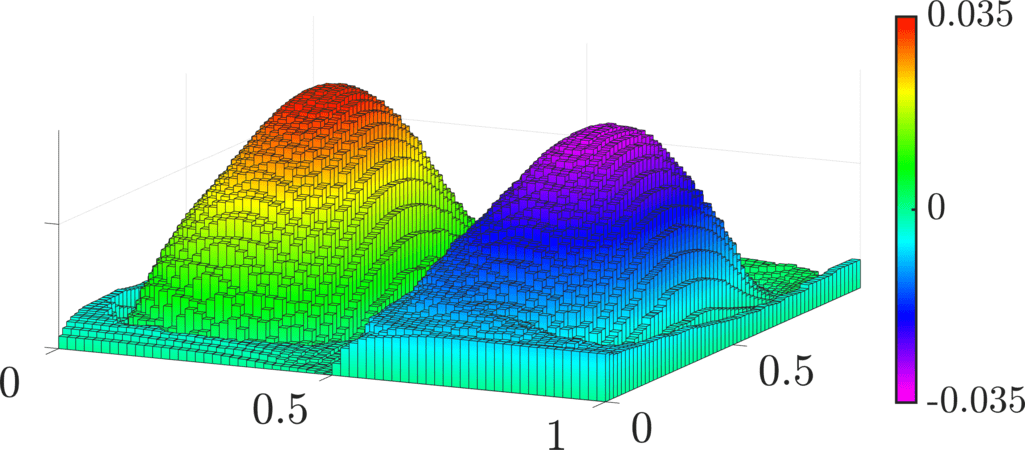}
\subcaption{Transitory double gyre at time $t=0$.}\label{fig:AI_Vs_CS_Pics_c}
\end{minipage}
\hfil
\begin{minipage}[b]{0.49\textwidth}
\centering
\includegraphics[width=\columnwidth]{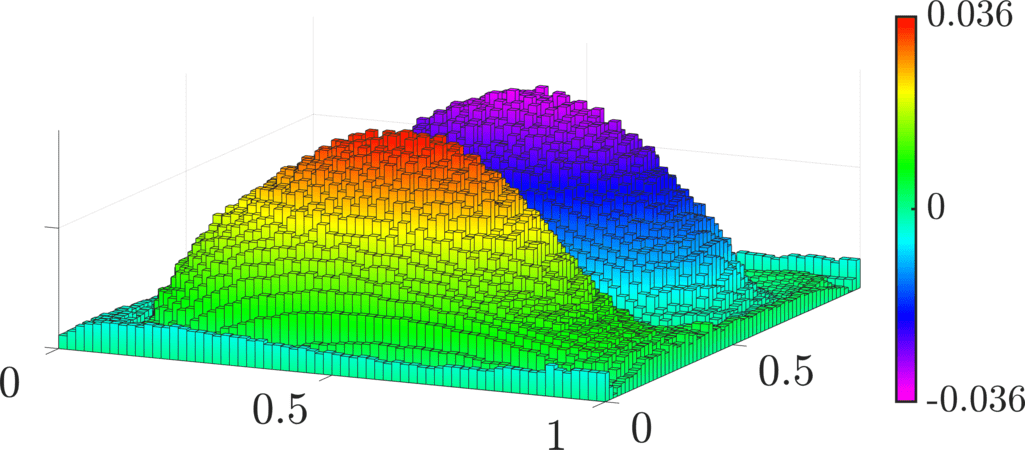}
\subcaption{Transitory double gyre at time $t=1$.}\label{fig:AI_Vs_CS_Pics_d}
\end{minipage}
\end{minipage}
\caption{Figures \ref{fig:AI_Vs_CS_Pics_a} and \ref{fig:AI_Vs_CS_Pics_b} show almost invariant structures as described by the (evolved) subdominant eigenvector of an Ulam matrix approximation to the transfer operator in the periodically driven double gyre flow, with parameters as in \cite{ShaddenShawnC2005Dapo}. 
Figures~\ref{fig:AI_Vs_CS_Pics_c} and \ref{fig:AI_Vs_CS_Pics_d} show
finite-time coherent structures as described by the (evolved) subdominant initial time singular vector of a composition of $10$ Ulam matrices describing the evolution of the transitory double gyre flow introduced by \cite{MosovskyBA2011TiTD}. See~\cite{FroyPG2014} for a thorough discussion of both models.}
\label{fig:AI_Vs_CS_Pics}
\end{figure}

The transfer operator point of view can be interpreted as tracking the evolution of an initial ensemble of trajectories, or a density, through time. These methods were first found to be useful in the identification of almost-invariant sets in the 1990s \cite{DellnitzMichael1999OtAo}.
Coherent structures are the time-dependent generalisation of almost-invariant sets. The latter do not move over time, as illustrated in the top row of
Figure~\ref{fig:AI_Vs_CS_Pics}; coherent structures, are shown in the bottom row. 
In applications, transfer operator based methods were initially investigated in the area of molecular dynamics \cite{DeuflhardDellnitzJungeetal.1999}, and later in the context of geophysical flows, starting with the work of Froyland, Padberg, England, and Treguier \cite{FroylandGary2007Doco}. This approach was later developed to identify and track time-varying structures, beginning with the works of Froyland, Lloyd and Quas on multiplicative ergodic theory in \cite{FroylandGary2010Csai,FLQ2}.
A survey of these techniques is provided in~\cite{GonzlezTokman2018MultiplicativeET}. 
In their most basic form, transfer operators, and also their adjoints, called composition or Koopman operators (see e.g.  \cite{BudisicMarko2012AK, WilliamsKevrekidisRowley2015}),
 provide a spectral approach for the study of autonomous dynamical systems. Indeed, transfer and Koopman operators are linear operators that encode the global behaviour of a dynamical system. Roughly speaking, their eigenvectors provide dynamically meaningful modes and the corresponding eigenvalues encode their rates of decay. This point of view has given rise to various methods that investigate transport and mixing in flows, see e.g. \cite{BanischKoltai},  \cite{KoltaiRenger},  \cite{FroylandRossSakeralliou} and references therein. Another, more geometric point of view, which is also used to handle truly time-dependent dynamics, is provided by the so-called Lagrangian Coherent Structures approach, in which 
key barriers to transport are sought \cite{HallerGeorge2015LCS,AllshouseMichaelR.2015Lbmf,HallerKarraschKogelbauer}. A review of the more commonly used methods in this direction is found in \cite{BalasuriyaSanjeeva2018GLcs}. There, the authors also present a general framework that seeks to better characterise the coherence of quantities that co-evolve with the vector field.

Multiplicative ergodic theory is concerned with existence and properties of spectral type decompositions for non-autonomous dynamical systems. That is, for systems whose evolution rules change over time. This was initially developed by Oseledets in the 1960's \cite{MR0240280} and was adapted and expanded to the semi-invertible setting in \cite{FroylandGary2010Csai,FLQ2,GTQuas1,GTQuas2}. This extension is crucial to the study of transfer operators of non-autonomous dynamical systems because it covers cases where the dynamics are not necessarily invertible. When the theory applies, it provides a (finite or countable) list of Lyapunov exponents, $0 = \lambda_1 > \lambda_2 > . . . $, for the transfer operator cocycle, which encodes the decay rates associated to the non-autonomous system over time\footnote{These Lyapunov exponents should not be confused with the Lyapunov exponents associated to trajectories in physical space, some of which may be positive in the context of chaotic systems.}.

Associated to each Lyapunov exponent $\lambda_i$, there is a finite-dimensional time-dependent space  $E_i(\omega)$. These are the so-called Oseledets spaces or modes which, in a hierarchical way, encode information about the system's coherent structures. 
Their time dependence may be associated, for example, with seasonal and random fluctuations in the system. It is also closely related to the underlying time dependence of the dynamics.

An idea going back to Raghunathan \cite{RaghunathanM.1979ApoO} is that Lyapunov exponents and Oseledets modes can be approximated using singular values and singular vectors arising from singular value decompositions (SVDs) that correspond to longer and longer evolution times.
At a general level, this approach is related to computational algorithms to approximate Oseledets modes and  so-called covariant Lyapunov vectors,  e.g. \cite{GinelliEtAl,FroylandGary2013CcLv,Noethen}.
In the context of transfer operators and matrix approximations thereof, SVDs were employed in \cite{FLS_2010, FSM_2010} to identify coherent structures and finite-time coherent sets. 
The result of \cite{RaghunathanM.1979ApoO} was extended in \cite{GTQuas2} to the infinite dimensional setting, showing that SVD type decompositions also provide approximations to Oseledets modes and Lyapunov exponents in the context of transfer operator cocycles.

This work develops and implements SVD based algorithms that build on multiplicative ergodic theory insights to extract detailed information about coherent structures in dynamical systems. 
In the numerical applications investigated here, infinite dimensional transfer operators are approximated by finite rank stochastic matrices. This is done using a popular Galerkin projection technique known as Ulam's method. The algorithms are used to rank and track structures whose location and boundaries shift over time, and to identify time windows where structural changes, such as merging and splitting events, occur. Furthermore, an equivariance test is introduced. This is used to assess the quality of pairings among structures as they are followed through time.

The algorithms are tested on three models. The first two models investigate a forced double well potential under periodic and quasi-periodic forcing, respectively, and for a range of time windows. 
The final case study relies on data from the European Centre for Medium-Range Weather Forecasts (ECMWF)
and investigates the splitting of the Southern (Antarctic) Polar Vortex in 2002. This splitting was directly related to the first observed major warming in the Southern Hemisphere and the division of the Antarctic ozone hole into two parts \cite{NewmanPaul2005TUSH, OrsoliniYvanJ.2005Aoso, Charltonetal_SPV}. In this latter model, we compare our numerical results to the method of normalised evolved singular vectors introduced in \cite{FSM_2010}. 
In all three cases, our methods give important information regarding the location and time windows within which coherent structures merge and/or separate in the presence of an underlying time-dependent and possibly chaotic flow.

After this paper was submitted, the related work \cite{ndour2020predicting} became available. In
\cite{ndour2020predicting} the authors develop a set-oriented
bifurcation analysis to better understand and identify spectral signals
associated with bifurcations of the almost invariant patterns
characterising an underlying autonomous dynamical system. 
\section{Background and framework}
\subsection{Non-autonomous dynamical systems}\label{ss:nonauto}
Non-autonomous dynamical systems are characterised by the fact that the system's evolution rule changes from one step to the next, depending on the environment.
To model these external influences, we consider a \textit{driving system} described by a
tuple $(\Omega,\mathcal{F},\mathbb{P},\sigma)$, where $(\Omega,\mathcal{F},\mathbb{P})$ is a probability space and $\sigma: \Omega \to \Omega$. 
Each $\omega\in \Omega$ corresponds to a possible state of the environment, and the map $\sigma$ dictates how the environment changes from one step in time to the next. For Section~\ref{ssec:MET}, we will also require $\sigma$ to be invertible, $\mathbb P$ preserving and ergodic. That is, $\mathbb P(\sigma^{-1}(E))=\mathbb P(E)$ for every $E\in \mathcal F$, and if $\sigma^{-1}(E)=E$ for some $E\in \mathcal F$, then either $\mathbb P(E)=0$ or $\mathbb P(E)=1$.

The evolution rule for the system is defined by a collection of maps $T_{\omega}: X \to X$, indexed by $\omega\in \Omega$, where $X$ is the configuration space.
The map $T_{\omega}$ can be thought of as describing the terminal location of particles $x \in X$, initialised in the environment $\omega$, 
after one step of the dynamics. 
The discrete time evolution of particle $x\in X$ in forward time ($n>0$) can thus be described by the following composition,
\begin{equation}
T^{(n)}_{\omega}(x) \coloneqq T_{\sigma^{n-1}\omega} \circ \cdots \circ T_{\sigma^{2}\omega} \circ T_{\sigma\omega} \circ T_{\omega}(x).
\end{equation} 
In what follows, we assume $X$ is a manifold and $T_\omega$ is non-singular for every $\omega \in \Omega$. That is, $\vol(T_{\omega}^{-1}(A))=0$ for all measurable $A \subset X$ such that $\vol(A)=0$, where $\vol$ denotes the Lebesgue measure on $X$.

\subsection{Transfer operators}
To each instance of the evolution rule $T_\omega$, there is an associated \textit{transfer} or \textit{Perron-Frobenius operator} $\mathcal{L}_{\omega}: L^{1}(X,\vol) \to L^{1}(X,\vol)$, where $ \mathcal{L}_{\omega} f$ is defined by the property that for every measurable $A\subset X$,
\begin{equation}\label{eq:PF}
\int_A \mathcal{L}_{\omega} f(x)\,d\vol(x) = \int_{ T_{\omega}^{-1}(A) } f(x)\,d\vol(x).
\end{equation}
If $f \in L^{1}(X,\vol)$ is the density of an ensemble of initial conditions in $X$, then $\mathcal{L}_{\omega}f$ describes the result of evolving this density under $T_{\omega}$. The Perron-Frobenius operator is Markovian in the sense that it is linear, $f \ge 0$ implies $ \mathcal{L}_{\omega} f \ge 0$  and $ \|\mathcal{L}_{\omega}f\|_{L^{1}} \le \|f\|_{L^{1}}$ \cite{LasotaAndrzej1994CFaN}.

As before, the $n$-step evolution of densities under the non-autonomous dynamics is described by the following composition, 
\begin{equation} \label{eq:PFTO}
\mathcal{L}_{\omega}^{(n)} \coloneqq \mathcal{L}_{\sigma^{n-1}\omega} \circ \cdots \circ \mathcal{L}_{\sigma^{2}\omega} \circ \mathcal{L}_{\sigma\omega} \circ \mathcal{L}_{\omega}.
\end{equation}

\subsection{Ulam's method and numerical approximations}\label{ss:Ulam}

In numerical investigations, a Galerkin projection known as \textit{Ulam's method} \cite{UlamStanislawM.1960Acom} is often employed to approximate the transfer operator. This method partitions $X$ into a pairwise disjoint collection of \textit{bins} $\{ B_{1}, \ldots, B_{m} \}$ of positive volume. For $i=1, \ldots ,m$, the indicator function on bin $B_{j}$ is denoted by $\mathbbm{1}_{B_{j}}$.
The Ulam approximation to $\mathcal{L}_{\omega}$  
is given by an $m\times m$ matrix $P(\omega)$ whose $ij$-th entry is obtained by computing the proportion of $Q$ uniformly distributed test points $x_{i,1},\ldots,x_{i,Q}\in B_{i}$ that move to $B_{j}$ after one step of the dynamics. That is,
\begin{equation}\label{eq:UlamsNum}
\boldsymbol{P}(\omega)_{i,j} = \frac{1}{Q} \sum\limits_{q=1}^{Q} \mathbbm{1}_{B_{j}}(T_{\omega} (x_{i,q})).
\end{equation}	
When the map $T_{\omega}$ arises from integration of a vector field,
one also approximates $T_{\omega}$ numerically, e.g. using Runge-Kutta numerical integration of the time dependent vector field which, if necessary, is interpolated linearly in space and time. 

To construct each Ulam matrix $\boldsymbol{P}({\omega})$ numerically as a sparse matrix, the software package GAIO \cite{Dellnitz2001} is used. These matrices are then combined to define a matrix cocycle approximating the $n$-step transfer operator $\mathcal L_\omega^{(n)}$ by
\begin{equation} \label{eq:TO_Ulam_illustration}
\boldsymbol{P}^{(n)}_{\omega} \coloneqq   
\boldsymbol{P}(\omega) \boldsymbol{P}(\sigma \omega) \cdots
\boldsymbol{P}(\sigma^{n-1} \omega)
,
\end{equation}
or equivalently,
\begin{equation} \label{eq:TO_Ulam_illustrationT}
\big(\boldsymbol{P}^{(n)}_{\omega} )^T\coloneqq \boldsymbol{P}(\sigma^{n-1} \omega) ^T\cdots \boldsymbol{P}(\sigma \omega)^T \boldsymbol{P}(\omega)^T, 
\end{equation}
where $P^T$ denotes the transpose of $P$.

\begin{figure}[H]
\includegraphics[width=\textwidth,center]{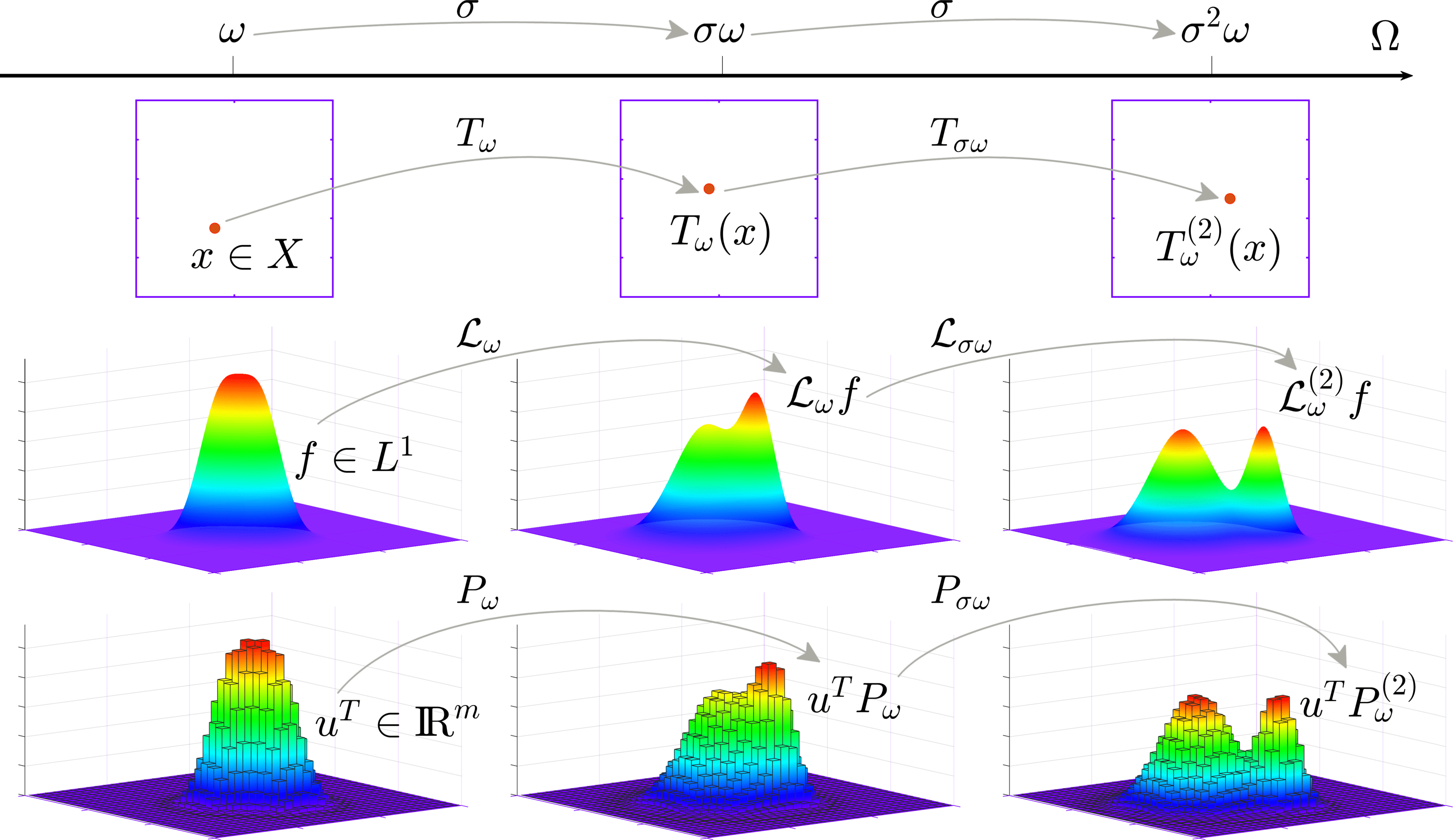}
\caption{Evolution in non-autonomous dynamical systems: driving system (above the arrow), particle evolution (2nd row), transfer operators (3rd row) and Ulam's method (bottom row).} \label{fig:TO_Ulam_Map}
\end{figure}

A visual summary of the concepts introduced in Sections~\ref{ss:nonauto}--\ref{ss:Ulam} is presented in Figure~\ref{fig:TO_Ulam_Map}. The bottom two rows present three-dimensional views of densities evolving under the dynamics.
For a more comprehensive perspective, two-dimensional visualisations will be employed in later figures. As here, the colour will reflect the value of the density at the corresponding location.
\vspace{-0.25cm}
In practice, one may also be interested to investigate systems where the domain and range are not necessarily the same space. 
For example, this is useful to track the evolution of particles initially seeded in some location of interest, as in \cite{FSM_2010}.
An extension of the above framework in this vein is possible by allowing $\omega$-dependent spaces $\{X_\omega\}_{\omega\in \Omega}$ and  maps $T_\omega: X_{\omega} \to X_{\sigma\omega}$. Further details regarding the application of Ulam's method to transfer operators can be found in \cite{DellnitzMichael2009OtAo, FP_2009, KlusStefan2015Otna, FroyPG2014,FroylandGary2013CcLv}.

\subsection{Singular value decompositions and multiplicative ergodic theory} \label{ssec:MET} 

The \textit{singular value decomposition} (SVD) algorithm decomposes a matrix $\boldsymbol{A} \in \mathbb{R}^{m \times m'}$ as $\boldsymbol{A}=\boldsymbol{U}\boldsymbol{S}\boldsymbol{V}^{T}$, where $\boldsymbol{U} \in \mathbb{R}^{m \times m}$ and $\boldsymbol{V} \in \mathbb{R}^{m' \times m'}$ are unitary matrices, whilst $\boldsymbol{S} \in \mathbb{R}^{m \times m'}$ is a diagonal matrix with $s_{1} \ge s_{2} \ge \ldots \ge s_{p} \ge 0$ on the diagonal and $p=\min(m,m')$. The entries in $S$, known as \textit{singular values} $\{s_{j}\}_{1\leq j \leq p}$, are uniquely determined by $\boldsymbol{A}$. The columns of $U$ are the corresponding \textit{left singular vectors} $\{u_{j}\}_{1\leq j \leq m}$ and those in $V$ are the \textit{right singular vectors} $\{v_{j}\}_{1\leq j \leq m'}$. These vectors satisfy the relation $u_j^T A =s_j v_j^T$, for $1\leq j \leq p$. Thus, vectors $\{u_{j}\}$ can be thought of as modes corresponding to the initial time which, under the application of $A$, evolve into multiples of the final time modes $\{v_{j}\}$.
For the case of square matrices where $m=m'$, if the $\{s_{j}\}$ are distinct then both $\{u_{j}\}$ and $\{v_{j}\}$ are uniquely determined up to a sign. 

The \textit{multiplicative ergodic theorem} (MET) provides a spectral type decomposition which allows one to investigate non-autonomous dynamical systems with spectral techniques. 
For example, in the context of matrix cocycles under right multiplication, Froyland, Lloyd and Quas show the following result, 
 allowing for both invertible and non-invertible matrices to be considered.

\begin{thm}\cite[Theorem 4.1]{FroylandGary2010Csai}
Let $\sigma$ be an ergodic, invertible measure-preserving transformation of the probability space $(\Omega,\mathcal{F}, \mathbb{P})$. Let $m\in \mathbb N$ and $P: \Omega \rightarrow \mathbb{R}^{m \times m}$  be a measurable family of matrices  satisfying $ \int \log^{+} {\|\boldsymbol{P}(\omega)\|} \, d\mathbb{P}(\omega) < \infty $ where $\log^{+}x\coloneqq \max\{ 0, \log x\}$.
Then there exist Lyapunov exponents $\lambda_{1} > \lambda_{2} > \ldots > \lambda_{\ell}\geq -\infty$,  numbers  $m_{1}, \ldots, m_{\ell}\in \mathbb N$ with $m_{1} + \ldots + m_{\ell} = m$ and measurable families of subspaces, called Oseledets spaces, $W_{j}(\omega) \subset \mathbb{R}^{m}$, $1\leq j \leq \ell$, such that for $\mathbb{P}$-almost every $\omega$ the following hold:

\begin{enumerate}
\item dim $W_{j}(\omega)=m_{j}$.
\item $\R^{m}=\bigoplus_{j=1}^{\ell} W_{j}(\omega)$.
\item $W_{j}(\omega)\boldsymbol{P}(\omega) = W_{j}(\sigma\omega)$ when $\lambda_j \neq -\infty$ and $W_{j}(\omega)\boldsymbol{P}(\omega) \subseteq W_{j}(\sigma\omega)$ in general.
\item $\lim_{n \to \infty} (1/n) \log \|u^T \boldsymbol{P}(\omega) \boldsymbol{P}(\sigma \omega) \cdots
\boldsymbol{P}(\sigma^{n-1} \omega)\|=\lambda_{j}$ for all $u^T \in W_{j}(\omega) \setminus \{ 0 \}$.
\end{enumerate}
\end{thm}
In short, the Oseledets spaces $W_j(\omega)$ provide a decomposition (splitting) of $\R^m$ into  $\omega$-dependent (time-varying) modes ordered by decay rate. The rate is determined by the Lyapunov exponent $\lambda_{j}$. These modes evolve according to point 3, which is referred to as the equivariance property.

The previous result involves limits as $n$ (time) approaches infinity.
The focus of this work is on extracting dynamical information from leading modes arising from singular value decompositions of matrix products  of the form $P_{\omega}^{(n)}=\boldsymbol{P}(\omega) \boldsymbol{P}(\sigma \omega) \cdots
\boldsymbol{P}(\sigma^{n-1} \omega)$, coming from Ulam matrix cocyles as in~\eqref{eq:TO_Ulam_illustration}.
The observation that these are related to objects arising from the MET goes back to Raghunathan's work \cite{RaghunathanM.1979ApoO}.
The Lyapunov exponents are approximated by the exponential growth rates of the associated singular values, that is $ \lambda_{j}=\lim_{n \to \infty} \frac{1}{n} \log s_{j}(\boldsymbol{P}^{(n)}_{\omega})$. This hints at a numerical means for calculating decay rates.
Furthermore, it follows from the approach of \cite{GTQuas2} that the left\footnote{Notice that matrices are multiplied on the right.} singular vectors $u_{j}^{T}$ of  $P_{\omega}^{(n)}$ approach  the Oseledets \textit{filtration} space $F_{k}(\omega) = \bigoplus_{1\leq i \leq k}W_i(\omega)$ as $n\to\infty$, where $k=1$ if $j\leq m_1$ and otherwise $1< k \leq \ell$ is such that $m_1+\dots + m_{k-1}< j \leq m_1+\dots + m_{k}$. 
\section{Algorithms} \label{ssec:NovelTechs}
The content of multiplicative ergodic theorems may be summarised by saying that leading Lyapunov exponents and Oseledets spaces associated to a matrix cocycle approximate the principal features of the underlying non-autonomous system, in the sense that they are the most persistent over time.
The algorithms introduced in this section aim at  investigating and tracking such features through time, by exploring two different approaches to matching modes across different time windows.

Algorithm~\ref{alg:RWinds} describes the creation of matrix products corresponding to \textit{rolling windows} following each other through time, and corresponding to blocks of equal length coming from an underlying matrix cocycle. Each of these matrices is decomposed into singular values and vectors. 
When each Lyapunov exponent in the MET has multiplicity one, one would expect to be able to track each mode separately over time, provided the rolling window length is large enough. 
However, in many interesting cases, and in most of the applications considered in this work, this separation is not achieved with reasonable window lengths. Thus, there may be intersections along the paths traced by singular values corresponding to rolling windows of equal length starting at different times.

Algorithm~\ref{alg:SValsDijk} presents a method to match or pair the singular values corresponding to different rolling windows as they evolve. 
 Algorithm~\ref{alg:L2VecDist} presents another such method, but uses only information from the corresponding singular vectors at neighbouring times to set the pairings. To illustrate how the results of Algorithms~\ref{alg:SValsDijk} and \ref{alg:L2VecDist} are related to the dynamics, Algorithm~\ref{alg:Movies} describes a means to visualise the evolution of structures in non-autonomous systems. Algorithm~\ref{alg:Eqvr} investigates the efficacy of the pairings from Algorithms~\ref{alg:SValsDijk} and~\ref{alg:L2VecDist} in approximating equivariant subspaces. 
\subsection{Rolling windows and singular value decompositions}\label{ssec:RWSVs}
The term \textit{rolling windows} refers to a collection of time intervals associated to specific initial time $t_{i}$, final time $t_{F}$, single step flow time $\tau$ and window length $n$. 
The initial and final times correspond to a given time interval $[t_i,t_F]$, where there is available data describing the evolution of a time-dependent vector field. The flow time $\tau$ is used to construct row stochastic Ulam matrices associated with transitions from $t$ to $t+\tau$, and denoted by $\boldsymbol{P}_{t}$, for $t=t_i, t_i+\tau, t_i+2\tau\dots$. 

Each time window $W_{t_{0}}^{(n)}$ covers $n$ time steps from time $t_{0}$ to time $t_{f}=t_{0}+n\tau$. The windows begin at a specified number of neighbouring times. For example, the window $W_{t_{i}}^{(n)}$ precedes $W_{t_{i}+\tau}^{(n)}$, which precedes $W_{t_{i}+2\tau}^{(n)}$,  and so on while $t_{f}+n\tau \le t_{F}$. For ease of explanation we set $\tau=1$ but altering $\tau$ for various flow times is straightforward to account for. For each time window, Algorithm~\ref{alg:RWinds} computes the SVD associated to the $\mathcal N$ largest singular values.
\begin{algorithm}[H]
\caption{Forward time rolling windows and singular value decompositions}
\label{alg:RWinds}
\begin{algorithmic}[1]
\State Set initial and final available times $t_{i}$ and $t_{F}$ 
\State Set number of singular values $\mathcal N$  
\State Choose time window length $n$
\For{$k \gets t_{i}$ to $t_{F}-n$}
\State $\boldsymbol{P}^{(n)}_{k} \gets \boldsymbol{P}_{k} \cdot \boldsymbol{P}_{k+1} \cdots \boldsymbol{P}_{k+n-1}$ 
\State $[U_{k}^{(n)}, \, S_{k}^{(n)}, \,V_{k}^{(n)}] \gets \textrm{SVD}(\boldsymbol{P}^{(n)}_{k}, \; \mathcal N)$
\EndFor
\Return $U_{K}^{(n)}, \, S_{K}^{(n)}, \,V_{K}^{(n)}$ where $K=\{t_{i},\ldots,t_{F}-n\}$
\end{algorithmic}
\end{algorithm}
\vspace{-0.25cm}
In this algorithm, 
$S_{k}^{(n)}$ gives the $\mathcal N$ largest singular values of $\boldsymbol{P}^{(n)}_{k}$, while
$U_{k}^{(n)}$ and $V_{k}^{(n)}$ are the associated collections of left and right singular vectors, respectively. 
\vspace{-0.25cm}
\subsection{Following coherent structures through time} \label{ssec:FCST}
In order to identify time windows associated with distinctive behaviour in the underlying dynamics, such as structural changes, we develop pairing techniques that attempt to track the different modes through time. This pairing process is difficult because the structures in the dynamical system may shift in coherence as they and their boundaries evolve and even interact. When a structure shifts in dynamical dominance it becomes comparatively more (less) coherent and is associated with a singular value that is ranked higher (lower) than the ranking of the singular value with which it is initially associated.
Algorithms~\ref{alg:SValsDijk} and~\ref{alg:L2VecDist} utilise either the path of singular values or singular vectors through time to track the evolution of dominant modes. Algorithm~\ref{alg:SValsDijk} looks to identify the path of various singular values by minimising the total change in singular values over neighbouring windows. Algorithm~\ref{alg:L2VecDist} focuses on pairing singular vectors from $U_{k}^{(n)}$ to their best match within those of $U_{k+1}^{(n)}$, accounting for the time evolution.
\subsubsection{Finding paths with singular values} \label{sssec:FMTT_SVal} \mbox{}
Algorithm~\ref{alg:SValsDijk} tracks the paths of the leading $\mathcal N$ singular values of sequential rolling windows obtained using Algorithm~\ref{alg:RWinds}. 
To pair singular values, Algorithm~\ref{alg:SValsDijk} begins with the construction of a directed, weighted graph $\boldsymbol{G}=(\mathcal{S}, \mathcal{E})$. The 
elements of $\mathcal{S}$ are given by $\{S_{k,j}^{(n)}\}_{t_i\leq k \leq t_F-n, 1\leq j \leq \mathcal N}$,
the collection of the $\mathcal{N}$ largest singular values of $\{\boldsymbol{P}^{(n)}_{k}\}_{t_i\leq k \leq t_F-n}$. Forward time linkages provide for ordered pairs of nodes, $\mathcal{E}$. The elements of $\mathcal{E}$ at time $k$ join all the neighbouring nodes $S_{k,j_{i}}^{(n)}$ and $S_{k+1,j_{f}}^{(n)}$ for  $j_{i},j_{f} \in \{1,\ldots,\mathcal{N}\}$. The weight on an edge between sequential nodes is the
distance between the corresponding points in a plot of $k$ vs $S_k^{(n)}$.

Dijkstra's algorithm is applied to $\boldsymbol{G}$ to find a path of minimum cost in terms of the distance between pairs of singular values over time. Once a minimal path is found, it is recorded and the associated nodes and edges are removed from $\boldsymbol{G}$. The graph is then re-consolidated by redefining a new $\mathcal{E'}\subset\mathcal{E}$ obtained by deleting all the edges contributing to the (removed) minimal path. One then iterates the method to find the next path of least cost. This method continues until no paths remain in $\boldsymbol{G}$.
\vspace{-0.25cm}
\begin{algorithm}[H]
\caption{Tracking modes through time with singular values}
\label{alg:SValsDijk}
\begin{algorithmic}[1]
\State Set $t_{i}$, $t_{F}$, $n$ and $\mathcal N$, and define $S_{K}^{(n)}$ as per Algorithm~\ref{alg:RWinds}
\State Create directed graph $\boldsymbol{G} \gets (\mathcal{S}, \mathcal{E})$ as described in the preamble~\ref{sssec:FMTT_SVal}
\For{$k \gets t_{i}$ to $t_{F}-n-1$}
\State Calculate edge weights $M(j_{i},j_{f},k)$ between each initial node $j_{i}$ at time $k$ \par and each final node $j_{f}$ at time $k+1$ as \par
$M(j_{i},j_{f},k) \gets \sqrt{(S_{k,j_{i}}^{(n)}-S_{k+1,j_{f}}^{(n)})^{2}+1}$ 
\EndFor
\State $\vars{path} \gets 1$
\While {$\mathcal{E} != \emptyset$ } 
\State Apply Dijkstra's algorithm on $\boldsymbol{G}$ to identify path of minimum total cost
\State Define  $\hat{\mathcal{S}}^{(\vars{path})}$, $\hat{\mathcal{E}}^{(\vars{path})}$ by the graph that defines the path of least cost
\State Remove relevant edges $\mathcal{E}' \gets \mathcal{E} \setminus \hat{\mathcal{E}}^{(\vars{path})}$ and vertices $\mathcal{S}' \gets \mathcal{S} \setminus \hat{\mathcal{S}}^{(\vars{path})}$
\State Re-consolidate $\boldsymbol{G}\gets(\mathcal{S}', \mathcal{E}')$
\State $\vars{path} \gets \vars{path} + 1$
\EndWhile
\State Denote the average value along each of the $\mathcal{N}$ paths by $\bar{S}^{(\{1,\ldots,\mathcal{N}\})}$
\State Reorder paths by average value \par $[ \, \sim \, , \; \{\vars{sorted\_mode\_order}\}] \gets \textrm{sort}(\bar{S}, \; \textit{`descending'})$
\State Define the path of modes tracked by singular values as \par $S^{(\{1,\ldots,\mathcal{N}\})}_{S} \gets \hat{S}^{(\{\vars{sorted\_mode\_order}\})}$
\end{algorithmic}
\end{algorithm}
\vspace{-0.5cm}
Paths tracked by $S_{S}$ identify the movement of modes through time and can indicate the occurrence of distinctive dynamical behaviour. Paths of interest may be signalled by a quick succession of crossings of singular value paths, which is associated with a switch in the comparative dominance of the associated structures, or by qualitative changes such as peaks in singular values, indicating a transition between phases of increasing and decreasing coherence.

\subsubsection{Finding paths with singular vectors} \label{sssec:FMTT_SVec} \mbox{}
Algorithm~\ref{alg:L2VecDist} tracks sorted, paired paths of singular values for neighbouring time windows by minimising the distance between singular vectors in a relevant metric. It concentrates on minimising the Euclidean norm of the difference between two neighbouring singular vectors from $U_{k}^{(n)}$ and $U_{k+1}^{(n)}$, taking into account the one step evolution from $P_k$. Algorithm~\ref{alg:L2VecDist} iteratively minimises $U_{k,j_{i}}^{(n)}$ with respect to some neighbouring vector $U_{k+1,j_{f}}^{(n)}$. The pairing vector $U_{k+1,j_{f}}^{(n)}$ is chosen to minimise the distance over all (remaining) choices of pairs as both $j_{i}$ and $j_{f}$ are unique at each time step. The path of a mode associated with some $U_{k,j_{i}}^{(n)}$ at time $k$ is thus permitted to shift in dominance to another $U_{k+1,j_{f}}^{(n)}$ at the neighbouring time $k+1$. 

\begin{algorithm}[H]
\caption{Tracking modes through time using singular vectors}
\label{alg:L2VecDist}
\begin{algorithmic}[1]
\State Set $t_{i}$, $t_{F}$, $n$ and $\mathcal N$, and define
$S_{K}^{(n)}$ and $U_{K}^{(n)}$ as per Algorithm~\ref{alg:RWinds}
\State Create the initial mode association $\hat{S}^{(j)}_{t_{i}} \gets S_{t_{i},j}^{(n)}$ for $j\in\{1,\ldots,\mathcal{N}\}$
\For {$k \gets t_{i}$ to $t_{F}-n-1$}
\State Define initial sets characterising all possible transitions $j',j'' \gets \{1, \ldots, \mathcal{N} \}$
\While {$j' != \emptyset$ } 
\State \vars{min\_dist} $ \gets \min_{j',j''}{ \left( \; { \left\| {U_{k,j'}^{(n)T}\boldsymbol{P}_{k}/}{{ \| U_{k,j'}^{(n)T}\boldsymbol{P}_{k} \| }_{2}}\pm U_{k+1,j''}^{(n)T} \; \right\| }_{2}  \; \right) } $ \label{lst:line:oneStep}
\State \vars{min\_modes} $ \gets \arg\min_{j',j''}\left(\vars{min\_dist}\right)$
\State Set $j' \gets j' \setminus \vars{min\_modes}(1)$ and $j'' \gets j'' \setminus \vars{min\_modes}(2)$
\State Create the new mode association $\hat{S}^{(\vars{min\_modes}(1))}_{k+1} \gets S_{k+1,\vars{min\_modes}(2)}^{(n)}$
\EndWhile
\EndFor
\State Characterise each path by average value as \par $\bar{S}^{(j)} \gets \frac{1}{t_{F}-n-1-t_{i}} \sum_{k=t_{i}}^{t_{F}-n-1} \hat{S}^{(j)}_{k}$
\State Re-sort in descending order as \par $[ \, \sim \, , \; \{\vars{{sorted\_modes}}\}] \gets \textrm{sort}(\bar{S}, \; \textit{`descending'})$

\State Define the path of modes tracked by left singular vectors as \par ${S}_{U}^{(\{1,\ldots,\mathcal{N}\})} \gets \hat{S}^{(\text{\{{sorted\_modes}\}})}$ \label{lst:line:oneStepS}
\end{algorithmic}
\end{algorithm}

\subsection{Visualising the evolution of coherent structures}
Algorithm~\ref{alg:Movies} describes a way to visualise the evolution of structures associated with relevant singular vectors. Recall from Section~\ref{ssec:RWSVs} that the time window $W_{k}^{(n)}$ describes the evolution of the dynamical system over the time period $k$ to $k+n$. In Algorithm~\ref{alg:Movies}, we utilise this fact and set the $\tilde{n}$-th frame in the animation to be a plot of the singular vector evolved for $\tilde{n}$ steps where $0 \le \tilde{n} \le n$. 

\begin{algorithm}[H]
\caption{Animating tracked modes over time}
\label{alg:Movies}
\begin{algorithmic}[1]
\State Set $t_{i}$, $t_{F}$, $n$ and $\mathcal N$  as per Algorithm~\ref{alg:RWinds}
\State Set $1\leq j \leq \mathcal N$,  the mode to be displayed
\State Define $S_{K}^{(n)}$ and $U_{K}^{(n)}$ as per Algorithm~\ref{alg:RWinds}
\State Set $\tilde{S}$ using either \par ${S}_{S}$ of Algorithm~\ref{alg:SValsDijk} or \par
 ${S}_{U}$ from Algorithm~\ref{alg:L2VecDist}
\State Let $u^{(n)}_{k,j}$ be the column of $\tilde{U}_{k}$ associated with $\tilde{S}^{(j)}$ for a chosen $k \in K$
\For{$\tilde{n} \gets 0$ to $n$}
\State Define $\boldsymbol{P}^{(\tilde{n})}_{k}$ as per Algorithm~\ref{alg:RWinds} noting that $\boldsymbol{P}^{(0)}\coloneqq \boldsymbol{Id}$ 
\State $\hat{u}_{k,j}^{(\tilde{n},n)} \gets u^{(n)T}_{k,j}\boldsymbol{P}^{(\tilde{n})}_{k} \, / \, {\| u^{(n)T}_{k,j}\boldsymbol{P}^{(\tilde{n})}_{k} \|}_{2}$
\EndFor
\end{algorithmic}
\end{algorithm}
\vspace{-0.25cm}
To best visualise the dynamics over time, the colour scale limits of each frame in an animated time window are defined as $\pm\max(|u^{(n)}_{k,j}|)$. One other way to visualise the evolution of modes consists of fixing a window length $n$ and simply plotting the subsequent realisations of singular vectors over a time frame contained in  $[t_i, t_F-n]$. Any plots or animations that utilise Algorithm~\ref{alg:Movies} will mention this explicitly. If no algorithm is mentioned, then the plots or animations depict singular vectors of neighbouring windows.

\subsection{Assessing the equivariance of evolved modes}
The primary purpose of Algorithm~\ref{alg:Eqvr} is to assess how effective Algorithms~\ref{alg:SValsDijk} and \ref{alg:L2VecDist} are at pairing modes across rolling time windows.
In particular,  the final vectors $\tilde{V}_{k}$  of rolling time windows are compared with the initial vectors $\tilde{U}_{k+n}$ of adjacent time windows. 
For this, the algorithm relies on a measure of equivariance mismatch between vectors, denoted by $\varsigma \in [0,1]$, which is related to point 3 (equivariance) of the multiplicative ergodic theorem
of Section~\ref{ssec:MET}. When pairing is effective, one expects low values of equivariance mismatch. If pairing is not effective, one expects a maximum value of $1$, which is realised when the vectors are orthogonal.

\begin{algorithm}[H]
\caption{Equivariance mismatch in forward time}
\label{alg:Eqvr}
\begin{algorithmic}[1]
\State Set $t_{i}$, $t_{F}$, $n$ and $\mathcal{N}$ as per Algorithm~\ref{alg:RWinds}.
\State Define $U_{K}^{(n)}$, $S_{K}^{(n)}$ and $V_{K}^{(n)}$ as per Algorithm~\ref{alg:RWinds}.
\State Set $\tilde{S}$ using either \par ${S}_{S}$ of Algorithm~\ref{alg:SValsDijk} and thus set $T=S$, or \par ${S}_{U}$
from Algorithm~\ref{alg:L2VecDist} thus setting $T=U$
\For{$j \gets 1$ to $\mathcal{N}$}
\For{$k \gets t_{i}$ to $t_{F}-n$}
\State $\tilde{V}^{(j)}_{k}$ and $\tilde{U}^{(j)}_{k}$ are singular vectors corresponding to $\tilde{S}^{(j)}_{k}$
\State $\varsigma^{(j)}_{T,k} \gets { \min\left( \; {\|\tilde{V}^{(j)}_{k} + \tilde{U}^{(j)}_{k+n} \|}_{2}, {\|\tilde{V}^{(j)}_{k} - \tilde{U}^{(j)}_{k+n}\|}_{2} \; \right)}/{\sqrt{2}}$
\EndFor
\EndFor
\end{algorithmic}
\end{algorithm}
\vspace{-0.25cm}
It is worth pointing out that all the applications in this work are concerned with (normalised) singular vectors of dimension at least $2^{12}$, so low values of equivariance mismatch are unlikely to be observed by chance. Indeed, Monte Carlo experiments show that in this context, if two such unit vectors are chosen at random, the expected equivariance mismatch value is higher than $0.99$, with minimum observed values consistently over $0.95$. For this reason, even apparently high equivariance mismatch values, such as $0.7$ or even $0.9$, can still be a good indication of coherence.

\section{Models and results} \label{sec:DataModel}
In order to test these algorithms we employ three models. The first two models describe the evolution of a double well potential subjected to small, time dependent perturbations to the vector field, either periodic or quasi-periodic. These changes allow for the merging and separation of two wells over time. The final model employs vector field reanalysis data from the European Centre for Medium-Range Weather Forecasts (ECMWF) to examine a splitting of the Southern (Antarctic) Polar Vortex. For the Ulam scheme, all models seed $100$ test points per bin and integrate the vector field using standard Runge-Kutta, interpolating linearly in space and time.

\subsection{Periodically forced double well potential} \label{ssec:DWPmodel} 
Our exploration begins with an analysis of highly idealised mergers and separations. We model a double well potential where the centre of each well is shifted slowly over time. In this model we consider two wells of equal depth and separation from the origin over all time. This allows us to model the evolution of two coherent structures in phase space. These two structures will merge or separate as the centres of the wells get close or far from each other, respectively.

The first step in this analysis is to implement Ulam's method. We define a grid of $2^{12}$ bins of equal volume over $X=[-\pi,\pi]^2$. To gain a clearer understanding of the dynamics of the mergers and separations we consider this model over 5 periods, each of length 100. That is, we set $t_{i}=0$ and $t_{F}=500$ for $\tau=1$. The two wells in this model are initially centred at $x=\pm 2$. Due to dissipation we expect a number of trajectories to exit the phase space over a given time period. 

The general model is defined by the following system of differential equations
\begin{equation} \label{eqn:SimpleModel}
  \begin{aligned}
    \dot{x}(t) & = y(t) \\
    \dot{y}(t) & = x(t)\left(\frac{x(t)}{2} + \alpha(t)\right)\left(\alpha(t) - \frac{x(t)}{2}\right).
  \end{aligned}
\end{equation}
If $\alpha$ were constant we would be in the Hamiltonian regime of energy conservation. To introduce the merging and separation of wells over time, we consider time dependent changes to the vector field. This is done by varying $\alpha$ between $0$ and $1$ as follows,
\begin{equation} \label{eqn:SimpleModel_Forcing}
\alpha(t) =
\left\{
  \begin{aligned}
    & 1 & \text{if} \quad & 0 \le t \le 10 \\
    & cos^{2}\left((t-10)\frac{\pi}{60}\right) & \text{if} \quad & 10 \le t \le 40 \\
    & 0 & \text{if} \quad & 40 \le t \le 60 \\
    & cos^{2}\left((t-30)\frac{\pi}{60}\right) & \text{if} \quad & 60 \le t \le 90 \\
    & 1 & \text{if} \quad & 90 \le t \le 100
  \end{aligned}
\right.
\end{equation}
which is extended periodically in $t\pmod{100}$.

This hybrid function ensures a clear merging and separation of structures under the non-autonomous dynamics. When $\tau=1$ the forcing is repeated every $100$ steps. The first merger of our two structures will occur no later than time $40$. Each merger is followed by a separation of the two structures. That separation takes place over the following $30$ time steps. This pattern of forcing is repeated periodically.

\begin{figure}[H]
\centering
\begin{minipage}[b]{\textwidth}
\centering
\begin{minipage}[b]{0.32\textwidth}
\centering
\includegraphics[width=\columnwidth]{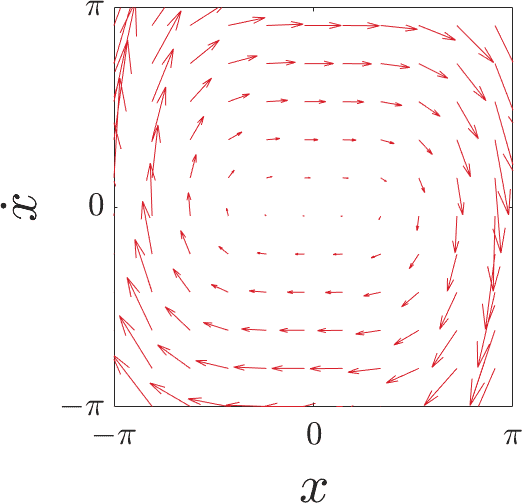}
\subcaption{$t=60$}\label{fig:DWP_t60}
\end{minipage}
\hfil
\begin{minipage}[b]{0.32\textwidth}
\centering
\includegraphics[width=\columnwidth]{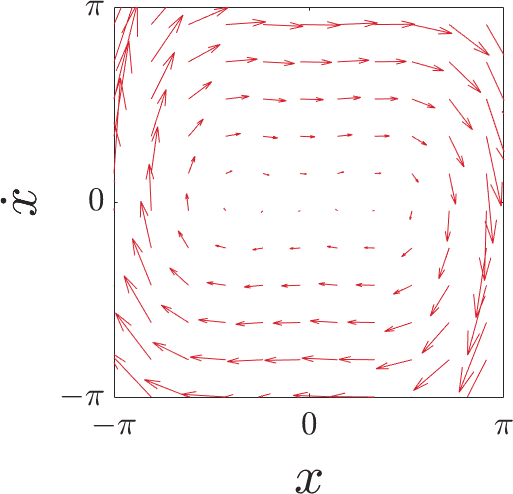}
\subcaption{$t=75$}\label{fig:DWP_t75}
\end{minipage}
\hfil
\begin{minipage}[b]{0.32\textwidth}
\centering
\includegraphics[width=\columnwidth]{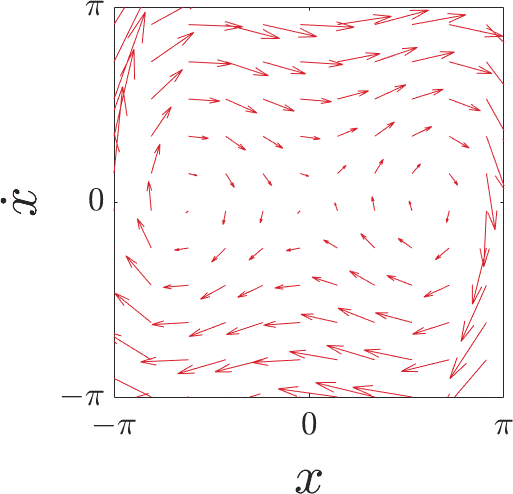}
\subcaption{$t=90$}\label{fig:DWP_t90}
\end{minipage}
\end{minipage}
\caption{Selected vector field instances for the periodically forced double well potential.}
\label{fig:DWPPics}
\end{figure}

Figures~\ref{fig:DWP_t60} and~\ref{fig:time_DWPV2} illustrate that the two wells have merged when $t$ is $60$, $160$ and so on. On the other hand, the two structures are separated around times when $t$ is close to $90$, $190$ and so forth, as shown in Figure~\ref{fig:DWP_t90}. Whilst it is not immediately clear from~\ref{fig:DWP_t75}, it is reasonable to assume the structures are experiencing some degree of separation by time $75$.
\vspace{-0.25cm}
In order to explore the behaviour of structures, we consider the forward time rolling windows of Algorithm~\ref{alg:RWinds}. In general, we consider time windows of length $n=50$ and $n=100$. The former allows time for one full merger (separation) to occur. The latter permits a merger and a separation, with both cores returning near to their initial positions. This is clear from the periodicity pattern in Figure~\ref{fig:time_DWPV2}. 

\subsubsection{Rolling time windows of length $n=50$} \label{sssec:50}
\begin{figure}[H]
\centering
{\includegraphics[width=\textwidth,center]{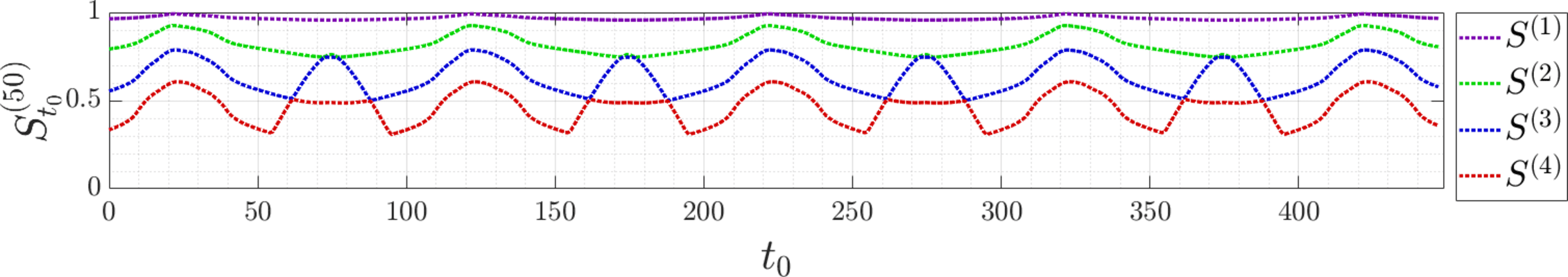}
\subcaption{Leading $4$ singular values of rolling windows for $n=50$ using Algorithm~\ref{alg:RWinds}.} 
\label{fig:delta50_full_4of4_a}
\includegraphics[width=\textwidth,center]{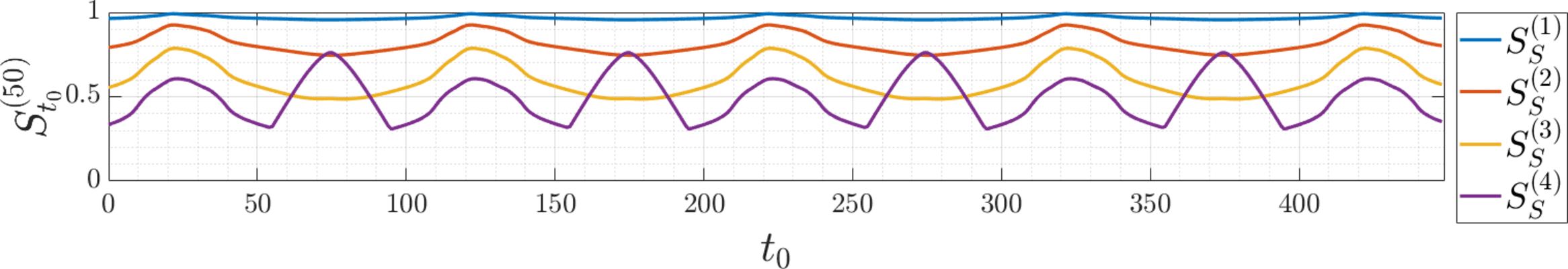}
\subcaption{Leading $4$ of $\mathcal{N}=4$ paths tracked by $S_{S}$ for $n=50$ using Algorithm~\ref{alg:SValsDijk}.} 
\label{fig:delta50_full_4of4_b}}
\caption{Tracking modes over rolling windows for the periodically forced double well potential.}
\label{fig:delta50_full_4of4}
\end{figure}

\paragraph*{\textit{Tracking structures over time}.} \label{sssec:TSOT} \vspace{-0.25cm}
We set $n=50$ and apply Algorithm~\ref{alg:RWinds} with $t_{i}$ and $t_{F}$ as in~\ref{ssec:DWPmodel}. This exploration is initially limited to the top $4$ modes and so we begin with $\mathcal{N}=4$. Figure~\ref{fig:delta50_full_4of4_a} illustrates the unpaired, unsorted singular values corresponding to rolling windows starting at times $t_{0}=t_i, t_{i}+1, \dots, t_{F}-n$ for these parameters. We apply Algorithm~\ref{alg:SValsDijk} to track modes through time. These results are presented in Figure~\ref{fig:delta50_full_4of4_b}. 
\vspace{-0.25cm}
Applying Algorithm~\ref{alg:L2VecDist} led to the same results as for Algorithm~\ref{alg:SValsDijk}. This occurs despite both algorithms relying on completely different methods to pair singular values through time. Algorithm~\ref{alg:SValsDijk} relies solely on the singular value structure whilst Algorithm~\ref{alg:L2VecDist} aims to match the associated singular vectors.

Figure~\ref{fig:delta50_full_4of4_b} shows that these techniques clearly identify three separated paths (blue, red and yellow in the electronic version). They also identify a fourth path which is characterised by much larger variations in singular value. The path of $S_{S}^{(4)}$ (purple in the electronic version) is associated with the fourth most dominant singular value at $t_0=55$. However, by  $t_0=75$ the singular value associated to this mode has risen in dominance to the second highest ranking. 
The timeframe of this rise and fall in the path of $S_{S}^{(4)}$ overlaps with the $4$ highest ranking crossings of singular values, which occur by times $t_0=62$, $73$, $78$ and $89$. As will be explained in the next section, it is the switching between increasing and decreasing phases of singular values, more so than the changes in dominance (crossings), which turns out to be related to the occurrence of fundamental changes in the structures associated with such modes.

\begin{figure}[H]
\centering
\begin{minipage}[b]{\textwidth}
\centering
\begin{minipage}[b]{0.49\textwidth}
\centering
\includegraphics[width=\textwidth,center]{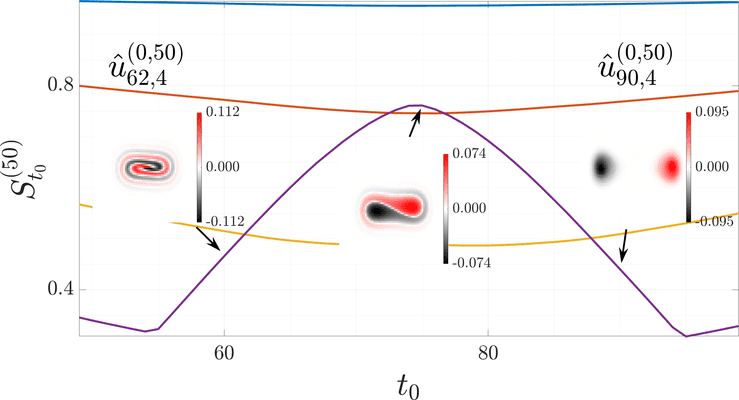}
\subcaption{Initial time singular vectors for indicated windows.} 
\label{fig:delta50_full_4of4_3plots_a}
\end{minipage}
\begin{minipage}[b]{0.49\textwidth}
\centering
\includegraphics[width=\textwidth,center]{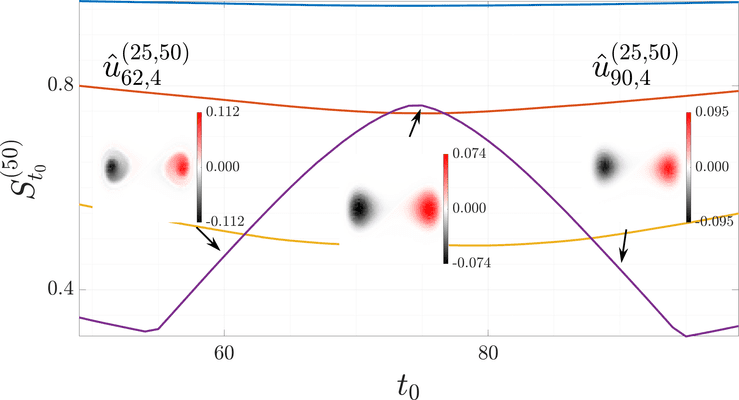}
\subcaption{Left singular vectors evolved for $25$ time steps.}
\label{fig:delta50_full_4of4_3plots_b}
\end{minipage}
\hfil \vfil
\begin{minipage}[b]{0.49\textwidth}
\centering
\includegraphics[width=\textwidth,center]{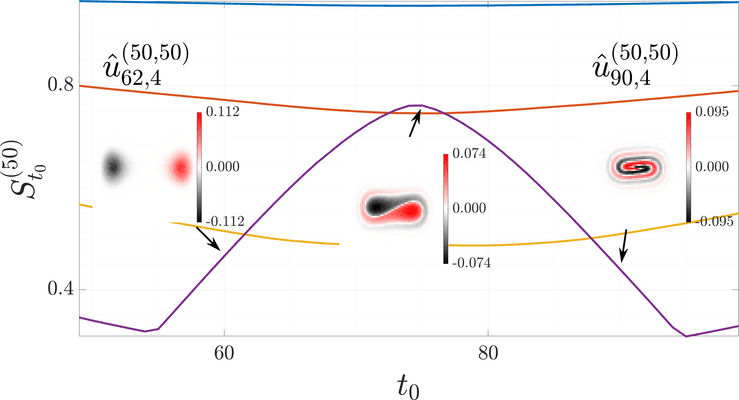}
\subcaption{Left singular vectors evolved for $50$ time steps.}
\label{fig:delta50_full_4of4_3plots_c}
\end{minipage}
\hfil
\begin{minipage}[b]{0.49\textwidth}
\centering
\includegraphics[width=0.8\textwidth,center]{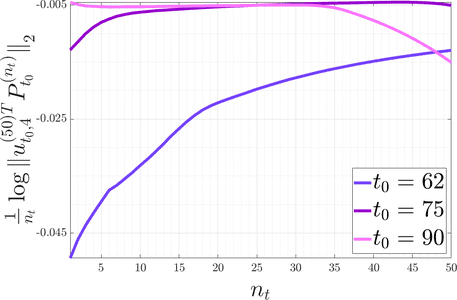}
\subcaption{$\frac{1}{n_{t}} \log{{\|u^{(50)T}_{t_0,4}\boldsymbol{P}_{t_0}^{(n_{t})}\|}_{2}}$ where $1\le n_{t}\le 50$ for $n_{t} \in \mathbb{Z}$.}
\label{fig:delta50_logPlot}
\end{minipage}
\caption{Tracking modes for time windows of length $n=50$, evolved using Algorithm~\ref{alg:Movies}.} 
\label{fig:delta50_full_4of4_3plots}
\end{minipage}
\end{figure}

Algorithm~\ref{alg:Movies} is now employed to explore the evolution of $S_{S}^{(4)}$ for time windows initialised at $t_0=62,75$ and $90$. Our results are presented in Figure~\ref{fig:delta50_full_4of4_3plots}. The initial time for each window is indicated by a black arrow. Figures~\ref{fig:delta50_full_4of4_3plots_a},~\ref{fig:delta50_full_4of4_3plots_b} and~\ref{fig:delta50_full_4of4_3plots_c} illustrate the singular vectors at three stages of evolution: initial time, mid-evolution and final time. Figure~\ref{fig:delta50_logPlot} indicates the coherency of structures in each of the associated windows, as the dynamics evolve. Negative values closer to zero suggest a greater degree of coherency for the associated structures as less mass is lost over time.
\vspace{-0.25cm}
Let us first consider the time window $W^{(50)}_{62}$. 
The upward trajectory of $S^{(4)}_S$ in this region may be an indication that an underlying structure is becoming increasingly more coherent as time progresses.
Indeed, in Figure~\ref{fig:delta50_full_4of4_3plots}, two structures associated with the evolution of $u^{(50)}_{62,4}$ begin entwined but evolve into clearly distinct cores. 
Each core has a well defined boundary by the time this window ends, as shown in Figure~\ref{fig:delta50_full_4of4_3plots_c}. 

We now consider $u^{(50)}_{75,4}$. The second most dominant mode for time window $W^{(50)}_{75}$ occurs as $S_{S}^{(4)}$ peaks. Such a peak occurs when the mode shifts from indicating increasing to decreasing coherency over time; thus indicating an increase in leakiness of boundaries of the associated structures. A peak such as this clearly identifies the time window over which the associated structures are most coherent.  

At this point one notes that $\alpha(0)=1$ in \eqref{eqn:SimpleModel_Forcing}, and thus the two cores are separated at $t=0$. At $\alpha(60)=0$ they have obviously merged but at $t=62$, $\alpha(62)=0.0109$ and so they have begun to separate again. Most noticeably, $\alpha(75)=0.5$ indicates the process is mid way. We also note that $\alpha(90)=1$ with $\alpha(89)=0.9973$, which indicates the initiation of movement towards merging.

By  $t_0=90$ the mode $S^{(4)}_{S}$ has moved from peaking to falling in modal dominance. The evolution of $u^{(50)}_{90,4}$ is also shown in Figure~\ref{fig:delta50_full_4of4_3plots}. This mode is initially characterised by two components with well defined boundaries that are separated in space. After $50$ time steps this mode evolves to a state of much lower coherency, as indicated by Figure~\ref{fig:delta50_logPlot}. The associated components are no longer well separated, their respective boundaries are now entwined and mass is not clearly contained within a well defined boundary. This is in agreement with the decreasing trajectory of $S^{(4)}_{S}$ around this time. 

\begin{figure}[h]
\setlength{\tabcolsep}{-0.5pt} 
\begin{tabularx}{\columnwidth}{|p{0.5cm}| *3{>{\Centering}X}| *3{>{\Centering}X}|}\hline 
&
\multicolumn{3}{c|}{\cellcolor{tn1!30}}&
\multicolumn{3}{c|}{\cellcolor{tn2!30}}\\[-1.1em]
$\: \: n$&
\multicolumn{3}{c|}{
\cellcolor{tn1!30} $51$}& 
\multicolumn{3}{c|}{
\cellcolor{tn2!30} $54$}\\
\hline
\multicolumn{1}{c|}{}&
\multicolumn{3}{c|}{
\begin{minipage}[t][][b]{0.45\textwidth}
\centering
\includegraphics[width=\columnwidth,center]{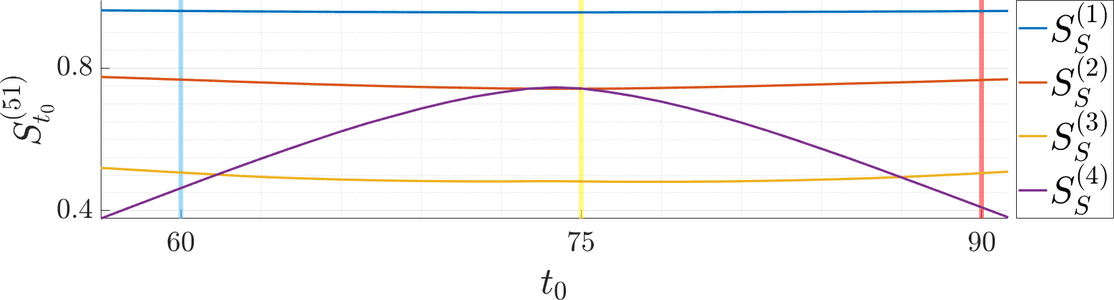}
\vspace{-0.3cm}
\end{minipage}
}&
\multicolumn{3}{c|}{
\begin{minipage}[t][][b]{0.45\textwidth}
\centering
\includegraphics[width=\columnwidth,center]{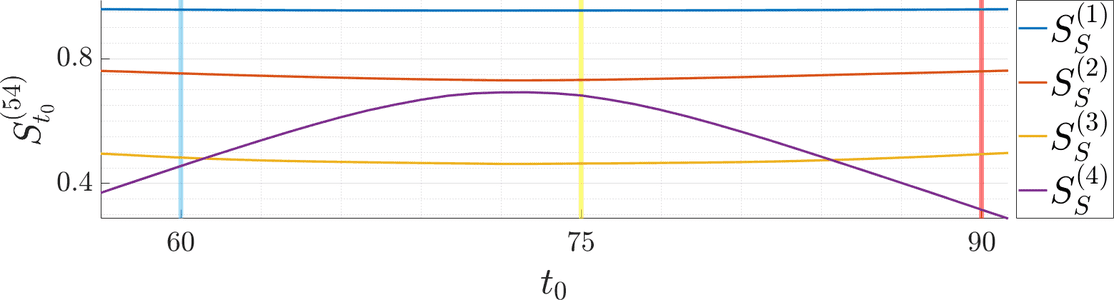}
\vspace{-0.3cm}
\end{minipage}
} 
\\ 

\hline
&
\cellcolor{t1!30} & 
\cellcolor{t2!30} & 
\cellcolor{t3!30} &           
\cellcolor{t1!30} & 
\cellcolor{t2!30} & 
\cellcolor{t3!30}\\[-1.1em]
$\: \: t_{0}$ &
\cellcolor{t1!30} $60$& 
\cellcolor{t2!30} $75$& 
\cellcolor{t3!30} $90$&           
\cellcolor{t1!30} $60$& 
\cellcolor{t2!30} $75$& 
\cellcolor{t3!30} $90$\\
\hline \hline 
\begin{minipage}[t][][b]{0.03\textwidth}
\centering
\cellcolor{S2!60}
\vspace{0.5cm}
\footnotesize{$\: \: S^{(2)}_{S}$} 
\end{minipage}
& 
\begin{minipage}[t][][b]{0.13\textwidth}
	\centering
	\vspace{-0.25cm}
	\includegraphics[height=1.8cm,center]{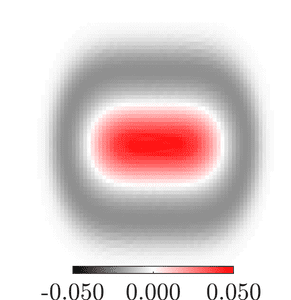}
	\vspace{-0.3cm}
\end{minipage}
& 
\begin{minipage}[t][][b]{0.13\textwidth}
	\centering
	\vspace{-0.25cm}
	\includegraphics[height=1.8cm,center]{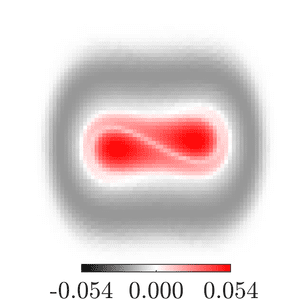}
	\vspace{-0.3cm}
\end{minipage}
& 
\begin{minipage}[t][][b]{0.13\textwidth}
	\centering
	\vspace{-0.25cm}
	\includegraphics[height=1.8cm,center]{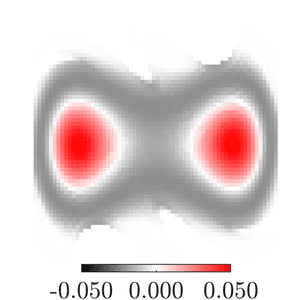}
	\vspace{-0.3cm}
\end{minipage}
&
\begin{minipage}[t][][b]{0.13\textwidth}
	\centering
	\vspace{-0.25cm}
	\includegraphics[height=1.8cm,center]{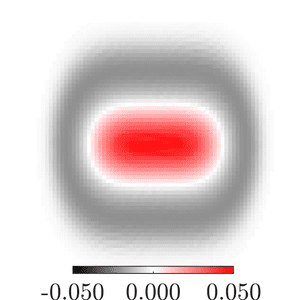}
	\vspace{-0.3cm}
\end{minipage}
& 
\begin{minipage}[t][][b]{0.13\textwidth}
	\centering
	\vspace{-0.25cm}
	\includegraphics[height=1.8cm,center]{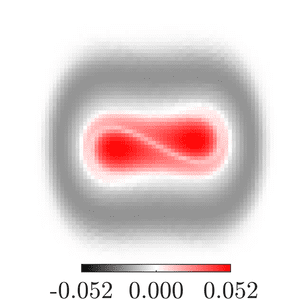}
	\vspace{-0.3cm}
\end{minipage}
& 
\begin{minipage}[t][][b]{0.13\textwidth}
	\centering
	\vspace{-0.25cm}
	\includegraphics[height=1.8cm,center]{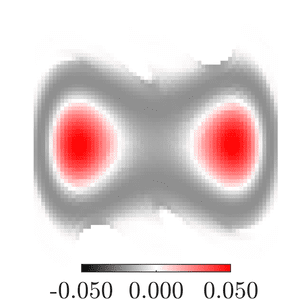}
	\vspace{-0.3cm}
\end{minipage}
\vspace{-0.2cm}
\\ \hline 
\begin{minipage}[t][][b]{0.03\textwidth}
\centering
\cellcolor{S4!60}
\vspace{0.5cm}
\footnotesize{$\: \: S^{(4)}_{S}$}
\end{minipage}
& 
\begin{minipage}[t][][b]{0.13\textwidth}
	\centering
	\vspace{-0.25cm}
	\includegraphics[height=1.8cm,center]{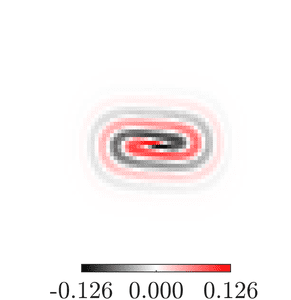}
	\vspace{-0.3cm}
\end{minipage}
& 
\begin{minipage}[t][][b]{0.13\textwidth}
	\centering
	\vspace{-0.25cm}
	\includegraphics[height=1.8cm,center]{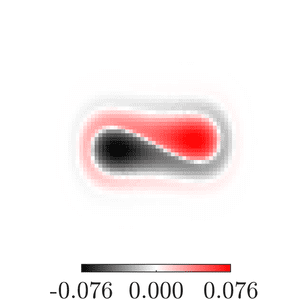}
	\vspace{-0.3cm}
\end{minipage}
& 
\begin{minipage}[t][][b]{0.13\textwidth}
	\centering
	\vspace{-0.25cm}
	\includegraphics[height=1.8cm,center]{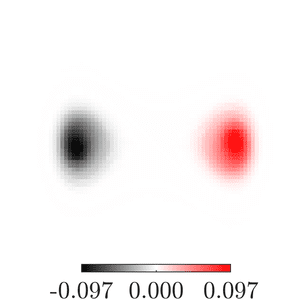}
	\vspace{-0.3cm}
\end{minipage}
&
\begin{minipage}[t][][b]{0.13\textwidth}
	\centering
	\vspace{-0.25cm}
	\includegraphics[height=1.8cm,center]{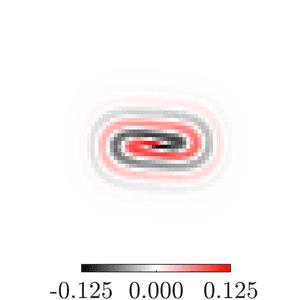}
	\vspace{-0.3cm}
\end{minipage}
& 
\begin{minipage}[t][][b]{0.13\textwidth}
	\centering
	\vspace{-0.25cm}
	\includegraphics[height=1.8cm,center]{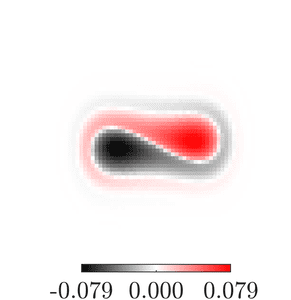}
	\vspace{-0.3cm}
\end{minipage}
& 
\begin{minipage}[t][][b]{0.13\textwidth}
	\centering
	\vspace{-0.25cm}
	\includegraphics[height=1.8cm,center]{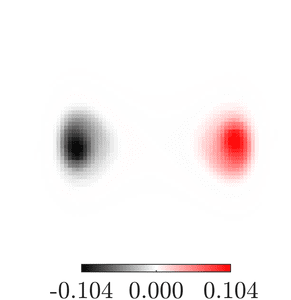}
	\vspace{-0.3cm}
\end{minipage}
\tabularnewline \hline
\end{tabularx}
\caption{Crossing introduced by shifting from $n=54$ to $n=51$ for the periodically forced double well potential.}
\label{fig:cross}
\end{figure}

\paragraph*{\textit{Analysing robustness of the pairings}.} \label{sssec:DNCs}
Let us examine the robustness of the pairings given by Algorithms~\ref{alg:SValsDijk} and \ref{alg:L2VecDist}, by analysing a situation where a small change in time window length $n$ introduces a crossing in the singular value plots. A crossing is said to occur when at least two modes intersect and shift in comparative dominance. Figure~\ref{fig:cross} illustrates a crossing that is introduced by shifting from $n=54$ to $n=51$ whilst keeping other parameters as above. Here modes are tracked using Algorithm~\ref{alg:SValsDijk} but Algorithm~\ref{alg:L2VecDist} produces similar results. Aside from the change introduced about $t_0=75$, the results for $n=54$ and $n=51$ are similar to those for $n=50$, as shown in Figure~\ref{fig:delta50_full_4of4_b}.

Of interest is how modes of time windows initialised at various times, say $t_{0}=60,75$ and $90$, behave for both $n=51$ and $54$. In particular we are interested in how the tracked mode $S_{S}^{(4)}$ relates to coherent structures in this system. Results for this case are presented in Figure~\ref{fig:cross}.
The left half of Figure~\ref{fig:cross} shows the initial time position of structures in the tracked modes associated with $S^{(2)}_{S}$ and $S^{(4)}_{S}$ for $n=51$. The right half illustrates the same for $n=54$. It is clear that the tracked mode $S_{S}^{(4)}$ of $n=54$ shifts through the third and fourth most dominant modes for $t_0 \in [58,90]$. This path never reaches the subdominant mode. The left columns show how the tracked mode $S_{S}^{(4)}$ of $n=51$ shifts all the way from the fourth most dominant position to the second over the same time frame. 

\vspace{-0.25cm}
Separating cores are associated with time windows in the leftmost column of each half, whilst merging cores are associated with those in the rightmost columns of each half of the table. The middle subcolumns capture the period over which the two main components are the most individually coherent, as neither core separates nor merges as dramatically as structures in surrounding time windows. This time window occurs as $S_{S}^{(4)}$ (Figure~\ref{fig:cross}, purple modes in electronic version) peaks in singular value. Structures characterised by singular vectors associated with singular values closer to one are expected to evolve coherently, as the larger the singular value, the longer the associated structure is expected to survive. It can been seen in Figure~\ref{fig:delta50_full_4of4_b} that the peak in $S_{S}^{(4)}$ occurs around the time that $S_{S}^{(2)}$ acheives a global minimum. 

Singular vectors associated with $S_{S}^{(4)}$ are linked by a propensity to separate into similarly weighted positive and negative elements. The location of these elements is associated with the position of two cores of interest. Background noise is associated with vector elements close to zero. This appears white in Figure~\ref{fig:cross}. Merging and separating cores can be tracked over time in a way that isolates two distinct components, even as their boundaries are fundamentally altered. 

Singular vectors associated with $S_{S}^{(2)}$ are linked by a propensity to unify the core elements as two halves of the one component. The center of this structure (red in the electronic version) is least coherent around time $75$, when its single boundary begins to weaken. It is most coherent around time windows initialised about time $25$. Whilst this is not shown in Figure~\ref{fig:cross}, it is evident in the analogous example provided by Figure~\ref{fig:delta50_full_4of4_b} for time windows of length $n=50$. When considering a single component, one expects this to be most coherent following a merger of the two previously seperate cores that define it. Indeed, it is clear from Figure~\ref{fig:time_DWPV2} that time windows of lengths around $50$ that are initialised at $t_{0}=25$ are most associated with a merging event. 

\begin{figure}[H]
\centering
\begin{minipage}{\textwidth}
\includegraphics[width=\textwidth,center]{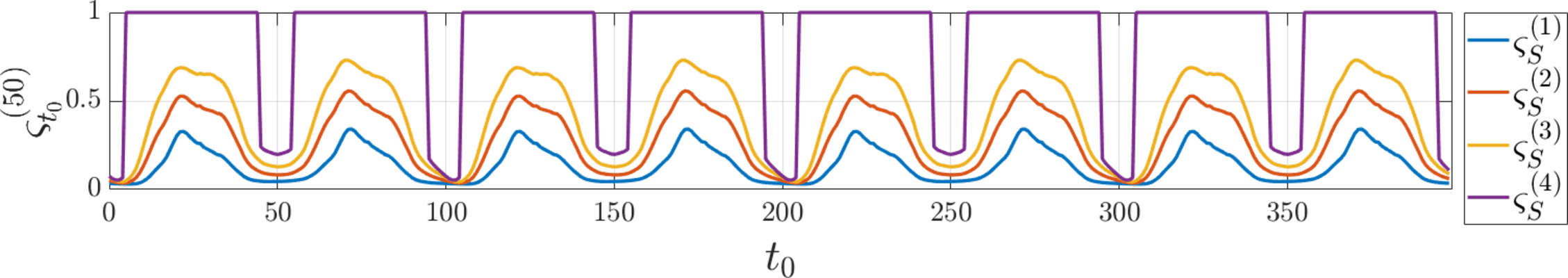}
\subcaption{Equivariance of leading $4$ of $\mathcal{N}=4$ paired modes $S_{S}$ as per Algorithms~\ref{alg:Eqvr} and~\ref{alg:SValsDijk}.}
\label{fig:deltaT50_Eq_Paired4of4}
\end{minipage}
\begin{minipage}{\textwidth}
\includegraphics[width=\textwidth,center]{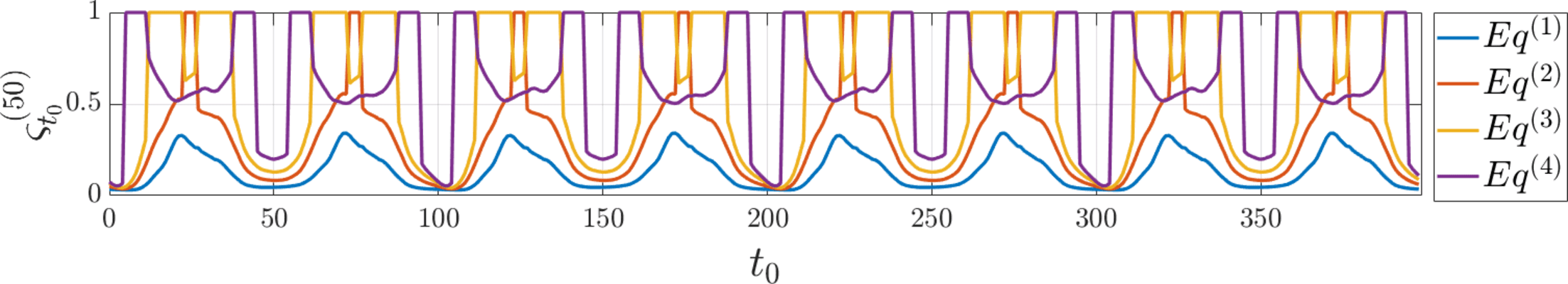}
\subcaption{Equivariance of leading $4$ of $\mathcal{N}=4$ unpaired modes.}
\label{fig:deltaT50_Eq_Unpaired}
\end{minipage}
\caption{Equivariance mismatch for the periodically forced double well potential when $n=50$.}\label{fig:deltaT50_Eq}
\end{figure}

Thus, while changes in dominance (crossings) are sensitive to changes in window length, the patterns of increasing and decreasing phases of singular values, as well as their associated singular vectors, appear robust. These phases are therefore a more reliable feature to explore when investigating fundamental aspects of coherent structures. 
\newline
\paragraph*{\textit{Testing equivariance}.}\label{sssec:ET} The accuracy of mode pairings from Algorithms~\ref{alg:SValsDijk} and \ref{alg:L2VecDist} can be assessed using the Equivariance test described in Algorithm~\ref{alg:Eqvr}. For time windows of length $n=50$, the results for pairing via Algorithm~\ref{alg:SValsDijk} are presented in Figure~\ref{fig:deltaT50_Eq_Paired4of4}. In this case results are the same as for Algorithm~\ref{alg:L2VecDist}. 
Recall that equivariance mismatch values away from one indicate effective pairing. Periods over which equivariance jumps rapidly between low and high values indicate an inconsistency of pairing between neighbouring time windows. 
\vspace{-0.5cm}
The efficacy of our pairing algorithms is evident from a comparison of Figures~\ref{fig:deltaT50_Eq_Paired4of4} and \ref{fig:deltaT50_Eq_Unpaired}. Indeed,  Algorithms~\ref{alg:SValsDijk} and \ref{alg:L2VecDist} are able to correct for the mismatch that occurs when mode dominance is assumed constant through time. In Figure~\ref{fig:deltaT50_Eq_Paired4of4} it is only the fourth mode that takes a value close to $1$. This mode is only well matched over very small windows. Such windows begin around the time when the two core structures shift from having merged to being separate or vice versa.  For example, one period of effective pairing for $S_{S}^{(4)}$ coincides with the mode rising and falling through other modes as shown in Figure~\ref{fig:delta50_full_4of4}. 

\begin{figure}[H]
\centering
\begin{minipage}{\textwidth}
\centering
\includegraphics[width=\textwidth,center]{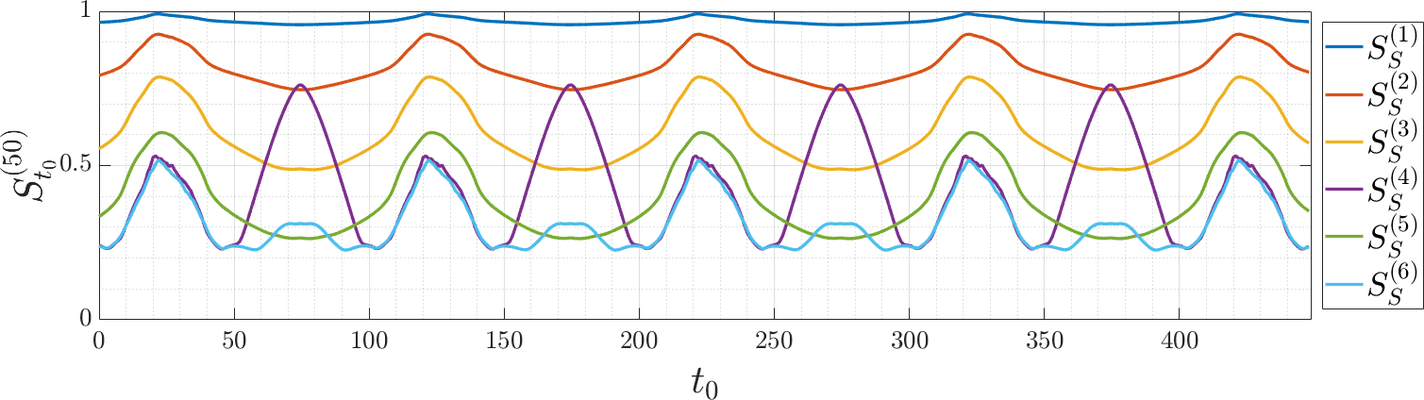}
\subcaption{Paths tracked by $S_{S}$ for $n=50$ using Algorithm~\ref{alg:SValsDijk}.}
\label{fig:delta50_top6_wind_full_a}
\end{minipage}
\begin{minipage}{\textwidth}
\centering
\includegraphics[width=\textwidth,center]{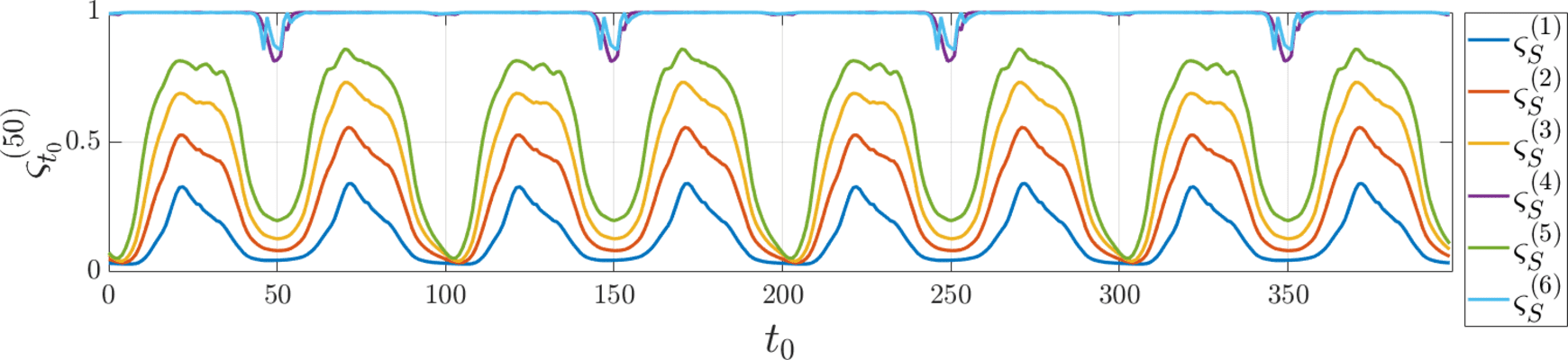}
\subcaption{Equivariance mismatch using Algorithm~\ref{alg:Eqvr} for paths tracked as per Figure~\ref{fig:delta50_top6_wind_full_a}.}
\label{fig:delta50_top6_wind_full_b}
\end{minipage}
\caption{Leading $6$ of $\mathcal{N}=6$ modes for the periodically forced double well potential when $n=50$.}
\label{fig:delta50_top6_wind_full}
\end{figure}

In order to explore the possibility of better pairing options for modes, one may track a larger number of singular values. Consider tracking $\mathcal{N}=6$ modes using Algorithm~\ref{alg:SValsDijk}. The tracked paths of for these parameters, along with the corresponding equivariance measures, are shown in Figure~\ref{fig:delta50_top6_wind_full}.
\vfil \hfil \vspace{-1cm}
\vspace{-0.25cm}
In this case it is clear that a fourth mode has now been satisfactorily tracked, in the sense that the corresponding equivariance mismatch measure remains bounded away from one for all times. As expected, this mode partially matches that of Figure~\ref{fig:delta50_full_4of4}.
On the other hand, one notes that as paths take on extremely close singular values, it becomes more difficult to clearly separate them and track them through time. 
For example, consider the paths tracked by $S_{S}^{(4)}$ and $S_{S}^{(6)}$ in Figure~\ref{fig:delta50_top6_wind_full_a} over $t_0 \in [10,40]$ or $[110,140]$, where multiple changes in dominance occur. In such cases, it is not clear which path the mode follows. This is in agreement with Figure~\ref{fig:delta50_top6_wind_full_b} which shows that neither $S_{S}^{(4)}$ nor $S_{S}^{(6)}$ is well paired over the majority of time windows. 

In general, the occurrence of very close singular values becomes more prevalent as one considers larger values of $\mathcal N$ in Algorithms~\ref{alg:SValsDijk} and \ref{alg:L2VecDist}.
This may sometimes be due to numerical errors or low resolution, but it may also be due to the fact that the associated modes no longer correspond to meaningful dynamical features of the underlying system.

\subsubsection{Rolling time windows of length $n=100$} \label{sssec:TSOT_100}
\begin{figure}[H]
\centering
\begin{minipage}{\textwidth}
\centering
\includegraphics[width=\textwidth,center]{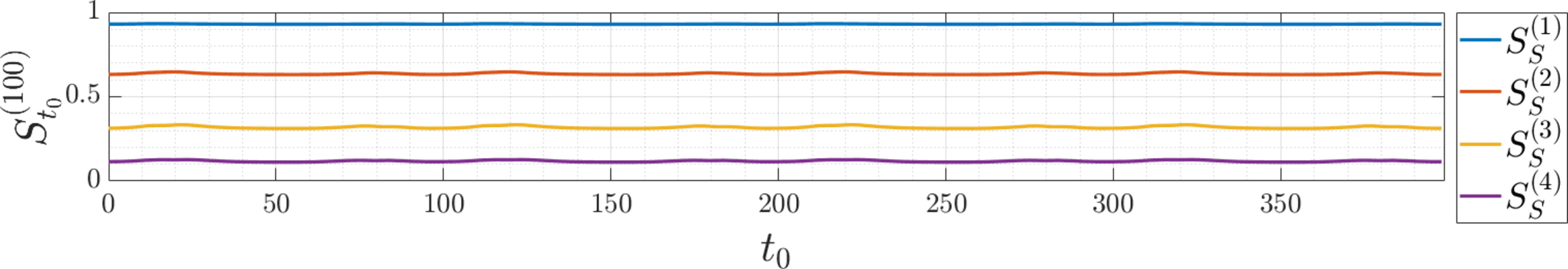}
\subcaption{Paths tracked by $S_{S}$ for $n=100$ using Algorithm~\ref{alg:SValsDijk}.}
\label{fig:delta100_wind399}
\end{minipage}
\begin{minipage}{\textwidth}
\centering
\includegraphics[width=\textwidth,center]{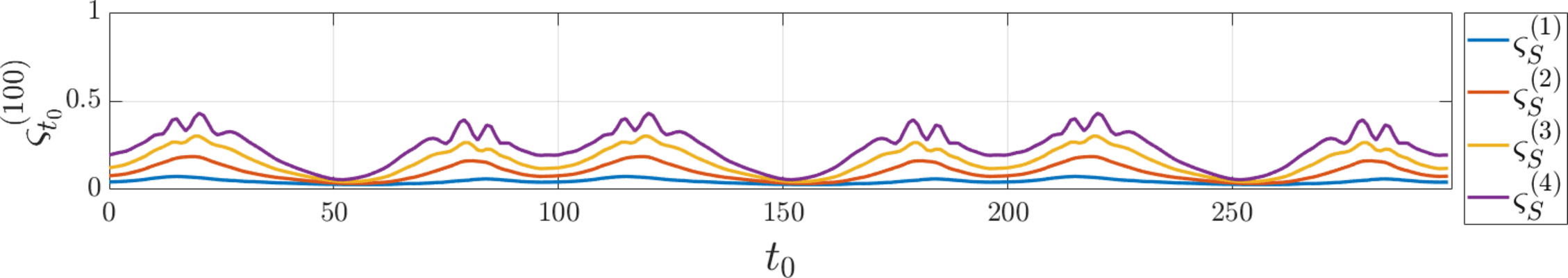}
\subcaption{Equivariance mismatch using Algorithm~\ref{alg:Eqvr} for paths tracked as per Figure~\ref{fig:delta100_wind399}.}
\label{fig:deltaT100_Eq_Paired_4of4}
\end{minipage}
\caption{Leading $4$ of $\mathcal{N}=4$ modes for the periodically forced double well potential when $n=100$.}
\label{fig:deltaT100_4of4}
\end{figure}

Increasing the value of $n$ in Algorithms~\ref{alg:SValsDijk} and \ref{alg:L2VecDist} incorporates information corresponding to longer time periods in each SVD calculation. This also makes it more likely that the singular values, approximating Lyapunov exponents for the cocycle, become more separated. In this model, shifting from $n=50$ to $n=100$ eradicates crossings among higher modes as the more transitory structures become less dynamically relevant. Figure~\ref{fig:delta100_wind399} plots the result of implementing Algorithm~\ref{alg:SValsDijk} on windows of length $n=100$ with the remaining parameters as per Section~\ref{sssec:TSOT}. The path of the leading singular values, $S_{S}^{(1)}$ in Figure~\ref{fig:delta100_wind399} is nearly constant at approximately $0.93$.

Each of the first four associated modes is clearly separated from the preceding mode, and because there are no crossings, the paired and unpaired modes give the same outcome. Despite their separation, these modes continue to exhibit the peaking behaviour characteristic of phases of increasing and decreasing coherence. Given a clear separation of modes, the equivariance test points to well paired modes by returning values away from one over time. Figure~\ref{fig:deltaT100_Eq_Paired_4of4} plots the outcomes of Algorithm~\ref{alg:Eqvr} in this instance. Both results presented in Figure~\ref{fig:deltaT100_4of4} are consistent with those utilising Algorithm~\ref{alg:L2VecDist}.

\vspace{-0.25cm}
In this case, for each time window, $n=100$ matrices of size $2^{12} \times 2^{12}$ are being multiplied, rather than just $n=50$ as in Section~\ref{sssec:ET}. Whilst a greater value of $n$ could increase the numerical error, the equivariance test being passed for larger $n$ conforms with the theoretical expectations associated to true coherent structures, associated to Oseledets modes, being present. Extending the number of singular values considered to $\mathcal{N}=5$ sees the maximal equivariance jump drastically, to $0.9180$. This suggests that, for longer time windows, important information regarding the more transitory dynamics is found in modes of lower dominance. 

\subsection{Quasi-periodically forced double well potential} \label{ssec:QPmodel}
In the previous section we have considered a periodic system where two core components merge and separate at periodically spaced times. Let us now introduce an additional quasi-periodic forcing to the vector field of Section~\ref{ssec:DWPmodel}. The two wells will still shift in and out over time. However, the time at which these changes occur will be less regular than before. 

\begin{figure}[H]
\centering
\includegraphics[width=1\textwidth,center]{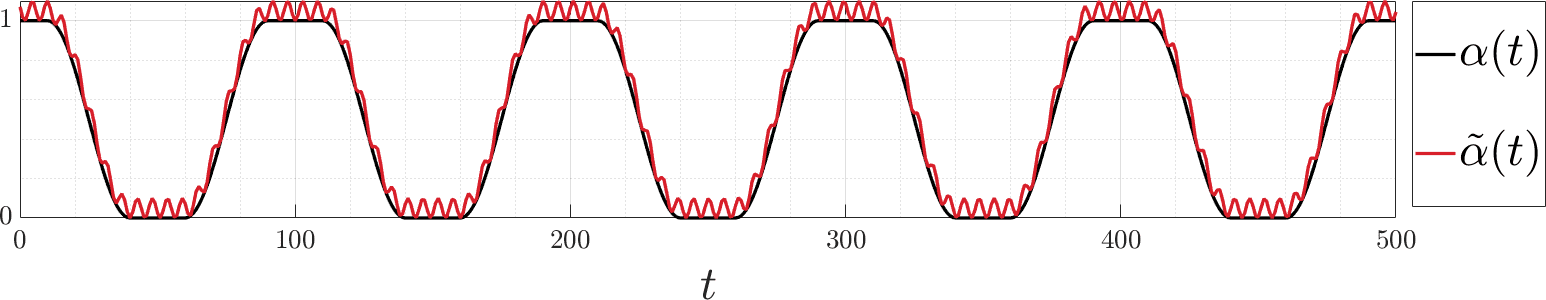}
\caption{An illustration of the behaviour of $\alpha(t)$ and $\tilde{\alpha}(t)$ over $5$ periods.}
\label{fig:time_DWPV2}
\end{figure}

For this, we replace $\alpha(t)$ in Equation~\eqref{eqn:SimpleModel} by $\tilde{\alpha}(t)=\alpha(t)+\gamma \cos^{2}(10t)$, with $\gamma=0.1$ and $\alpha(t)$ defined as in Equation~\eqref{eqn:SimpleModel_Forcing}. This driving force behaves in the manner depicted by $\tilde{\alpha}(t)$ in Figure~\ref{fig:time_DWPV2}. Under the quasi-periodic forcing given by $\tilde{\alpha}(t)$, the core of each structure moves closer to the boundaries of phase space than in the  periodic case.

\vspace{-0.25cm}
The quasi-periodic driving force introduces additional complexity by disturbing the location of the center of each core more irregularly over time. This additional complexity  affects the ability of our algorithms to effectively pair modes through time. In order to choose which pairing strategy is more effective, we rely on the measure of equivariance described in Algorithm~\ref{alg:Eqvr}. 

\begin{figure}[H]
\centering
\includegraphics[width=\columnwidth,center]{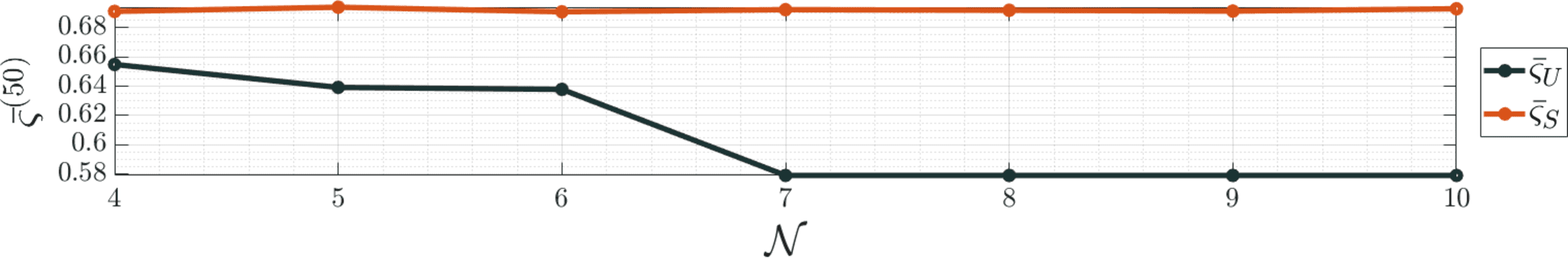}
\caption{Mean equivariance mismatch, as per Algorithm~\ref{alg:Eqvr}, for the leading $4$ of $\mathcal{N}$ modes using the two pairing methods given by Algorithms~\ref{alg:SValsDijk} ($\bar{\varsigma}_{S}$) and~\ref{alg:L2VecDist} ($\bar{\varsigma}_{U}$) for $n=50$.} \label{fig:QPDWP_comparison}
\end{figure}

Maintaining time windows of length $n=50$, the efficacy of Algorithms~\ref{alg:SValsDijk} and \ref{alg:L2VecDist} is compared. Figure~\ref{fig:QPDWP_comparison} presents a comparative summary of the average value of equivariance mismatch for the four leading modes averaged over all time. Here $\mathcal{N}$, the total number of singular values considered, varies but all other parameters remain as in Section~\ref{sssec:TSOT}.

\begin{figure}[H]
\centering
\begin{minipage}{\textwidth}
\centering
\includegraphics[width=\textwidth,center]{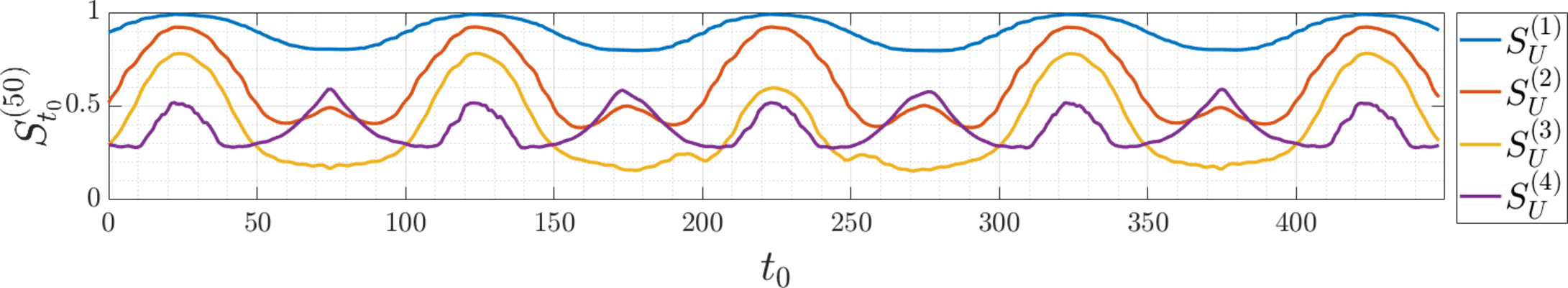}
\end{minipage}
\begin{minipage}{\textwidth}
\centering
\includegraphics[width=\textwidth,center]{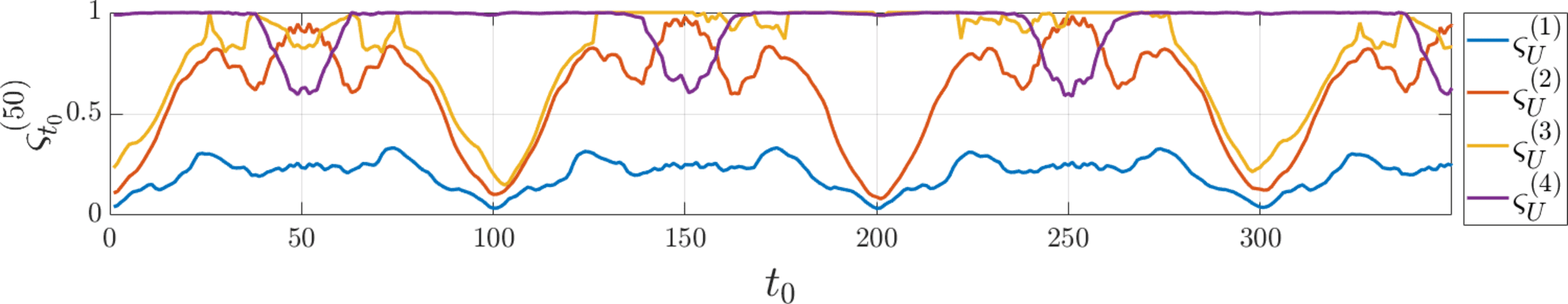}
\end{minipage}
\caption{Leading $4$ from a total $\mathcal{N}=5$ tracked modes for $n=50$ using Algorithms~\ref{alg:L2VecDist} (top) and~\ref{alg:Eqvr} (bottom).}
\label{fig:delta50_full_SVec4of5_QPDWP}
\end{figure}

\begin{figure}[H]
\centering
\begin{minipage}{\textwidth}
\centering
\includegraphics[width=\textwidth,center]{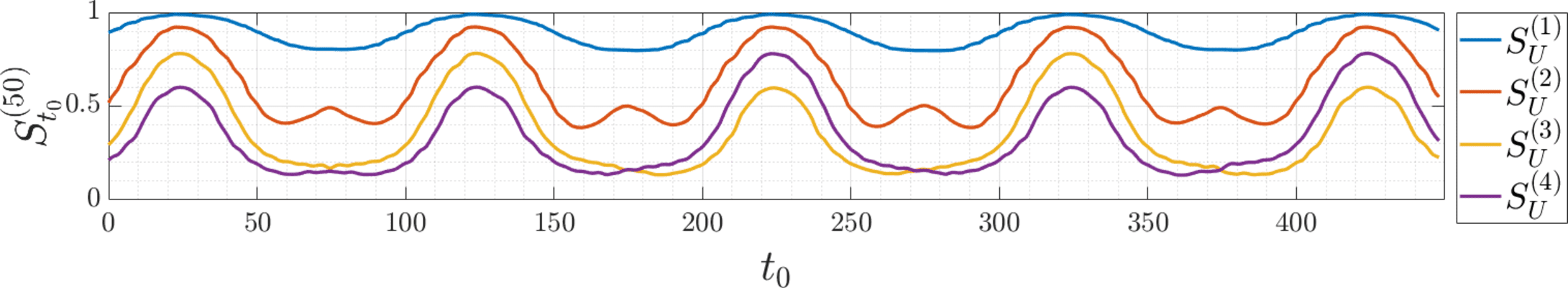}
\end{minipage}
\begin{minipage}{\textwidth}
\centering
\includegraphics[width=\textwidth,center]{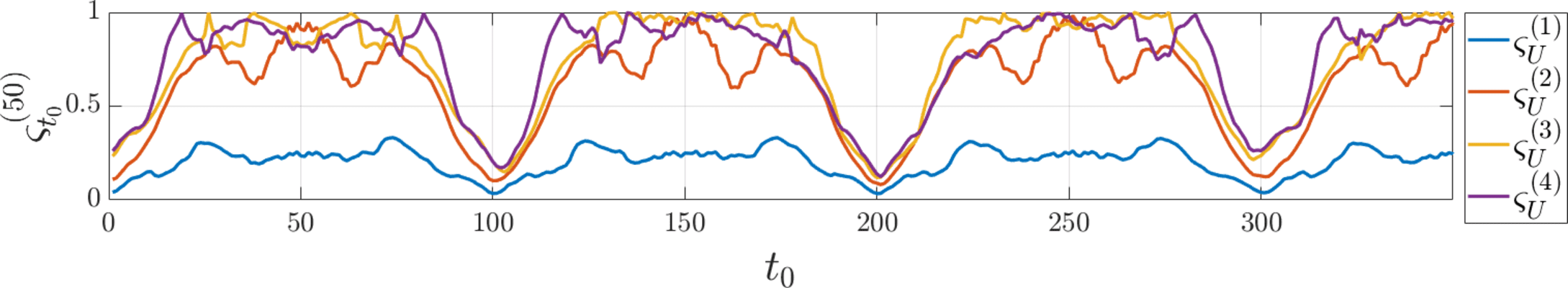}
\end{minipage}
\caption{Leading $4$ from a total $\mathcal{N}=7$ tracked modes for $n=50$ using Algorithms~\ref{alg:L2VecDist} (top) and~\ref{alg:Eqvr} (bottom).}
\label{fig:delta50_full_SVec4of7_QPDWP}
\end{figure}

\vspace{-0.25cm}
In this case, $\bar{\varsigma}_{U}$ consistently attains lower values than $\bar{\varsigma}_{S}$, with $\bar{\varsigma}_{U}$ initially plateauing out at $\mathcal{N}=5$. The mean equivariance mismatch $\bar{\varsigma}_{U}$, then moves to a lower minimum when $\mathcal{N} \ge 7$. Let us further explore those values of $\mathcal{N}$ that initiate each plateau. Results utilising $\mathcal{N}=5$ are presented in Figure~\ref{fig:delta50_full_SVec4of5_QPDWP}, those for $\mathcal{N}=7$ are presented in Figure~\ref{fig:delta50_full_SVec4of7_QPDWP}. 
\vspace{-0.25cm}
\vspace{-0.25cm}
In both Figures~\ref{fig:delta50_full_SVec4of5_QPDWP} and~\ref{fig:delta50_full_SVec4of7_QPDWP} the leading mode is paired consistently through time. This is indicated by an equivariance mismatch value far from $1$. The one anomaly in Figure~\ref{fig:delta50_full_SVec4of7_QPDWP} is a switching of modes $S_{U}^{(3)}$ and $S_{U}^{(4)}$. This occurs with increasing frequency as the singular values get closer to zero and the modes themselves are not well separated. 

As before, in the search for modes that indicate the occurrence of fundamental changes in a system, we turn our attention to peaks, corresponding to transitions between increasing and decreasing coherency phases. 

Let us investigate the case $\mathcal{N}=5$ in more detail. Pairings associated with $\bar{\varsigma}_{U}$ for $\mathcal{N}\ge 7$ include modes that are less well separated over time, such as $S_U^{(3)}$ and $S_U^{(4)}$ around time $175$ in Figure~\ref{fig:delta50_full_SVec4of7_QPDWP}. 
 
\begin{figure}[H]
\includegraphics[width=\columnwidth,center]{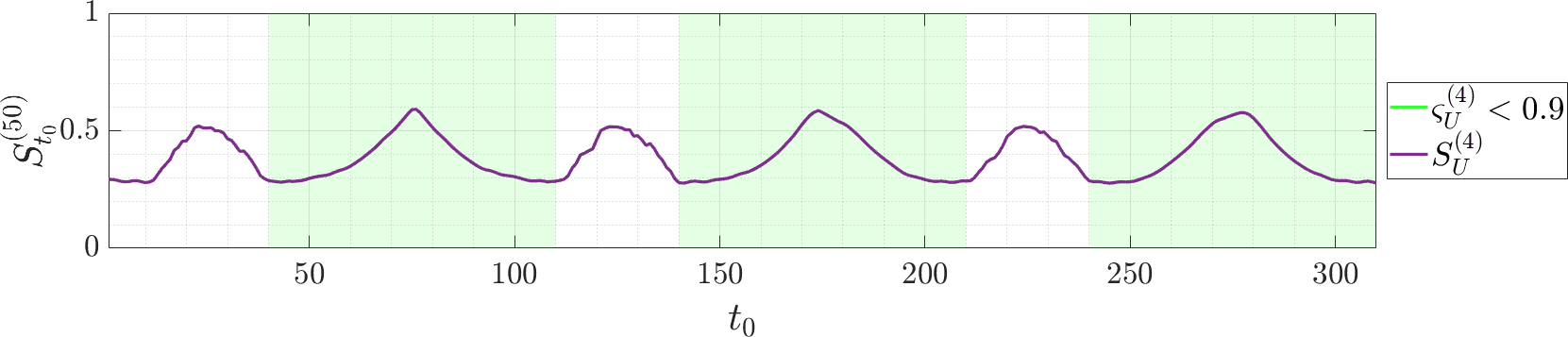}
\caption{Consecutive windows corresponding to reasonable equivariance for $S_{U}^{(4)}$ of Figure~\ref{fig:delta50_full_SVec4of5_QPDWP}.}
\label{fig:delta50_full_SVec5of5_DWPV2_withLines_Eq}
\end{figure}

\vspace{-0.25cm}
The fact that peaks in $S_{U}^{(4)}$ of Figure~\ref{fig:delta50_full_SVec4of5_QPDWP} occur at times less evenly spaced than those of Section~\ref{ssec:DWPmodel} reflects the quasi-periodic nature of this system. One notes that $S_{U}^{(4)}$ of Figure~\ref{fig:delta50_full_SVec4of5_QPDWP} is well paired over the period when peaks develop. This is illustrated in Figure~\ref{fig:delta50_full_SVec5of5_DWPV2_withLines_Eq}. 

\begin{figure}[H]
\setlength{\tabcolsep}{-0.5pt} 
\begin{tabularx}{\columnwidth}{|p{0.5cm}| *6{>{\Centering}X}|}\hline 
\multicolumn{7}{|c|}{
\begin{minipage}[t][][b]{\columnwidth}
\centering
\vspace{-0.1cm}
\includegraphics[width=0.99\columnwidth,center]{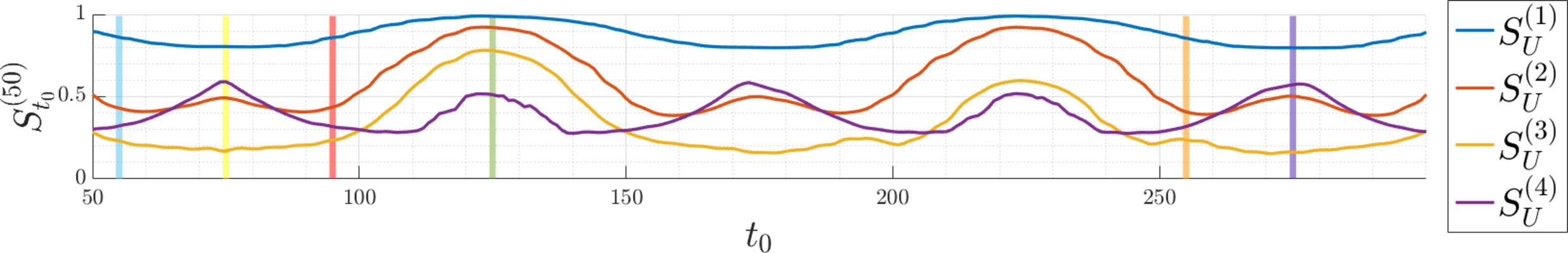}
\vspace{-0.3cm}
\end{minipage} }\\
 \hline \hline
&
\cellcolor{t1!30} & 
\cellcolor{t2!30} & 
\cellcolor{t3!30} &           
\cellcolor{t4!30} & 
\cellcolor{t5!30} & 
\cellcolor{t6!30} \\[-1.1em]
$\: \: t_{0}$ &
\cellcolor{t1!30} $55$& 
\cellcolor{t2!30} $75$& 
\cellcolor{t3!30} $95$&           
\cellcolor{t4!30} $125$& 
\cellcolor{t5!30} $255$& 
\cellcolor{t6!30} $275$ 
\\ 
\hline \hline 
\begin{minipage}[t][][b]{0.03\textwidth}
\centering
\cellcolor{S1!60} 
\vspace{0.5cm}
\footnotesize{$\: \: S^{(1)}_{U}$}
\end{minipage}
& 
\begin{minipage}[t][][b]{0.12\textwidth}
	\centering
	\vspace{-0.25cm}
	\includegraphics[height=1.8cm,center]{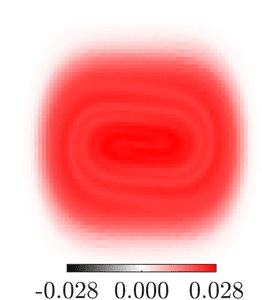}
	\vspace{-0.3cm}
\end{minipage}
& 
\begin{minipage}[t][][b]{0.12\textwidth}
	\centering
	\vspace{-0.25cm}
	\includegraphics[height=1.8cm,center]{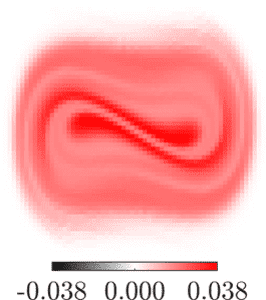}
\end{minipage}
& 
\begin{minipage}[t][][b]{0.12\textwidth}
	\centering
	\vspace{-0.25cm}
	\includegraphics[height=1.8cm,center]{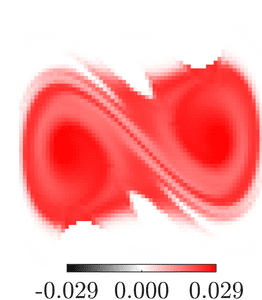}
\end{minipage}
& 
\begin{minipage}[t][][b]{0.12\textwidth}
	\centering
	\vspace{-0.25cm}
	\includegraphics[height=1.8cm,center]{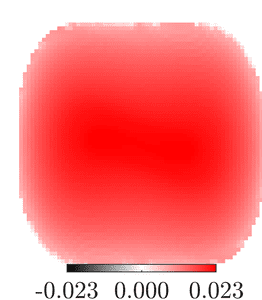}
\end{minipage}
& 
\begin{minipage}[t][][b]{0.12\textwidth}
	\centering
	\vspace{-0.25cm}
	\includegraphics[height=1.8cm,center]{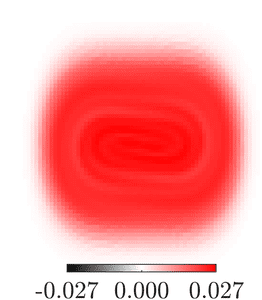}
\end{minipage}
& 
\begin{minipage}[t][][b]{0.12\textwidth}
	\centering
	\vspace{-0.25cm}
	\includegraphics[height=1.8cm,center]{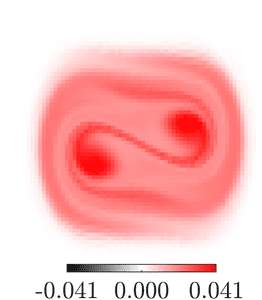}
\end{minipage}
\\  \hline 
\vspace{0cm}
\begin{minipage}[t][][b]{0.03\textwidth}
\centering
\cellcolor{S2!60}
\vspace{0.2cm}
\footnotesize{$\: \: S^{(2)}_{U}$}
\end{minipage}
& 
\begin{minipage}[t][][b]{0.12\textwidth}
	\centering
	\vspace{-0.25cm}
	\includegraphics[height=1.8cm,center]{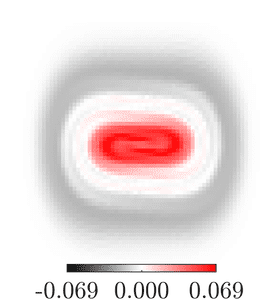}
	\vspace{-0.3cm}
\end{minipage} \hfil
& 
\begin{minipage}[t][][b]{0.12\textwidth}
	\centering
	\vspace{-0.25cm}
	\includegraphics[height=1.8cm,center]{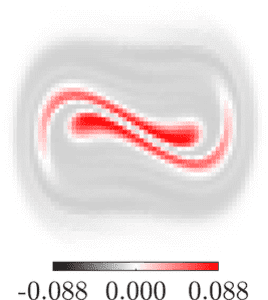}
\end{minipage}
& 
\begin{minipage}[t][][b]{0.12\textwidth}
	\centering
	\vspace{-0.25cm}
	\includegraphics[height=1.8cm,center]{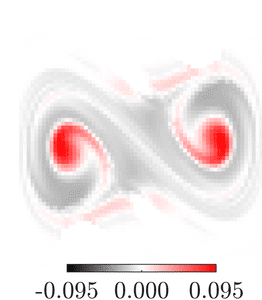}
\end{minipage}
& 
\begin{minipage}[t][][b]{0.12\textwidth}
	\centering
	\vspace{-0.25cm}
	\includegraphics[height=1.8cm,center]{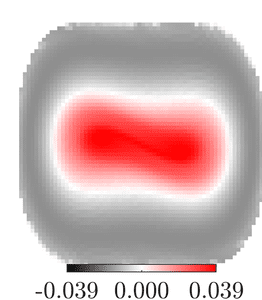}
\end{minipage}
&
\begin{minipage}[t][][b]{0.12\textwidth}
	\centering
	\vspace{-0.25cm}
	\includegraphics[height=1.8cm,center]{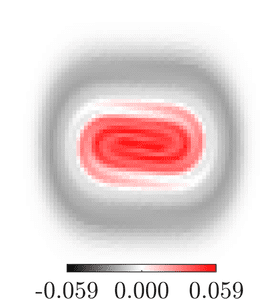}
\end{minipage}
& 
\begin{minipage}[t][][b]{0.12\textwidth}
	\centering
	\vspace{-0.25cm}
	\includegraphics[height=1.8cm,center]{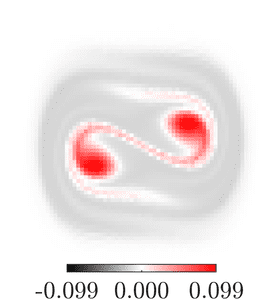}
\end{minipage}
\\ \hline 
\vspace{0cm}
\begin{minipage}[t][][b]{0.03\textwidth}
\centering
\cellcolor{S3!60}
\vspace{0.2cm}
\footnotesize{$\: \: S^{(3)}_{U}$} 
\end{minipage}
& 
\begin{minipage}[t][][b]{0.12\textwidth}
	\centering
	\vspace{-0.25cm}
	\includegraphics[height=1.8cm,center]{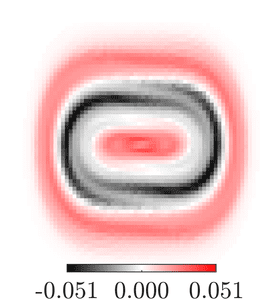}
	\vspace{-0.3cm}
\end{minipage}
& 
\begin{minipage}[t][][b]{0.12\textwidth}
	\centering
	\vspace{-0.25cm}
	\includegraphics[height=1.8cm,center]{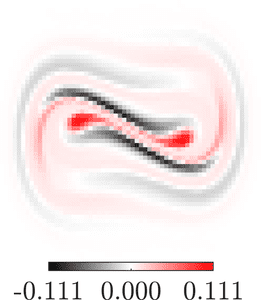}
\end{minipage}
& 
\begin{minipage}[t][][b]{0.12\textwidth}
	\centering
	\vspace{-0.25cm}
	\includegraphics[height=1.8cm,center]{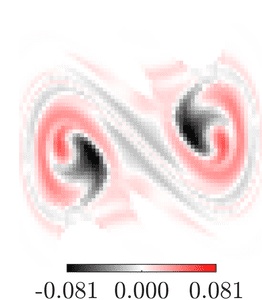}
\end{minipage}
& 
\begin{minipage}[t][][b]{0.12\textwidth}
	\centering
	\vspace{-0.25cm}
	\includegraphics[height=1.8cm,center]{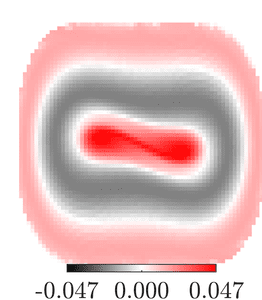}
\end{minipage}
&
\begin{minipage}[t][][b]{0.12\textwidth}
	\centering
	\vspace{-0.25cm}
	\includegraphics[height=1.8cm,center]{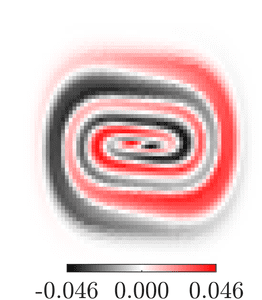}
\end{minipage}
& 
\begin{minipage}[t][][b]{0.12\textwidth}
	\centering
	\vspace{-0.25cm}
	\includegraphics[height=1.8cm,center]{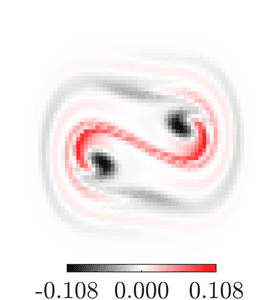}
\end{minipage}
\\  \hline 
\vspace{0cm}
\begin{minipage}[t][][b]{0.03\textwidth}
\centering
\cellcolor{S4!60}
\vspace{0.2cm}
\footnotesize{$\: \: S^{(4)}_{U}$}
\end{minipage}
& 
\begin{minipage}[t][][b]{0.12\textwidth}
	\centering
	\vspace{-0.25cm}
	\includegraphics[height=1.8cm,center]{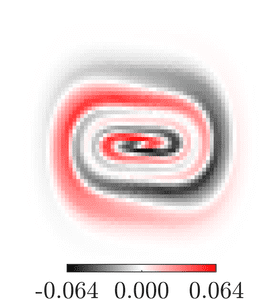}
	\vspace{-0.3cm}
\end{minipage}
& 
\begin{minipage}[t][][b]{0.12\textwidth}
	\centering
	\vspace{-0.25cm}
	\includegraphics[height=1.8cm,center]{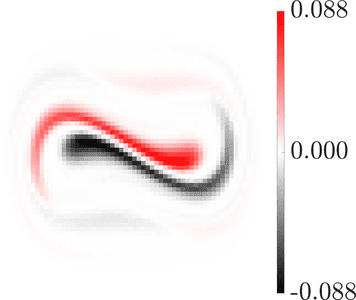}
\end{minipage}
& 
\begin{minipage}[t][][b]{0.12\textwidth}
	\centering
	\vspace{-0.25cm}
	\includegraphics[height=1.8cm,center]{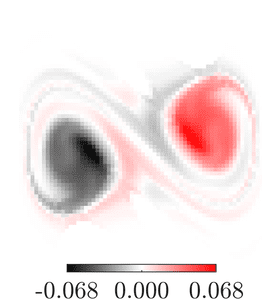}
\end{minipage}
& 
\begin{minipage}[t][][b]{0.12\textwidth}
	\centering
	\vspace{-0.25cm}
	\includegraphics[height=1.8cm,center]{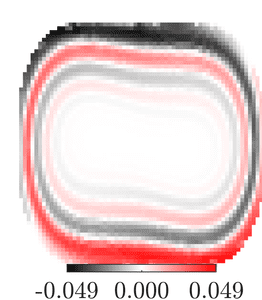}
\end{minipage}
&
\begin{minipage}[t][][b]{0.12\textwidth}
	\centering
	\vspace{-0.25cm}
	\includegraphics[height=1.8cm,center]{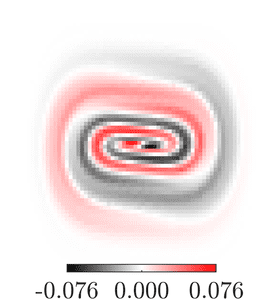}
\end{minipage}
& 
\begin{minipage}[t][][b]{0.12\textwidth}
	\centering
	\vspace{-0.25cm}
	\includegraphics[height=1.8cm,center]{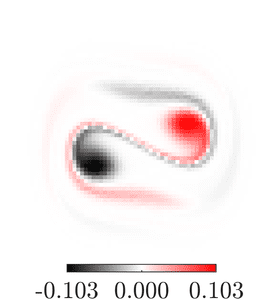}
\end{minipage}
\tabularnewline \hline
\end{tabularx}
\caption{Initial time singular vectors corresponding to rolling windows initialised at the various $t_{0}$ indicated by colour coded bars and column headings. These are paired according to the paths illustrated in Figure~\ref{fig:delta50_full_SVec4of5_QPDWP}.} \label{fig:main_table_QPDWP}
\end{figure}

\vspace{-0.25cm}
Let us further explore the behaviour of modes illustrated in Figure~\ref{fig:delta50_full_SVec5of5_DWPV2_withLines_Eq} by examining the initial time singular vectors for a variety of initial times $t_{0}$. A selection of relevant findings for $t_{0}=55,75,95,125,255$ and $275$ are presented in Figure~\ref{fig:main_table_QPDWP}.
\vspace{-0.25cm}
As in Section~\ref{ssec:QPmodel}, $S_{U}^{(1)}$ of Figure~\ref{fig:main_table_QPDWP} is generally well separated from lower modes. An examination of the corresponding initial time singular vectors, presented in the row of Figure~\ref{fig:main_table_QPDWP} associated with $S_{U}^{(1)}$, shows that the leading mode identifies the general location of dominant structures in the system. 

Figure~\ref{fig:delta50_full_SVec4of5_QPDWP} shows that $S_{U}^{(2)}$ achieves an equivariance value that is generally less than $0.8$ over the considered initial times. An examination of the left singular vectors corresponding to $S_{U}^{(2)}$ serves to reinforce the notion that this pairing is satisfactory. The qualitative consistency of singular vectors illustrated in the row corresponding to $S_{U}^{(2)}$ in Figure~\ref{fig:main_table_QPDWP} demonstrates the efficacy of this pairing. It is clear from this figure that the mode associated with $S_{U}^{(2)}$ isolates the core as one component rather than two distinct entities. As in~\ref{sssec:DNCs}, this structure is surrounded by a constantly mixing shell, even as it separates in two. One notes that $\tilde{\alpha}(125)=0.5952$ with $\tilde{\alpha}(150)=0.0201$ and $\tilde{\alpha}(175)=0.5577$, thus this structure achieves peak coherency over time windows starting near to local maxima of singular values, similar to $W^{(50)}_{125}$.

To the left of $W^{(50)}_{125}$, the time window $W^{(50)}_{95}$ characterises a period of increasing coherency for the structure associated with $S_{U}^{(2)}$. One notes that $\tilde{\alpha}(95)=1.006$ but $\tilde{\alpha}(145)=0.0445$ as the boundary is restored and time progresses. To the right, $W^{(50)}_{255}$ defines a full seperation event. This is characterised by a period of lower coherency where the boundaries of the structure core are manipulated and stretched through time, eventually resulting in a division into separate components.

Singular vectors associated with $S_{U}^{(3)}$ attain a lower equivariance than those of $S_{U}^{(4)}$ but this pairing is less smooth across the various $t_{0}$. For example, $\varsigma_{U}^{(3)}$ in Figure~\ref{fig:delta50_full_SVec4of5_QPDWP} peaks rather sharply at $t_0=63,75$ and $127$. Singular vectors associated with $S_{U}^{(3)}$, as illustrated in the corresponding row of Figure~\ref{fig:main_table_QPDWP}, show that these peaks could indicate an undetected change in pairing regime. Singular vectors associated with $S_{U}^{(3)}$ identify smaller substructures in this system. These are less consequential for the global dynamics and more often associated with the dynamics within each core rather than with the interactions between cores.

The mode associated with $S_{U}^{(4)}$ is of special interest. Of note is the tendency for it to separate phase space into near equally weighted positive and negative components in a way that higher modes do not. As Figure~\ref{fig:delta50_full_SVec5of5_DWPV2_withLines_Eq} made clear, this mode is not well paired over all time windows. However, it is well paired over periods associated with peaking behaviour. In Figure~\ref{fig:main_table_QPDWP} this would include $t_{0}=75,173$ and $276$. For time windows initialised when the two components tend towards mixing, $S_{U}^{(4)}$ falls in dominance and the mode is difficult to track.

The various time windows considered in Figure~\ref{fig:main_table_QPDWP} include those about two fundamentally distinct peaks in $S_{U}^{(4)}$. The sharper of these peaks occurs at $t_{0}=75$ whilst a more rounded and longer lasting peak occurs at $t_{0}=276$. The sharper peak attains the higher maximum of $S_{U}^{(4)}={0.5908}$ at $t_{0}=75$. This peak is occurs at a time located mid-way between the two highest crossings that surround it. Over the associated time window the half-distance between each core and the origin, as measured by $\tilde\alpha$, shifts from ${0.5931}$ towards one and then back to ${0.5952}$. 

\begin{figure}[H]
\centering
\begin{minipage}[t][][b]{0.8\columnwidth}
\begin{minipage}[t][][b]{1.1\columnwidth}
\begin{minipage}[t][][b]{0.19\textwidth}
	\centering
	\includegraphics[width=\textwidth,center]{QP_U_75_to_125_2_75.png}
	\subcaption*{$\hat{u}^{(0,50)}_{75,4}$}
	\label{fig:QPDWP_res12_n50_startAtPeak_75_a}
\end{minipage}
\begin{minipage}[t][][b]{0.19\textwidth}
	\centering
	\includegraphics[width=\textwidth,center]{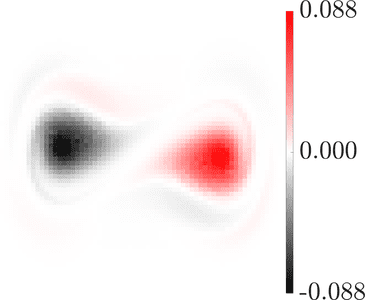}
	\subcaption*{$\hat{u}^{(7,50)}_{75,4}$}
	\label{fig:QPDWP_res12_n50_startAtPeak_75_b}	
\end{minipage}
\begin{minipage}[t][][b]{0.19\textwidth}
	\centering
	\includegraphics[width=\textwidth,center]{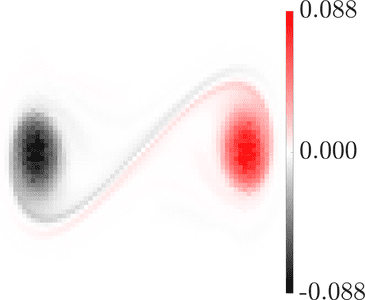}
	\subcaption*{$\hat{u}^{(17,50)}_{75,4}$}
	\label{fig:QPDWP_res12_n50_startAtPeak_75_c}	
\end{minipage}
\begin{minipage}[t][][b]{0.19\textwidth}
	\centering
	\includegraphics[width=\textwidth,center]{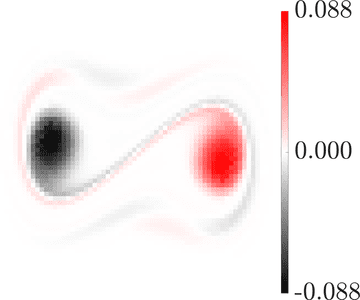}
	\subcaption*{$\hat{u}^{(45,50)}_{75,4}$}
	\label{fig:QPDWP_res12_n50_startAtPeak_75_d}	
\end{minipage}
\end{minipage}
\label{fig:QPDWP_res12_n50_startAtPeak_75}
\hfil
\begin{minipage}[t][][b]{1.1\columnwidth}
\begin{minipage}[t][][b]{0.19\textwidth}
	\centering
	\includegraphics[width=\textwidth,center]{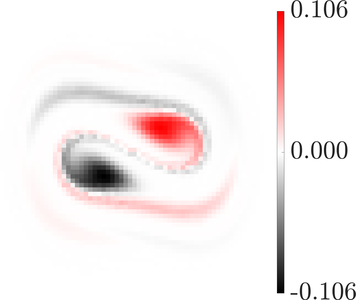}
	\subcaption*{$\hat{u}^{(0,50)}_{274,4}$}
	\label{fig:QPDWP_res12_n50_startAtPeak_274_a}
\end{minipage}
\begin{minipage}[t][][b]{0.19\textwidth}
	\centering
	\includegraphics[width=\textwidth,center]{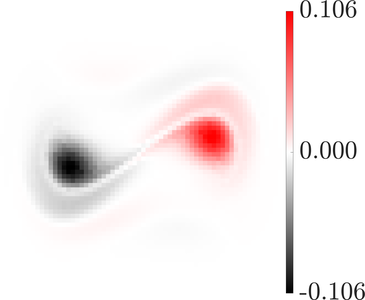}
	\subcaption*{$\hat{u}^{(7,50)}_{274,4}$}
	\label{fig:QPDWP_res12_n50_startAtPeak_274_b}	
\end{minipage}
\begin{minipage}[t][][b]{0.19\textwidth}
	\centering
	\includegraphics[width=\textwidth,center]{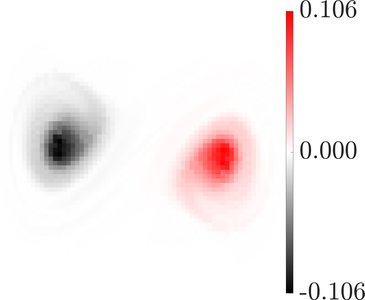}
	\subcaption*{$\hat{u}^{(17,50)}_{274,4}$}
	\label{fig:QPDWP_res12_n50_startAtPeak_274_c}	
\end{minipage}
\begin{minipage}[t][][b]{0.19\textwidth}
	\centering
	\includegraphics[width=\textwidth,center]{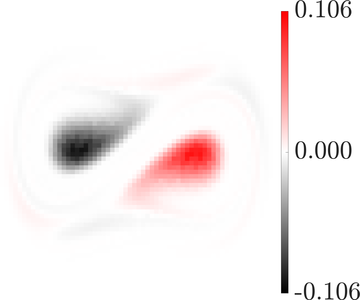}
	\subcaption*{$\hat{u}^{(45,50)}_{274,4}$}
	\label{fig:QPDWP_res12_n50_startAtPeak_274_d}	
\end{minipage}
\end{minipage}\hfil
\end{minipage}\hfil
\begin{minipage}[t][][b]{0.2\columnwidth}
\centering
\hspace{-1.375cm}
\begin{minipage}[t][][b]{1.5\columnwidth}
\vspace{0.7cm}
\includegraphics[width=\columnwidth]{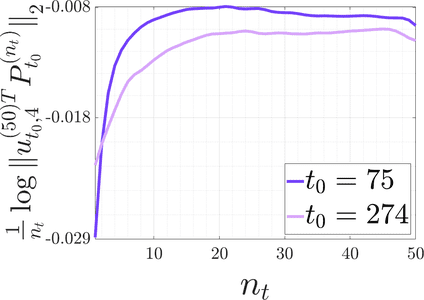}
\subcaption*{\hspace{0.7cm}$\frac{1}{n_{t}} \log{{\|u^{(50)T}_{t_{0},4}\boldsymbol{P}_{t_{0}}^{(n_{t})}\|}_{2}}$}\label{fig:PDWP_anim_centering_log}	
\end{minipage}
\end{minipage}
\caption{Evolved $u^{(50)}_{75,4}$ (top) and $u^{(50)}_{274,4}$ (bottom) of mode $S_{U}^{(4)}$ in Figure~\ref{fig:main_table_QPDWP}, evolved as per Algorithm~\ref{alg:Movies}.}
\label{fig:QPDWP_res12_n50_startAtPeak_274}
\end{figure}

The rounder peak at $t_{0}=276$ reaches $S_{U}^{(4)}={0.5772}$ but remains closer to this maximum for a comparatively longer period. This peak is located to the left of the time window centred between the two highest crossings that surround it. We concentrate on the centred time window $W_{274}^{(50)}$ in order to better characterise the nature of this peak. Over the time window $W_{274}^{(50)}$ the half-distance between each core and the origin shifts from to ${0.4678}$ to ${0.5523}$. Given the previous distance, characterising $W_{75}^{(50)}$, is slightly less than that which characterises $W_{274}^{(50)}$; this rounder peak is associated with a time window where structures experience a larger separation whilst also traversing a greater distance. Appropriate stills from animations of $S_{U}^{(4)}$ using Algorithm~\ref{alg:Movies} for time windows $W_{75}^{(50)}$ and $W_{274}^{(50)}$ are shown in Figure~\ref{fig:QPDWP_res12_n50_startAtPeak_274}. 

\vspace{-0.25cm}
Here we see two cores isolated from background noise. Each core is defined almost exclusively by the positive and negative elements in the singular vectors. Column two of Figure~\ref{fig:QPDWP_res12_n50_startAtPeak_274} shows the left singular vectors after an evolution of $7$ steps. In the case of the sharper peak ($\hat{u}^{(7,50)}_{75,4}$) the division between the two components is better defined than in the case of the flatter, longer lasting peak ($\hat{u}^{(7,50)}_{274,4}$). This occurs despite the fact that the cores associated with the flatter peak experience a greater separation over the full time window.

Column four shows the left singular vectors after an evolution of $45$ steps. In the case of the sharper peak ($\hat{u}^{(45,50)}_{75,4}$) a thin ribbon between the two indicates that mixing has commenced. This is not the case for the flatter, longer lasting peak ($\hat{u}^{(45,50)}_{274,4}$), as the cores are still well separated. Column five shows that the coherency of both of the associated structures initially increases over time. It also indicates that whilst the sharper peak experiences an initially rapid increase in coherency, it is also the first to shift noticeably towards decreasing average coherency, as the time window closes. Taken together, these findings support the notion that the rounder peak is associated with a time window that exhibits comparatively less mixing of the cores over the full time window, even though it concludes with a smaller distance between centres.

\subsection{Splitting of the Southern Polar Vortex} \label{ssec:SPVmodel}
Our final model examines the splitting of the Southern Polar vortex (SPV) in the middle to upper stratosphere in late $2002$. The SPV forms over the austral autumn and breaks apart in the spring. A unique splitting of the SPV in September $2002$ was the first observed major stratospheric warming in the Southern Hemisphere (SH). In $2002$ the stratospheric SPV is understood to have been weakening as early as $21$ September. Splitting is said to have begun by $24$ September. The SPV had separated, at the level of $10$ hPa, by $26$ September \cite{Charltonetal_SPV}. This model was also investigated in \cite{LekienRoss} using finite time Lyapunov exponents and Lagrangian coherent structures. 

\begin{figure}[H]
\centering
\begin{minipage}[b]{\textwidth}
\centering
\begin{minipage}[b]{0.3\textwidth}
\centering
\includegraphics[width=\columnwidth]{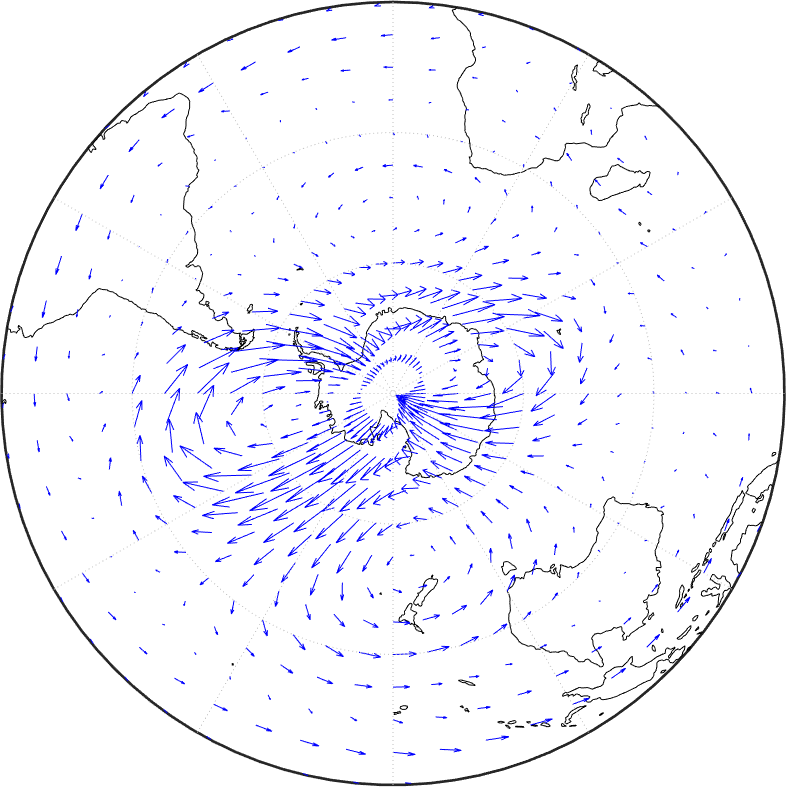}
\subcaption{$0000$ on $7$ Sep.}
\end{minipage}
\hfil
\begin{minipage}[b]{0.3\textwidth}
\centering
\includegraphics[width=\columnwidth]{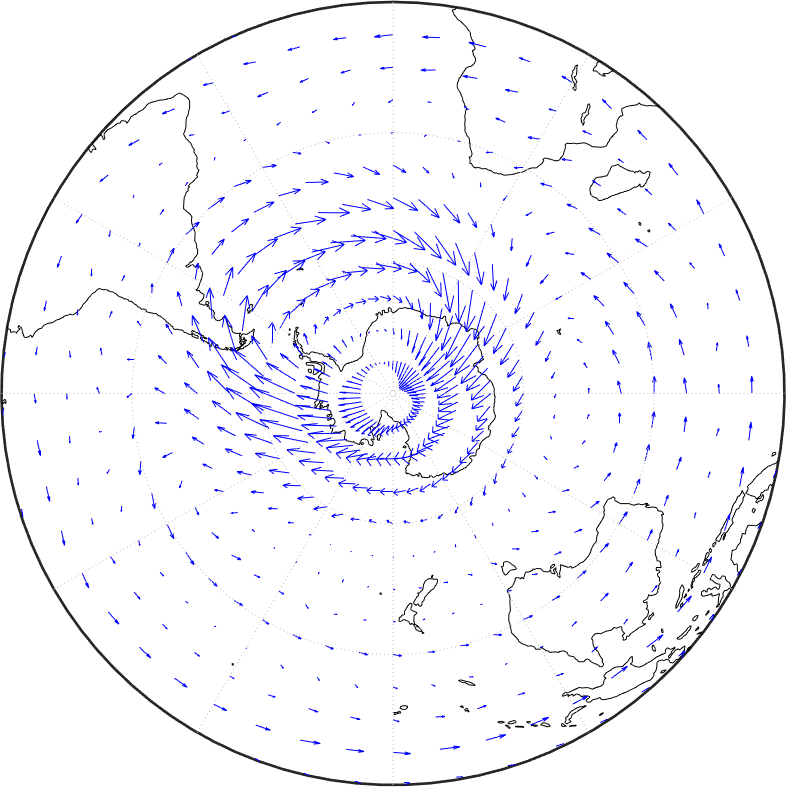}
\subcaption{$0600$ on $16$ Sep.}
\end{minipage}
\hfil
\begin{minipage}[b]{0.3\textwidth}
\centering
\includegraphics[width=\columnwidth]{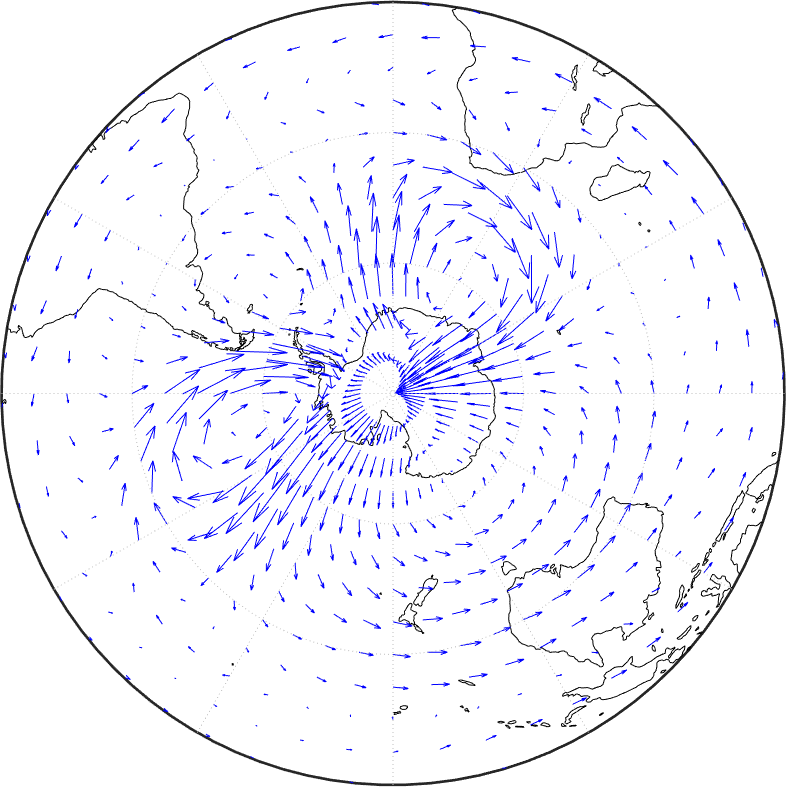}
\subcaption{$1200$ on $25$ Sep.}
\end{minipage}
\end{minipage}
\caption{Southern hemisphere wind speed (easterly and northerly) on the $850$ K isentropic surface.}
\label{fig:SPV_VF}
\end{figure}

To examine the splitting we employ ECMWF Re-Analysis data on the isentropic surface defined by a potential temperature of $850$ K (near $10$ hPa) for a large portion of the SH \cite{DeeD2011TErc}. In this case the vector field is defined by the horizontal speed of air moving to the east and north. The ECMWF provides this data at $6$ hourly intervals ($0000$, $0600$, $1200$ and $1800$ UTC) in the temporal direction at a spatial resolution of $0.75^{\circ}$. This vector field is characterised by up to $480 \times 121$ data points in the longitude and latitude directions (for the full SH) and up to $368$ in the temporal direction (August through October 2002). Figure~\ref{fig:SPV_VF} indicates this wind speed at three times of interest for the full SH. One notes that times mentioned in this context all refer to Coordinated Universal Time (UTC).
\vspace{-0.5cm}
For this model we employ two methods. We either consider $X=S^{1}\times[-90^\circ,-30^{\circ}]$ where $S^{1}$ is the circle parameterised from $0^{\circ}$ to $360^{\circ}$ or we consider the full SH but only seed bins with centres at a latitude lower than $50^{\circ}$S such that $X_{t_{0}}=S^{1}\times[-90^{\circ},-50.27^{\circ}]$. New bins are included over time, as they are occupied by the advected particles. This approach is similar to that of \cite{FSM_2010}. To examine the splitting of the SPV we choose time windows of length $n=56$. This is a similar length to that employed in Sections~\ref{ssec:DWPmodel} and~\ref{ssec:QPmodel}, and corresponds to two weeks. Over periods of this length, isentropic surfaces generally do not experience significant changes \cite{BernardBinson2002,FSM_2010}. 

Let us begin by employing Algorithm~\ref{alg:RWinds} to construct rolling windows for Ulam matrices of size $2^{14} \times 2^{14}$. In this case, we set $X=S^{1}\times[-90^{\circ},-30^{\circ}]$. This allows us to cover the Southern hemisphere in $2^{14}$ bins south of $30^{\circ}$S. This is done to capture the dynamics of both the SPV and any daughter vortices that might be generated.

\begin{figure}[H]
\centering
\includegraphics[width=\columnwidth,center]{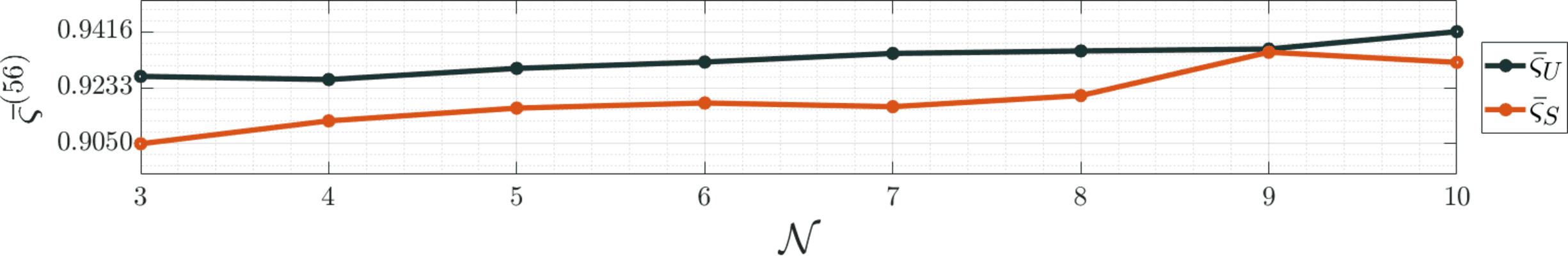}
\caption{Mean equivariance mismatch, as per Algorithm~\ref{alg:Eqvr}, for the leading $3$ of $\mathcal{N}$ modes using the two pairing methods given in Algorithms~\ref{alg:SValsDijk} and~\ref{alg:L2VecDist} with $n=56$ and $t_0 \in[0000 \: 1 \: \text{August}, 1800 \: 30 \: \text{September}]$. Here the Ulam matrices, describing transitions for the area south of $30^{\circ}$S, are of dimension $m \times m$ for $m=2^{14}$.} \label{fig:Eq_SPV_res14_n55_trackingMethods}
\end{figure}

\begin{figure}[H]
\begin{minipage}[b]{\textwidth}
\centering
\begin{minipage}[b]{\textwidth}
\centering
\includegraphics[width=1\textwidth,center]{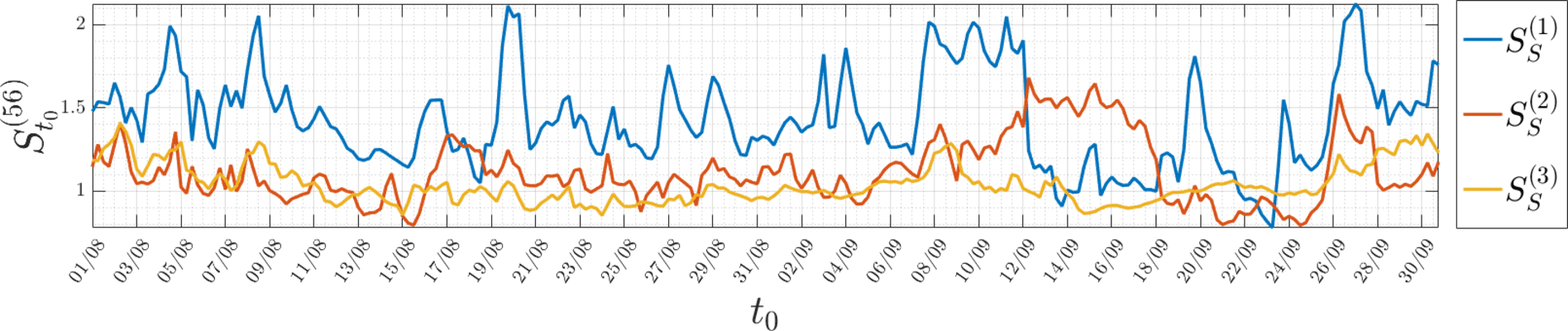}
\subcaption{Rolling windows for compositions of square matrices considering the area south of $30^{\circ}$S.}
\label{fig:SPV_res14_n56_sqMat}
\end{minipage}
\begin{minipage}[b]{\textwidth}
\centering
\includegraphics[width=1\textwidth,center]{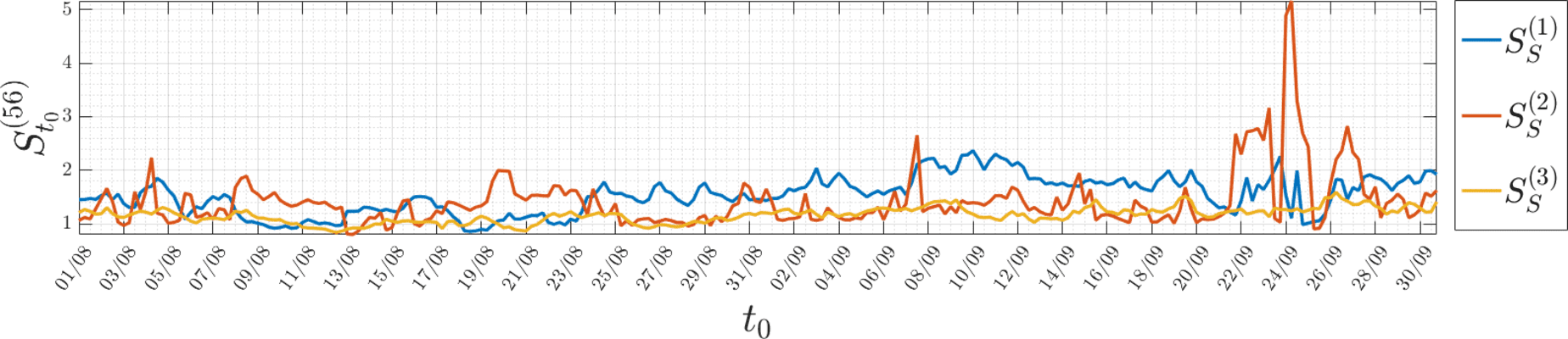}
\subcaption{Rolling windows for the SH where only the area south of $50.27^{\circ}$S is initially seeded for each $t_{0}$.}
\label{fig:SPV_res14_n56_rectMat}
\end{minipage}
\end{minipage}
\caption{Leading $3$ of $\mathcal{N}=3$ tracked paths of singular values of rolling windows paired using Algorithm~\ref{alg:SValsDijk} for $n=56$.}
\label{fig:SPV_res14_n56_trackings}
\end{figure}

We use this initial setting to choose an appropriate number of singular vectors to consider. We then examine the average value of equivariance mismatch for our two pairing methods, with $t_{0}$ set to consider all initial times in August and September $2002$. Figure~\ref{fig:Eq_SPV_res14_n55_trackingMethods} presents results that look to identify an effective pairing through time for the leading $3$ of $\mathcal{N}$ modes. 
\vfil \hfil 
\vspace{-0.25cm}
In both cases the lowest equivariance occurs when $\mathcal{N}=3$. As a greater number of modes are considered $\bar{\varsigma}_{U}$ and $\bar{\varsigma}_{S}$ tend to increase. The lowest value of average equivariance mismatch for all pairings is attained by $\bar{\varsigma}_{S}$ at $\mathcal{N}=3$. In light of these results we set $\mathcal{N}=3$ and use Algorithm~\ref{alg:SValsDijk} to examine the pairing of modes through time. Results for these parameters are presented in Figure~\ref{fig:SPV_res14_n56_sqMat}. 

As in \cite{FSM_2010}, it makes sense to also focus on a smaller region where the vortex splitting occurs.
In this setting, the evolution of mass is tracked through time conditional on an initial seeding. Here we consider the full southern hemisphere but at each instance of $t_{0}$ we only seed bins with centres at a latitude lower than $50^{\circ}$S such that $X_{t_{0}}=S^{1}\times[-90^{\circ},-50.27^{\circ}]$ initialises each rolling window. This area is seeded because it is known that the stratospheric polar night jet develops at latitudes of about $60^{\circ}$S during the austral winter. 

Given that our measure of equivariance is not applicable in this setting,
we choose to maintain $\mathcal{N}=3$ for comparability. The Ulam matrices are constructed by seeding $m=7,296$ bins at each initial time $t_{0}$. All time windows end with no more than $m'=11,776$ bins. Results for this case are presented in Figure~\ref{fig:SPV_res14_n56_rectMat}.  
\vspace{-0.25cm}
The most notable characteristic for either of the results presented in Figure~\ref{fig:SPV_res14_n56_trackings} is the striking peak on $24$ September in Figure~\ref{fig:SPV_res14_n56_rectMat}. This indicates that something is fundamentally different between the time windows ending before this peak arises and those starting once the peak has dissipated. Figure~\ref{fig:SPV_res14_n56_sqMat} offers a much more complex illustration of the dynamics. In this case what is being indicated regarding the occurrence of fundamental changes in the dynamics is much less clear.

The definitive peak in Figure~\ref{fig:SPV_res14_n56_rectMat} occurs at $0600$ on $24$ September. There is also an an earlier anomaly that begins to move upwards just prior to $22$ September. These peaks occur as the SPV is elongated and separates. Also of interest is the smaller sharp peak that begins around $6$ September and reaches a maximum on $7$ September. This peak follows an earlier occurrence that served to elongate the SPV. This event did not result in a splitting of the vortex \cite{Charltonetal_SPV}. 

\begin{figure}[H] 
\begin{minipage}[b]{0.245\textwidth}
\centering 
\includegraphics[width=\textwidth,center]{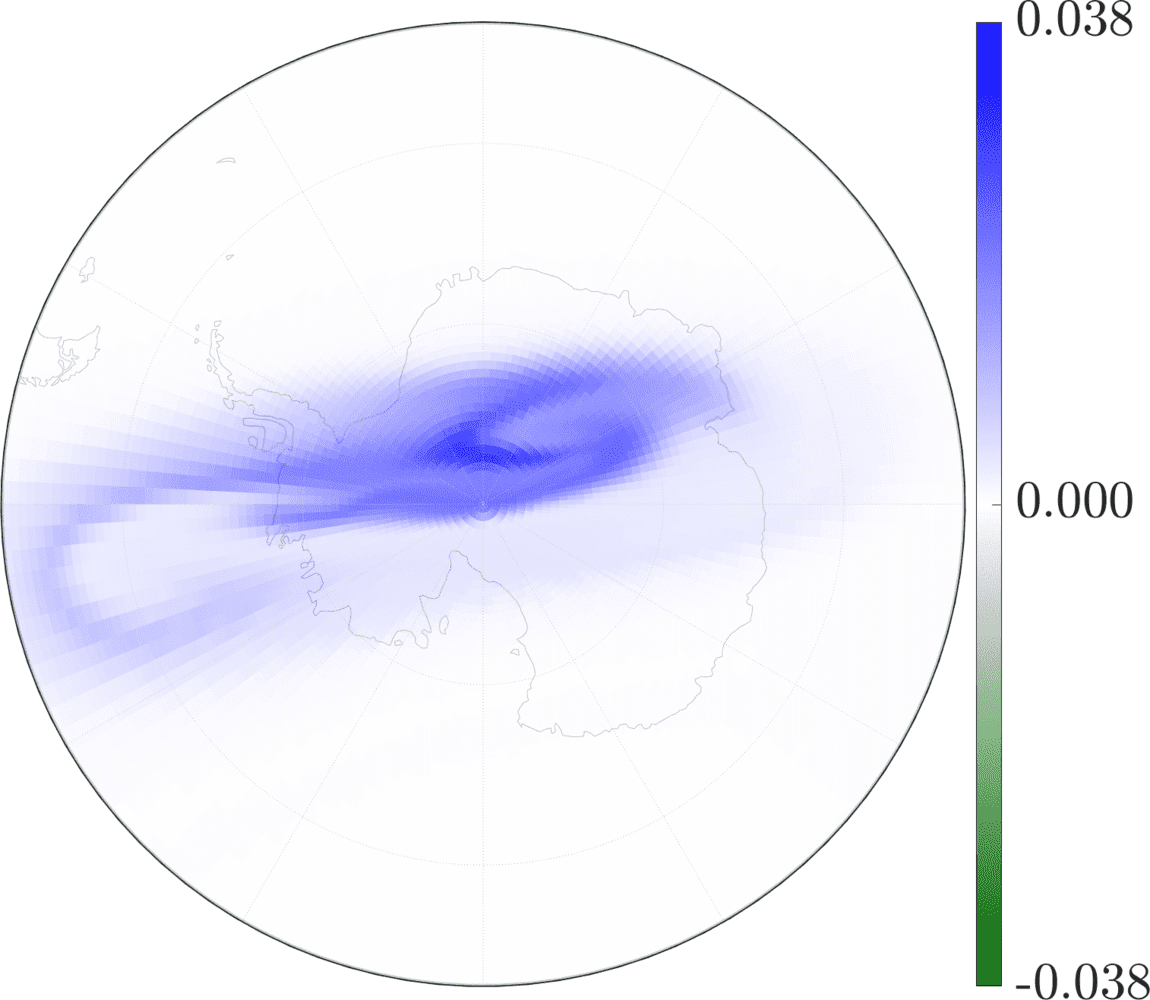}
\subcaption{$0000$ on $7$ Sep.}\label{fig:jetStreaks_b2}
\end{minipage}
\begin{minipage}[b]{0.245\textwidth}
\centering 
\includegraphics[width=\textwidth,center]{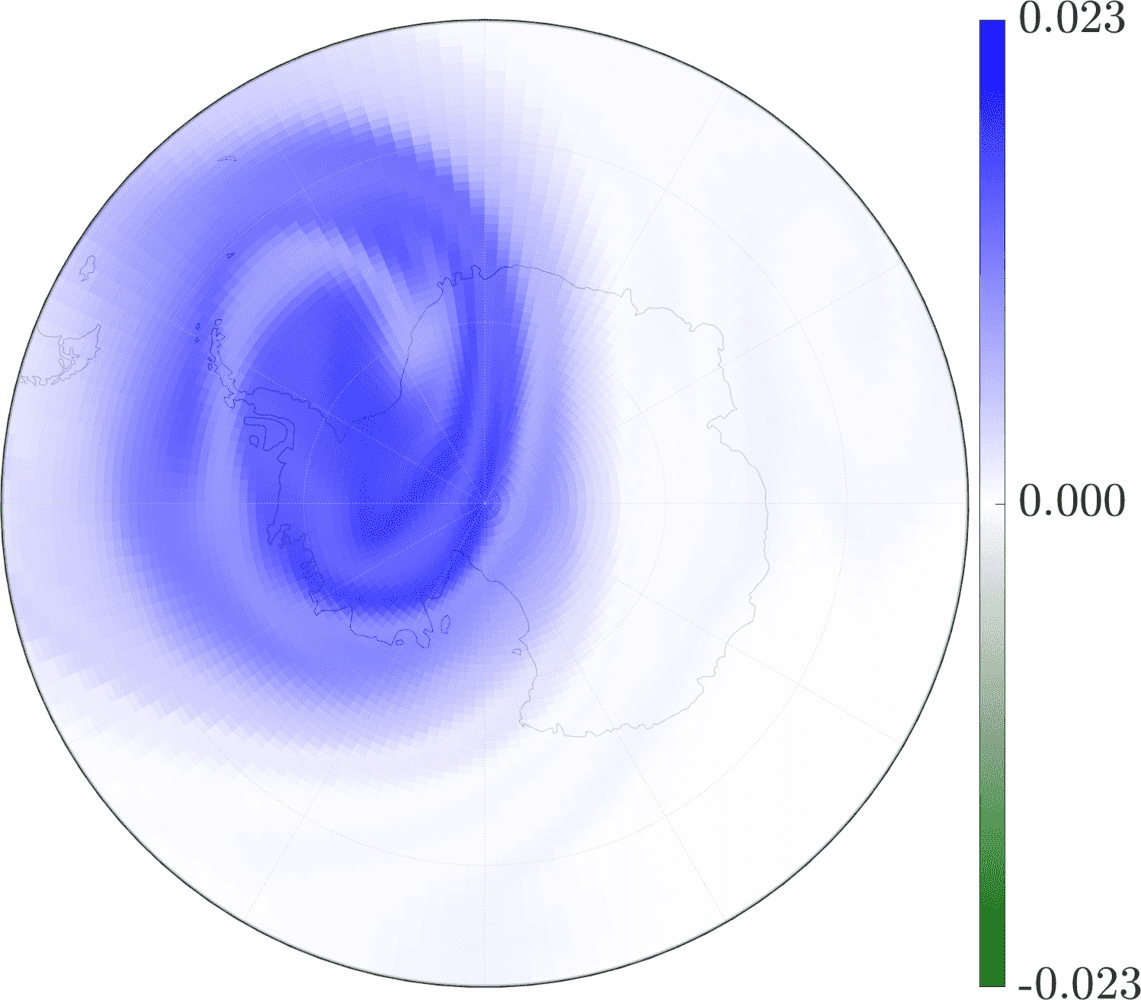}
\subcaption{$0600$ on $13$ Sep.}\label{fig:jetStreaks_a}
\end{minipage}
\begin{minipage}[b]{0.245\textwidth}
\centering
\includegraphics[width=\textwidth,center]{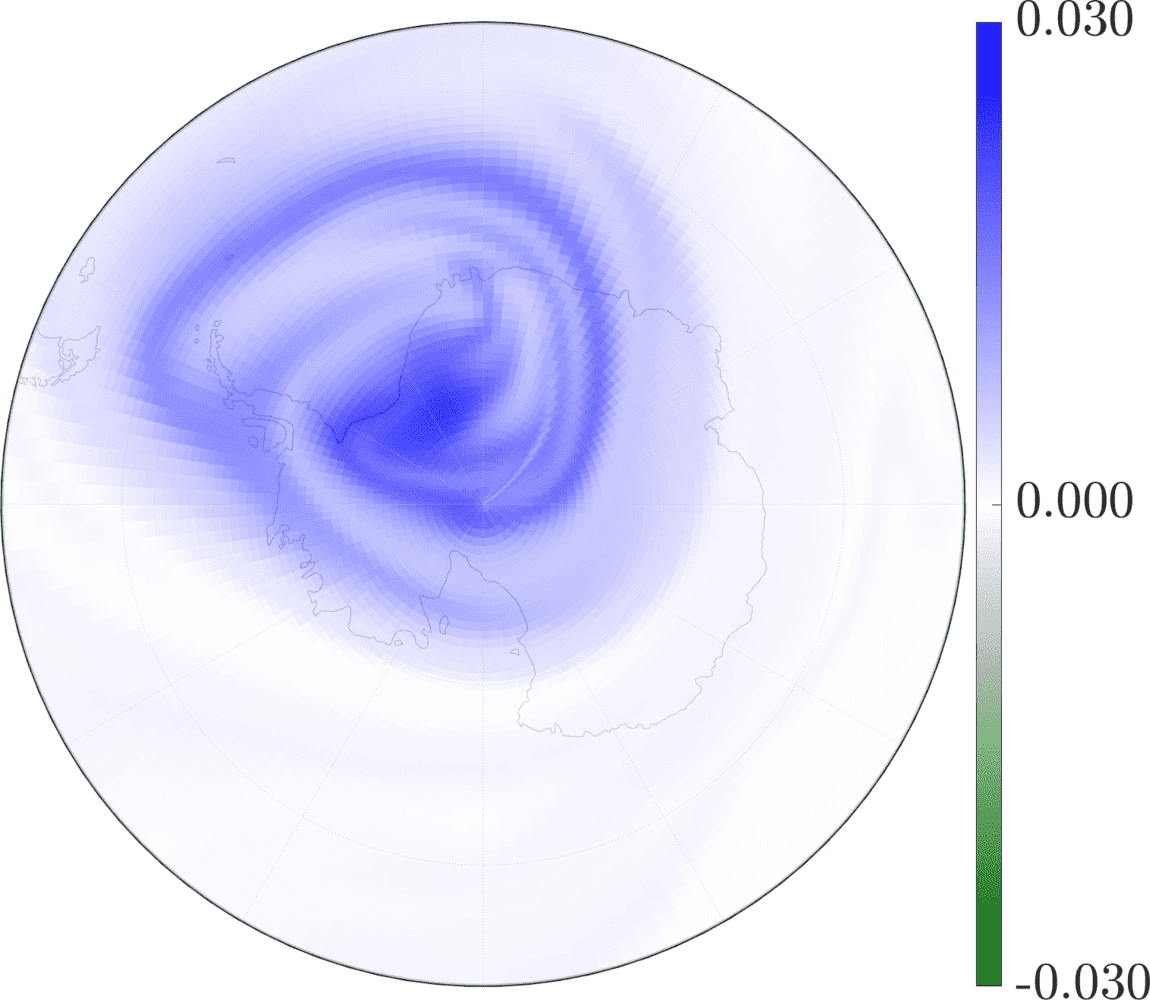}
\subcaption{$0600$ on $16$ Sep.}\label{fig:jetStreaks_a2}
\end{minipage}
\begin{minipage}[b]{0.245\textwidth}
\centering 
\includegraphics[width=\textwidth,center]{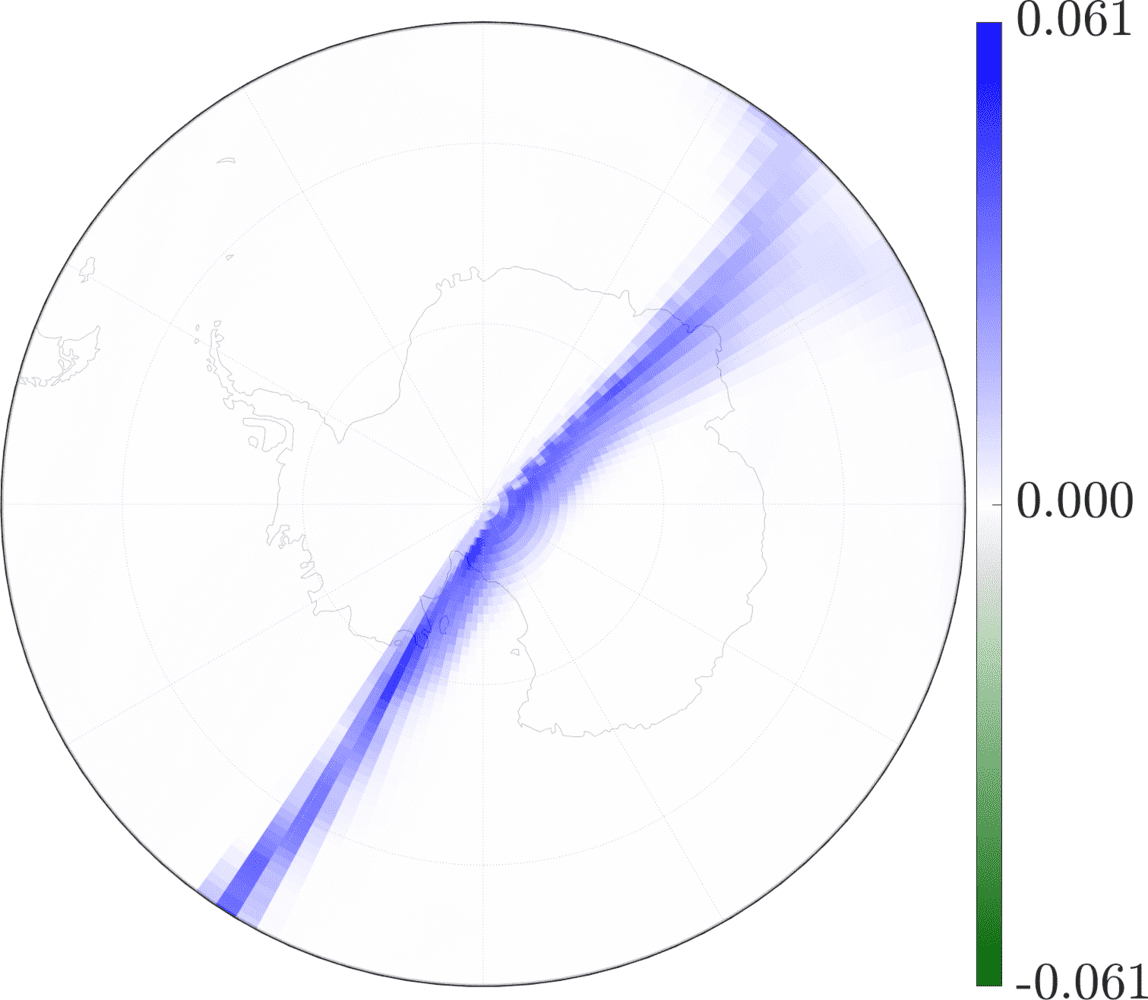}
\subcaption{$0600$ on $22$ Sep.}\label{fig:jetStreaks_b}
\end{minipage}
\begin{minipage}[b]{0.245\textwidth}			
\centering
\includegraphics[width=\textwidth,center]{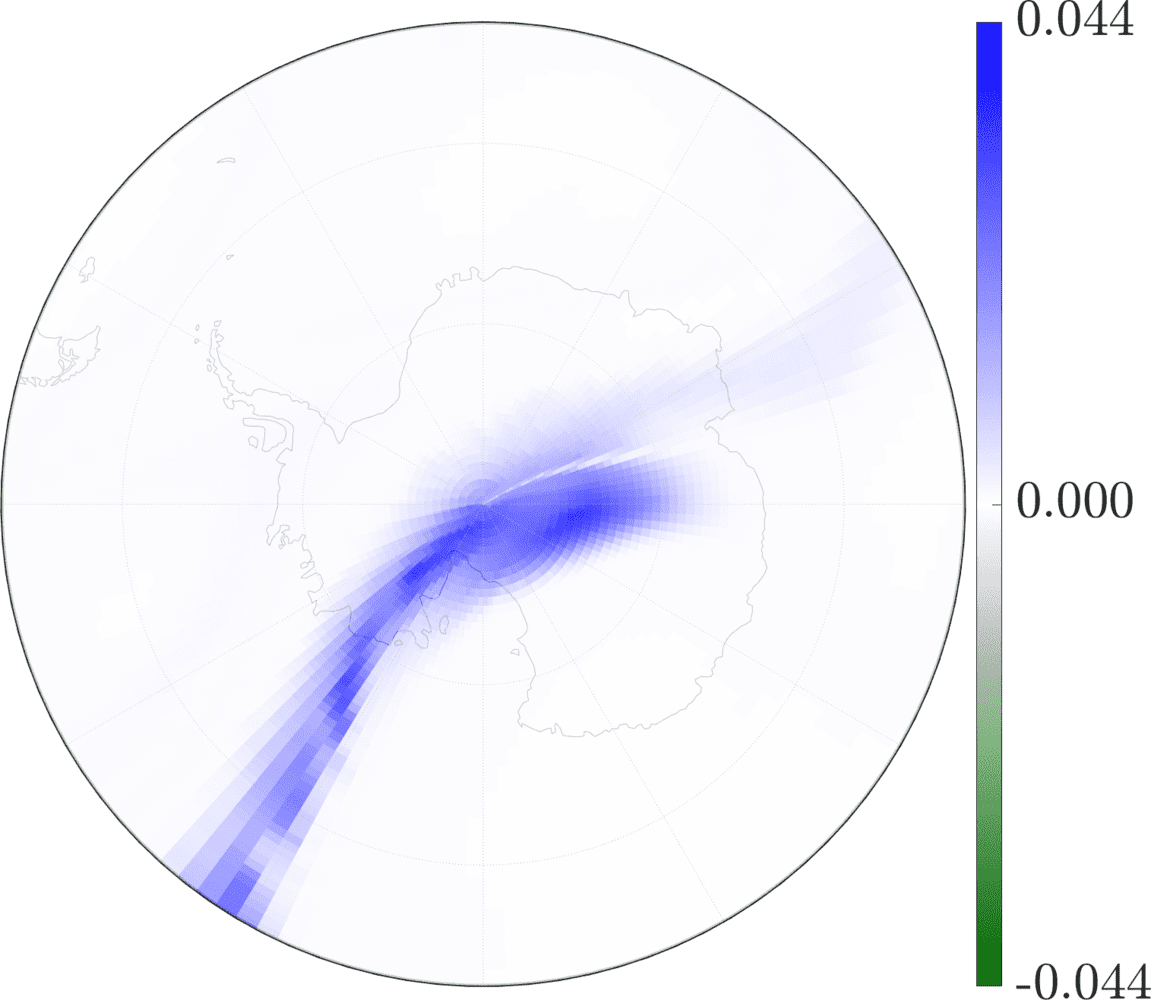}
\subcaption{$1800$ on $23$ Sep.}\label{fig:jetStreaks_c}
\end{minipage}
\begin{minipage}[b]{0.245\textwidth}			
\centering
\includegraphics[width=\textwidth,center]{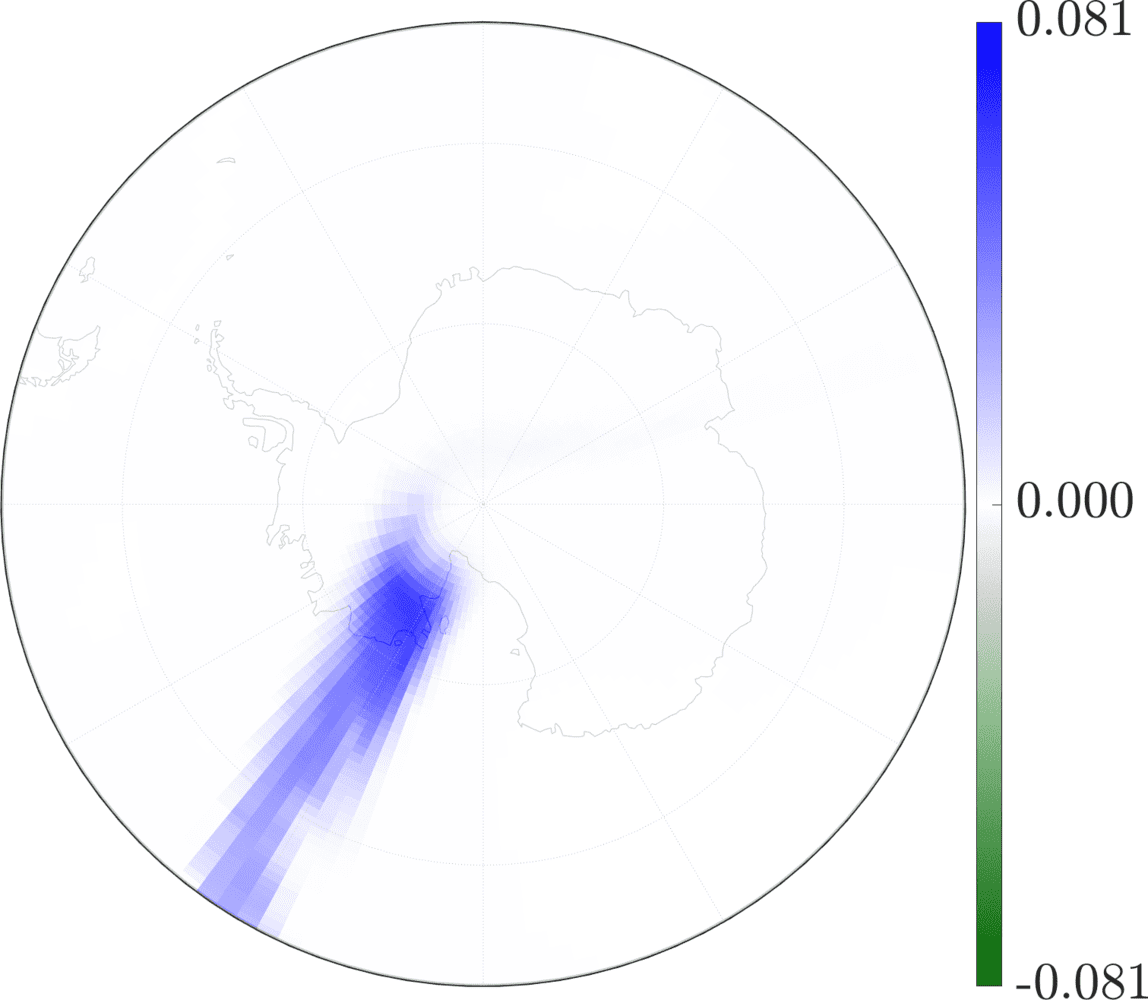}
\subcaption{$0600$ on $24$ Sep.}\label{fig:jetStreaks_c2}
\end{minipage}
\begin{minipage}[b]{0.245\textwidth}
\centering
\includegraphics[width=\textwidth,center]{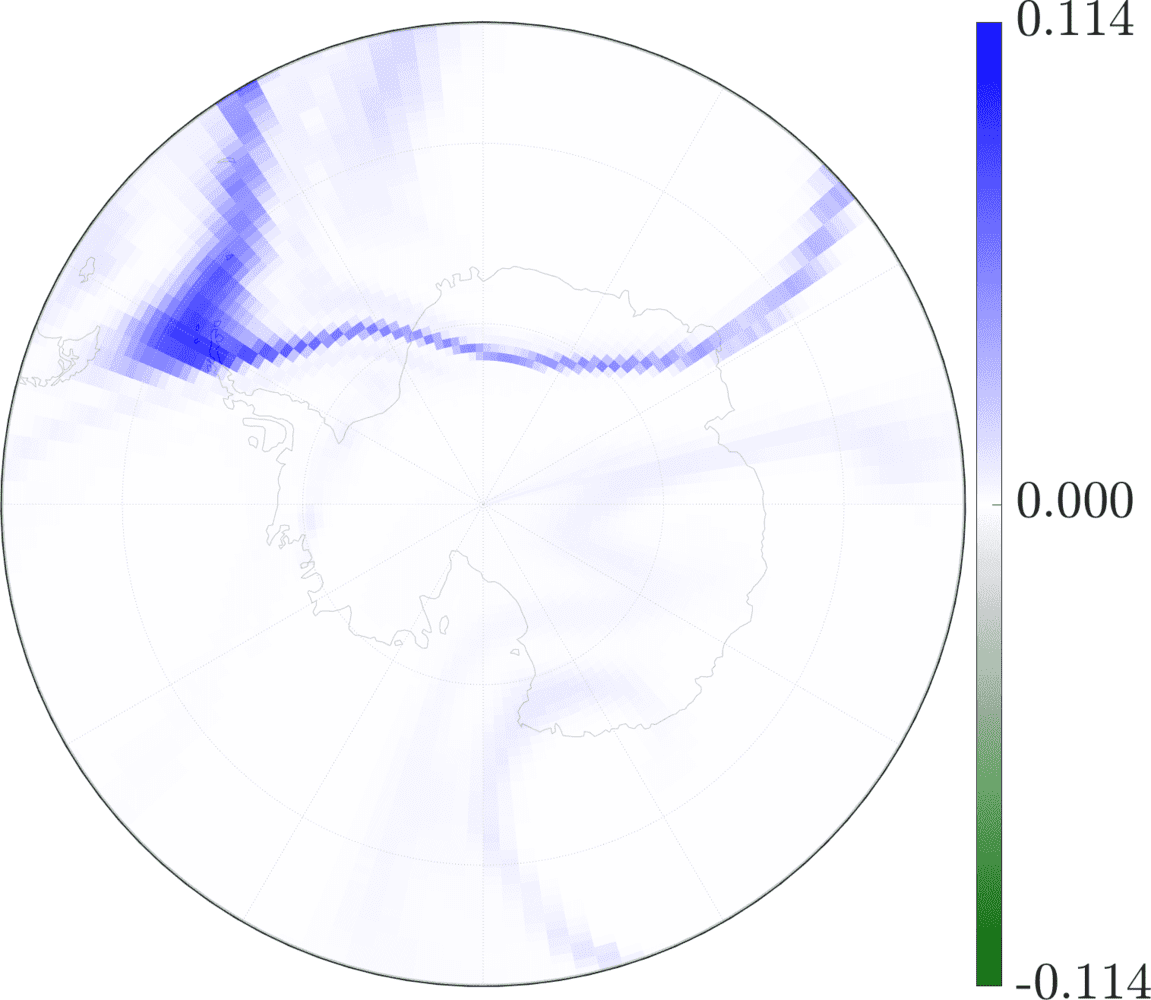}
\subcaption{$1200$ on $24$ Sep.}\label{fig:jetStreaks_d}
\end{minipage}
\begin{minipage}[b]{0.245\textwidth}
\centering
\includegraphics[width=\textwidth,center]{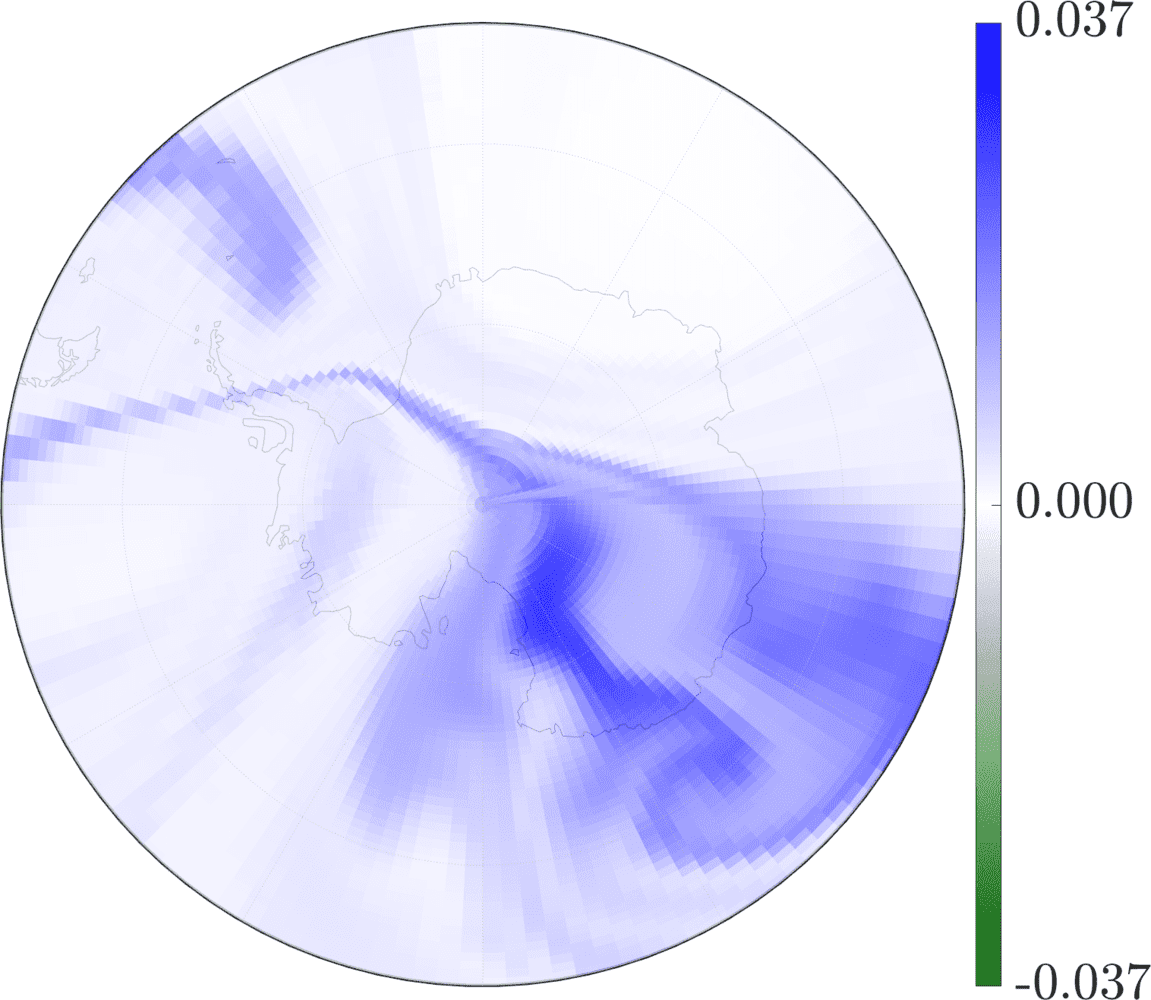}
\subcaption{$1800$ on $24$ Sep.}\label{fig:jetStreaks_d2}
\end{minipage}
\caption{Leading singular vectors, for various $t_{0}$, of matrix compositions associated with Figure~\ref{fig:SPV_res14_n56_rectMat} where time windows are of length $n=56$. The area illustrated is south of $50^{\circ}$S and the time given in the label is the relevant $t_{0}$ for that window.}
\label{fig:jetStreaks}
\end{figure}

A collection of initial time singular vectors for various time windows are illustrated in Figure~\ref{fig:jetStreaks}. The earlier elongation from which the SPV recovered is shown in Figure~\ref{fig:jetStreaks_b2}. Figures~\ref{fig:jetStreaks_a} and~\ref{fig:jetStreaks_a2} illustrate the fact that the polar vortex exists as a mass of cold air contained by thin, rapid flowing streams of air. The strongest of these streams serve to isolate the SPV from warmer surroundings and make up the polar jet. 

\vspace{-0.25cm}
By $22$ September the SPV was in an elongated state with two anticyclones to either side. A weaker anticyclone was circulating over the tip of South America whilst a stronger, quasi-stationary one had developed between Australia and Antarctica \cite{Charltonetal_SPV,oNeill2017vortex}. The ribbon of mass that divides the space in Figure~\ref{fig:jetStreaks_b} into two distinct halves appears to coincide with jet streams dividing the SPV from the area associated with the stronger anticyclone. 

Because the polar vortex is surrounded by the polar jet stream, any breakdown of the polar jet stream is directly related to the behaviour of the polar vortex through time. Figures~\ref{fig:jetStreaks_c},~\ref{fig:jetStreaks_c2} and~\ref{fig:jetStreaks_d} detect one such breakdown. Figure~\ref{fig:jetStreaks_c} captures a "buckling" of the polar jet stream as it weakens. The SPV subsequently separates in two and the stronger anticyclone extends between the SPV to connect with the weaker one near the tip of South America. This can be seen in Figures~\ref{fig:jetStreaks_d} and~\ref{fig:jetStreaks_d2}.

Given the fundamental differences in what is observed in time windows that begin at either side of the peak at $0600$ $24$ September, we explore time windows centred at this peak. For singular vectors corresponding to the paired modes illustrated in Figure~\ref{fig:SPV_res14_n56_rectMat} we examine windows that cover the period $[0600$ $17$ September, $0600$ $1$ October$]$. 
 
\begin{figure}[H]
\begin{minipage}[b]{0.245\textwidth}
\centering
\includegraphics[width=\textwidth,center]{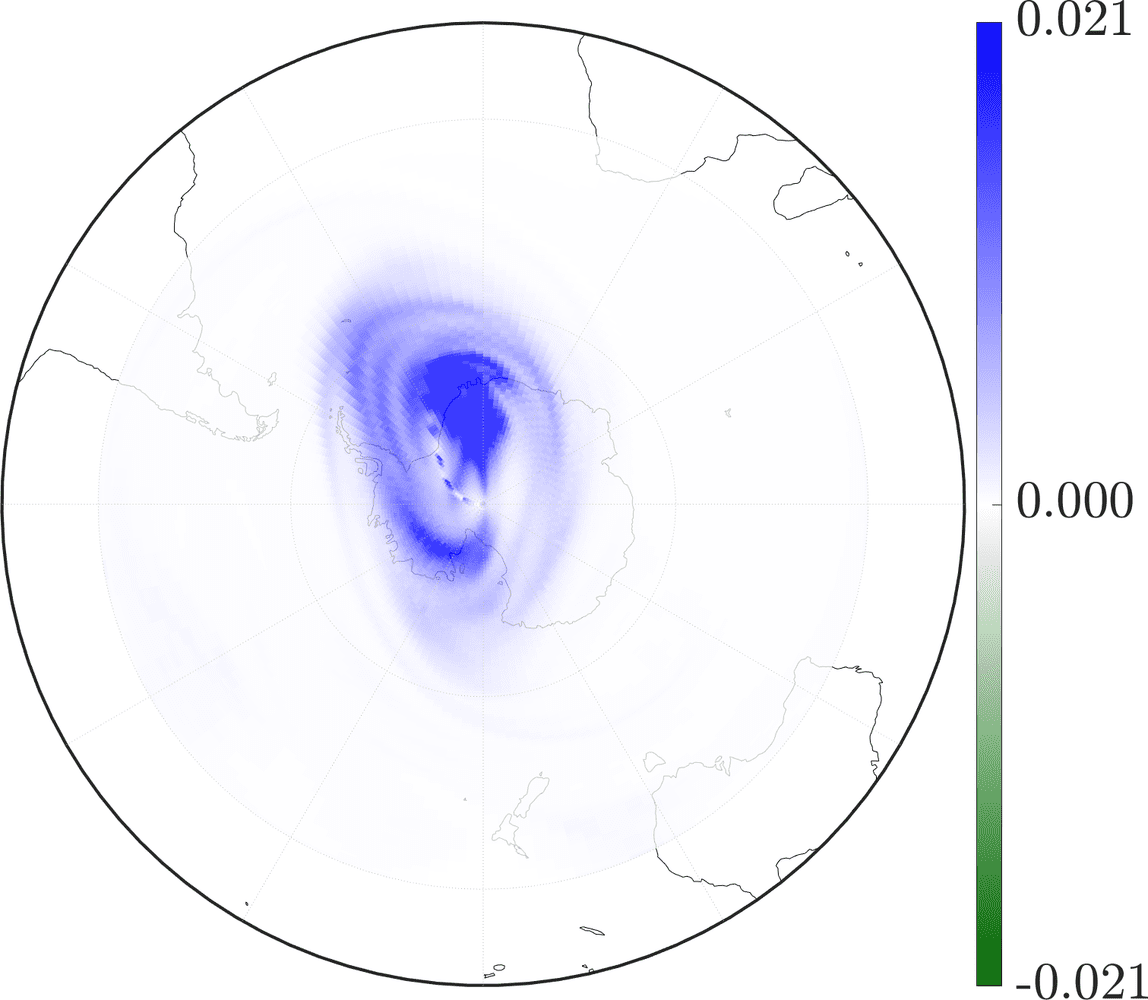}
\subcaption{$1800$ on $19$ Sep.}
\end{minipage}
\begin{minipage}[b]{0.245\textwidth}
\centering
\includegraphics[width=\textwidth,center]{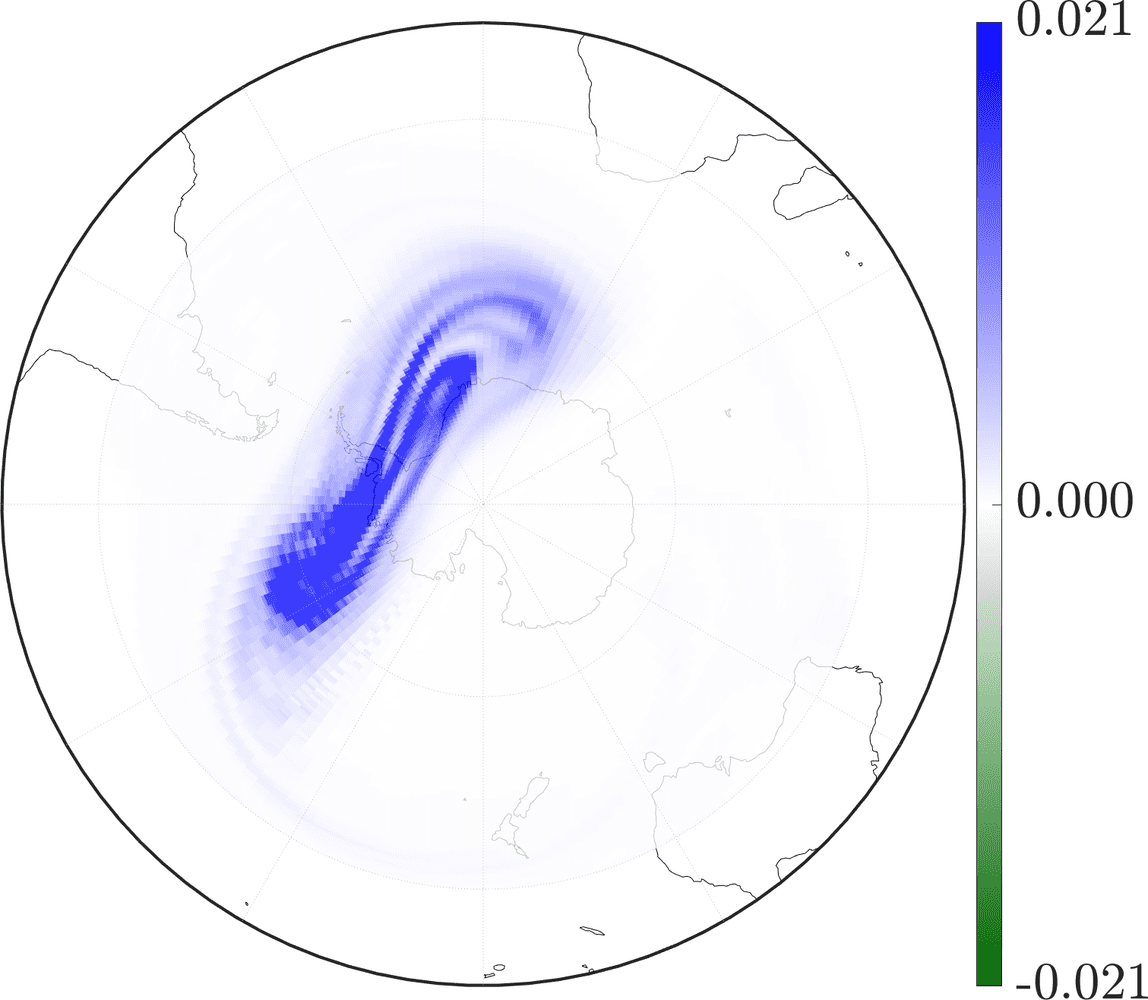}
\subcaption{$1200$ on $22$ Sep.}
\end{minipage}
\begin{minipage}[b]{0.245\textwidth}
\centering
\includegraphics[width=\textwidth,center]{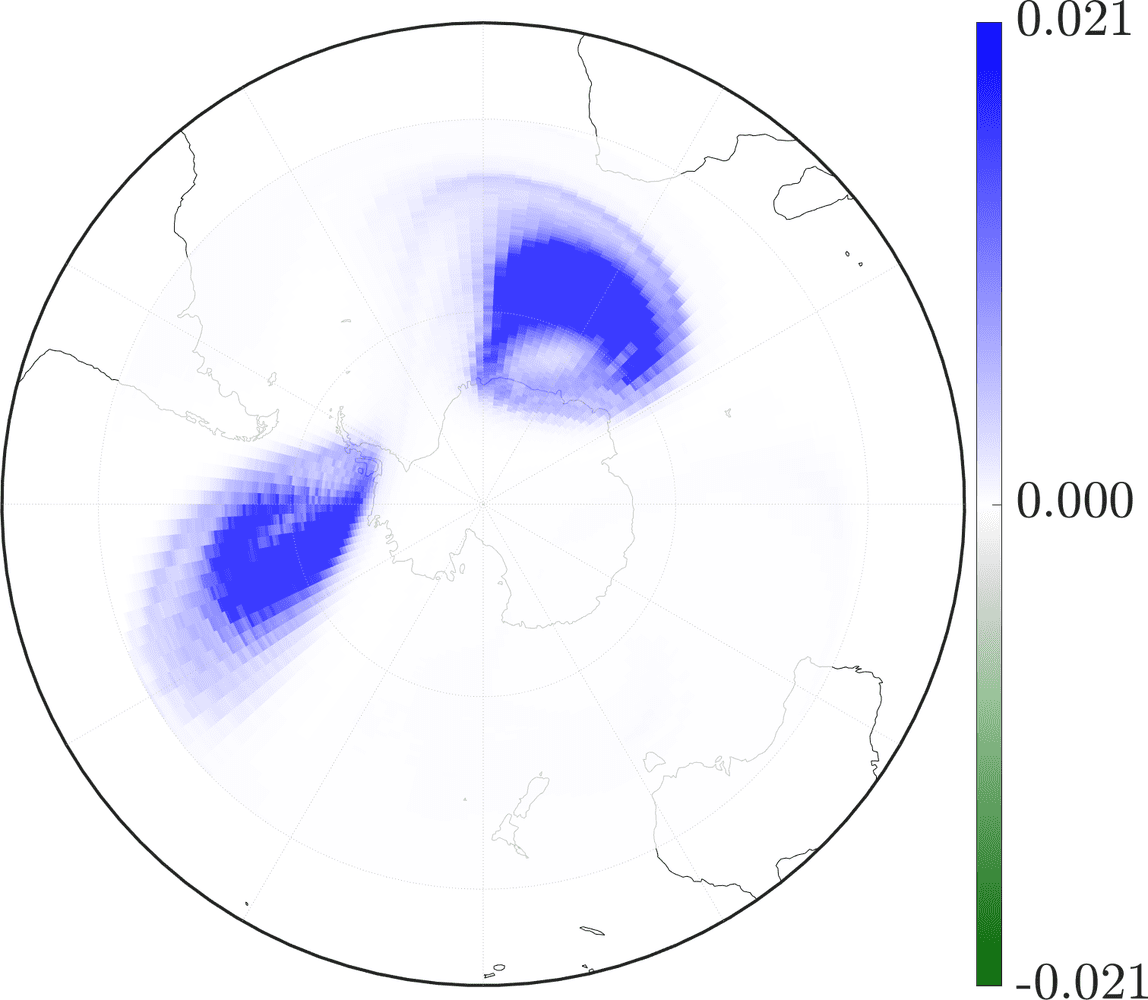}
\subcaption{$1200$ on $25$ Sep.}
\end{minipage}
\begin{minipage}[b]{0.245\textwidth}			
\centering
\includegraphics[width=\textwidth,center]{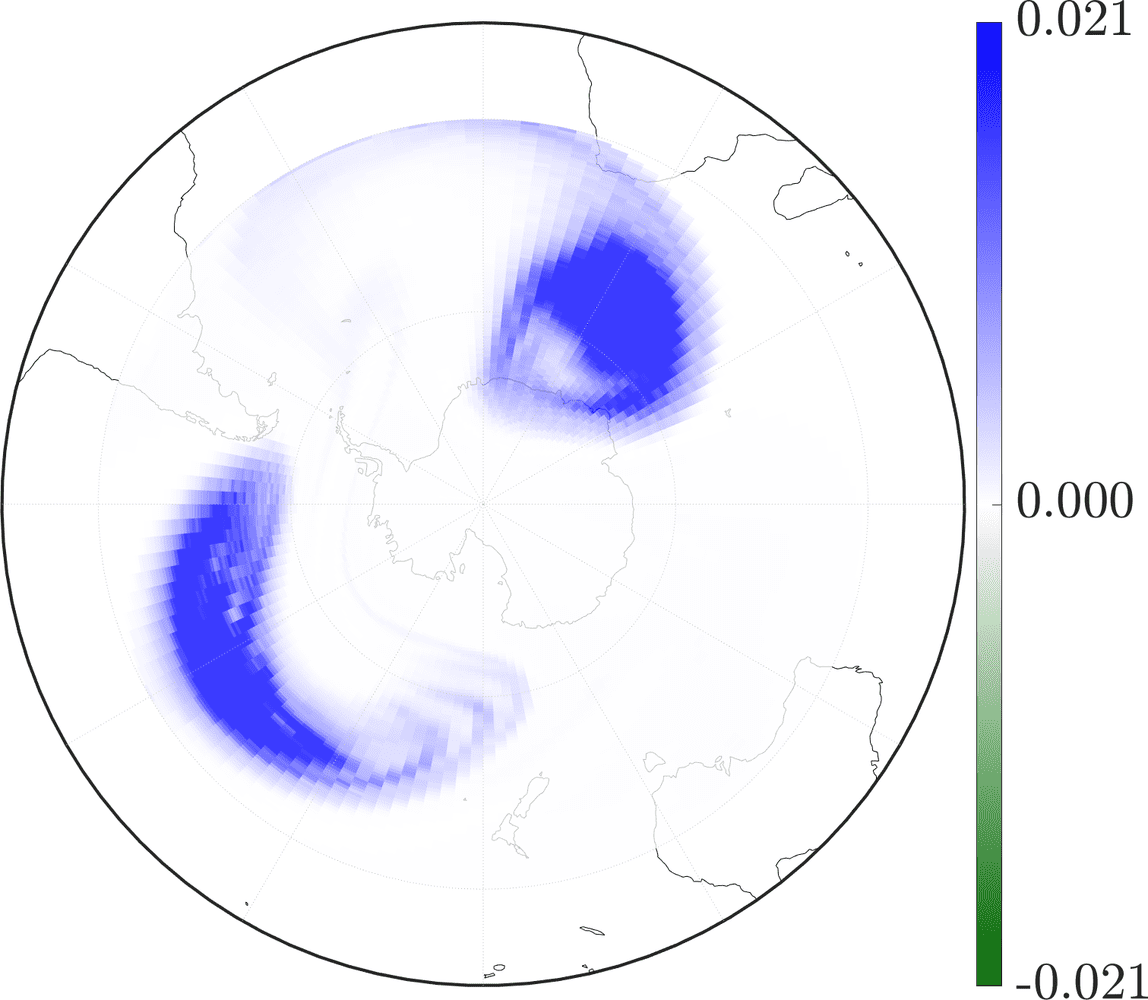}
\subcaption{$1200$ on $27$ Sep.}
\end{minipage}
\caption{Evolved leading mode associated with Figure~\ref{fig:SPV_res14_n56_sqMat} for a time window centred on the peak at $1800$ on $23$ Sep. This is illustrated on the area south of $15^{\circ}$S.}
\label{fig:res14Sq_south30}
\end{figure}

\begin{figure}[H]
\begin{minipage}[b]{0.245\textwidth}
\centering
\includegraphics[width=\textwidth,center]{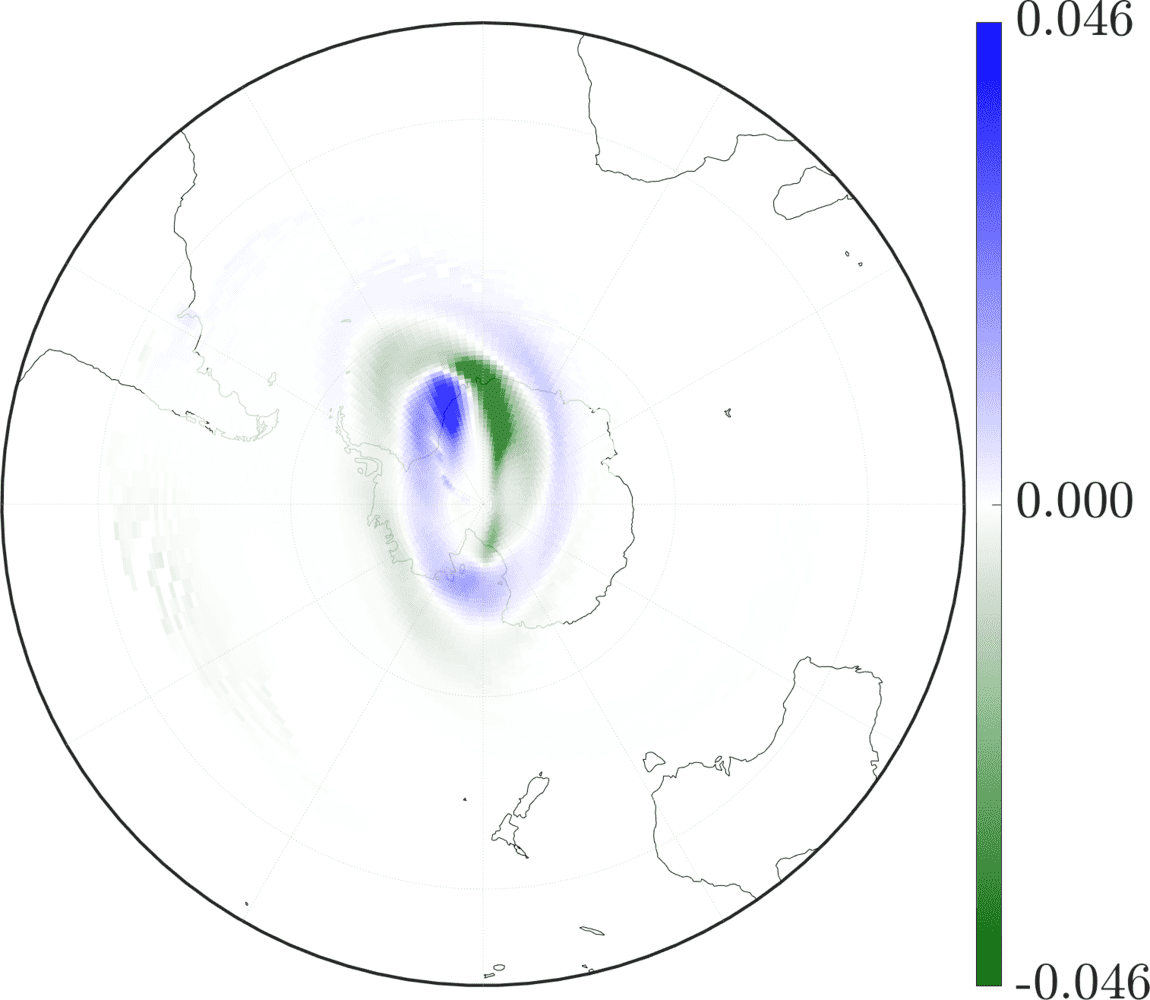}
\subcaption{$1800$ on $19$ Sep.}
\end{minipage}
\begin{minipage}[b]{0.245\textwidth}
\centering
\includegraphics[width=\textwidth,center]{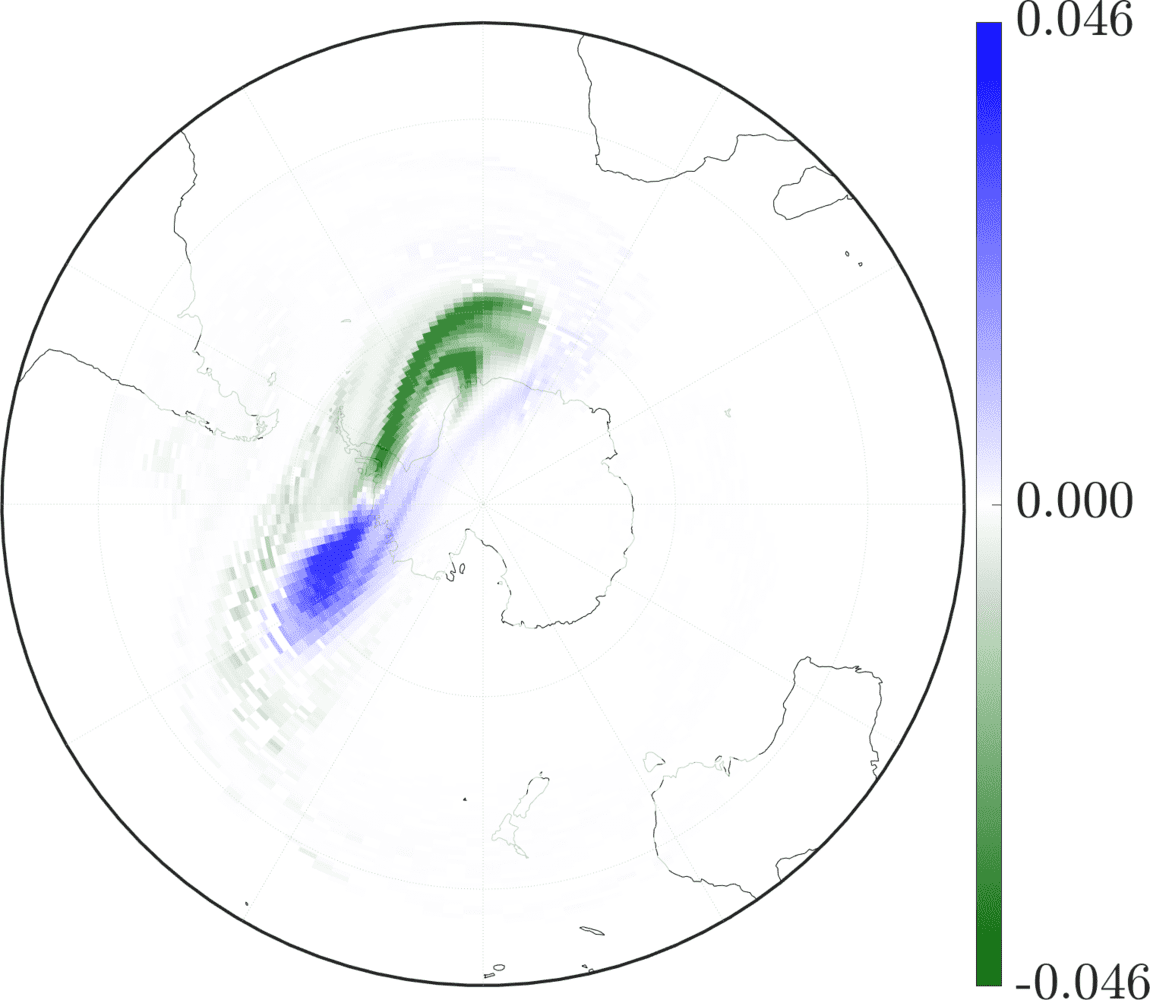}
\subcaption{$1200$ on $22$ Sep.}
\end{minipage}
\begin{minipage}[b]{0.245\textwidth}
\centering
\includegraphics[width=\textwidth,center]{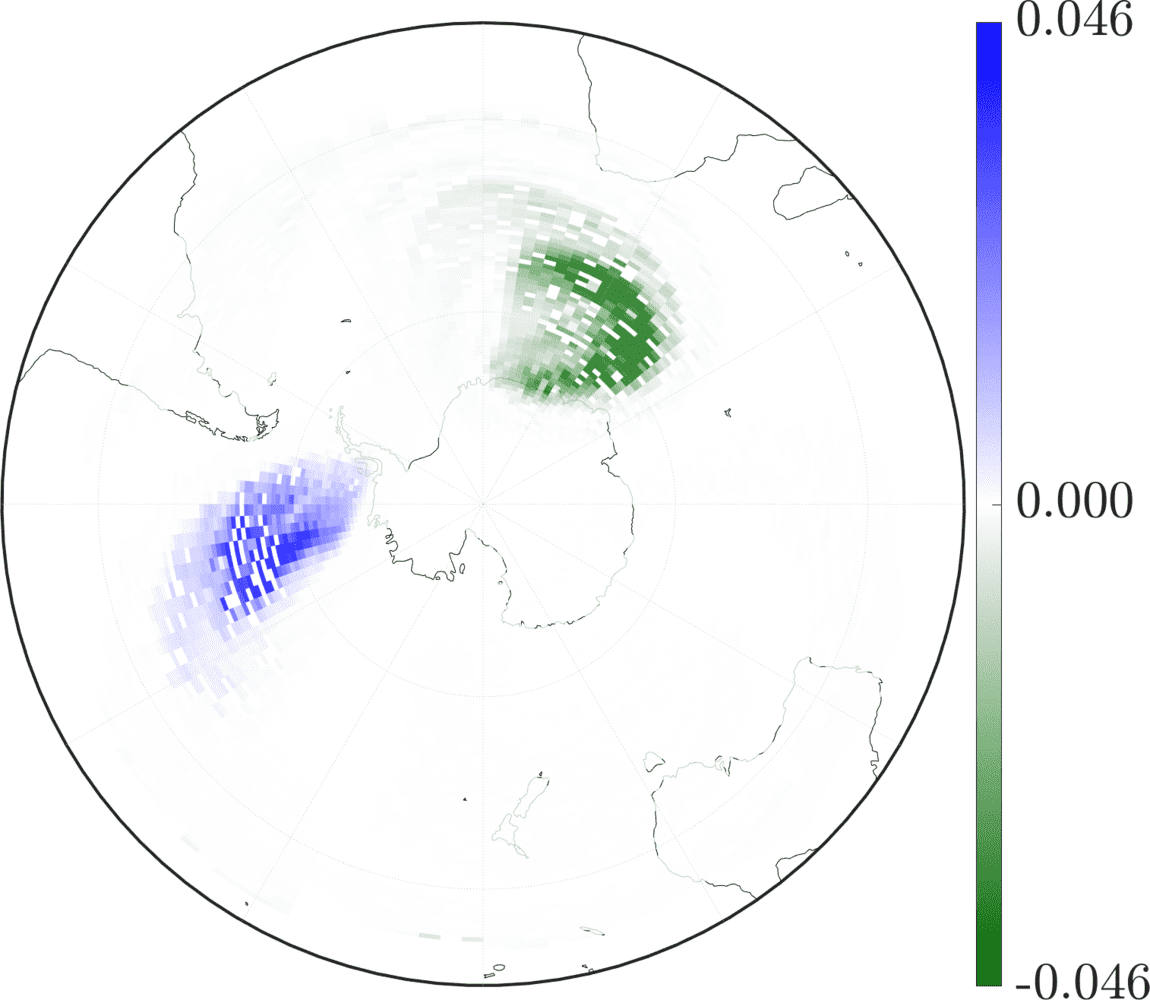}
\subcaption{$1200$ on $25$ Sep.}
\end{minipage}
\begin{minipage}[b]{0.245\textwidth}			
\centering
\includegraphics[width=\textwidth,center]{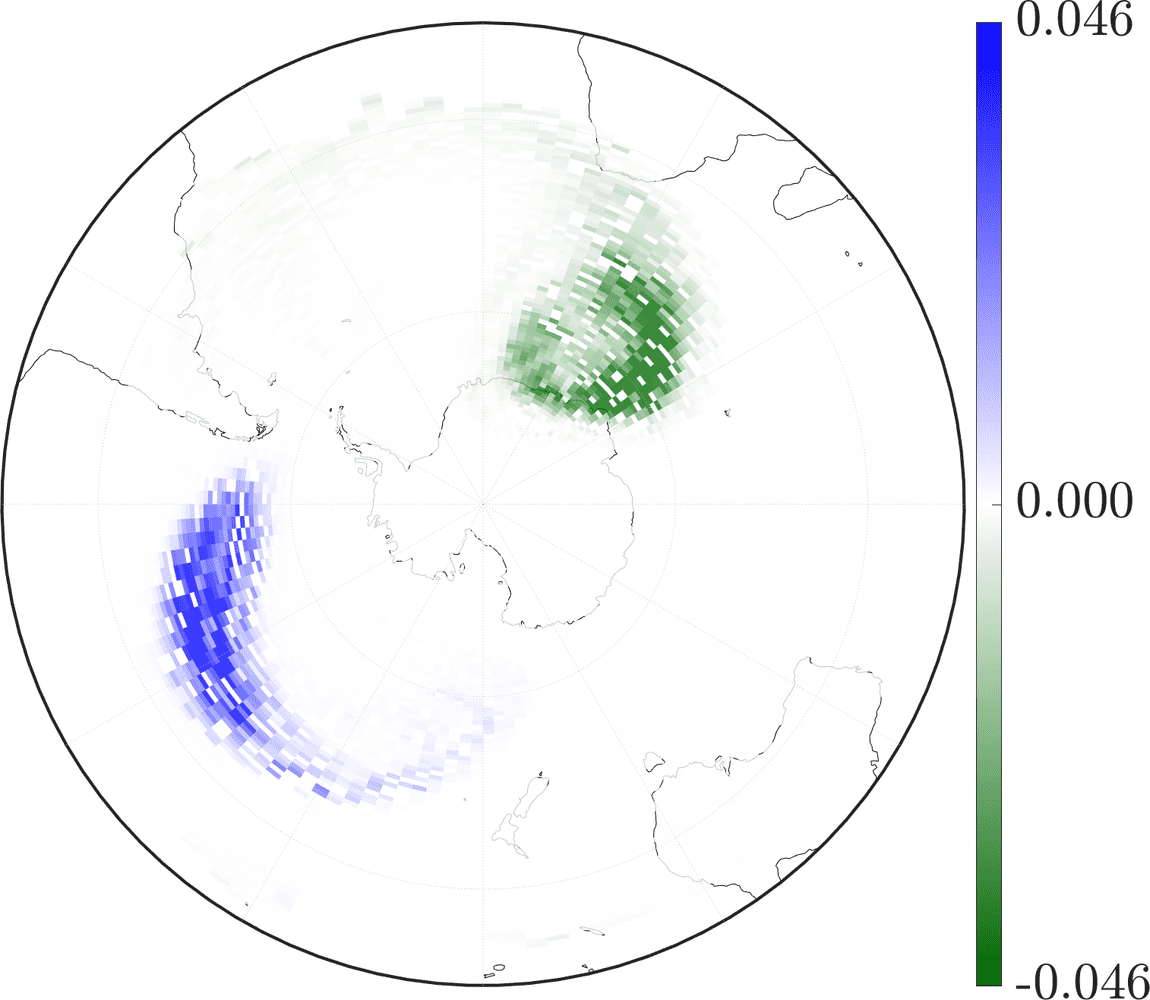}
\subcaption{$1200$ on $27$ Sep.}
\end{minipage}
\caption{Evolved subdominant mode associated with Figure~\ref{fig:SPV_res14_n56_rectMat} for a time window centred on the peak at $0600$ on $24$ Sep. This is illustrated on the area south of $15^{\circ}$S.}
\label{fig:SV2square}
\end{figure}

For singular vectors corresponding to the modes paired through time in Figure~\ref{fig:SPV_res14_n56_sqMat}, it is less clear which windows are appropriate. As such we utilise the pairings of Figure~\ref{fig:SPV_res14_n56_rectMat} and choose to explore the peak that best corresponds to this case. Thus, when exploring the paired modes of Figure~\ref{fig:SPV_res14_n56_sqMat} we choose a time window centred at $1800$ $23$ September. That is, in this case we consider the time window $[1800$ $16$ September, $1800$ $30$ September$]$. Figure~\ref{fig:res14Sq_south30} plots the leading mode for this time window using Algorithm~\ref{alg:Movies}. Comparable results for the time window centred at $0600$ $24$ September, with modes tracked as per Figure~\ref{fig:SPV_res14_n56_rectMat}, are presented in Figure~\ref{fig:SV2square}. To illustrate the precise splitting suggested by our algorithms, we show the subdominant singular vector in the latter case.

\vspace{-0.25cm}
\vspace{-0.25cm}
The noticeable pixellation in Figure~\ref{fig:SV2square} indicates that a non-optimal number of bins have zero mass. This suggests that the mass missing from these bins may be coming from areas beyond the initial seeding. This issue is addressed by extending the seeding to bins with a centre south of $20^{\circ}$S. Rather than consider a new set of rolling windows, we centre this time window mid way between the aforementioned cases, at $0000$ $24$ September.

\begin{figure}[H]
\begin{minipage}[b]{0.245\textwidth}
\centering
\includegraphics[width=\textwidth,center]{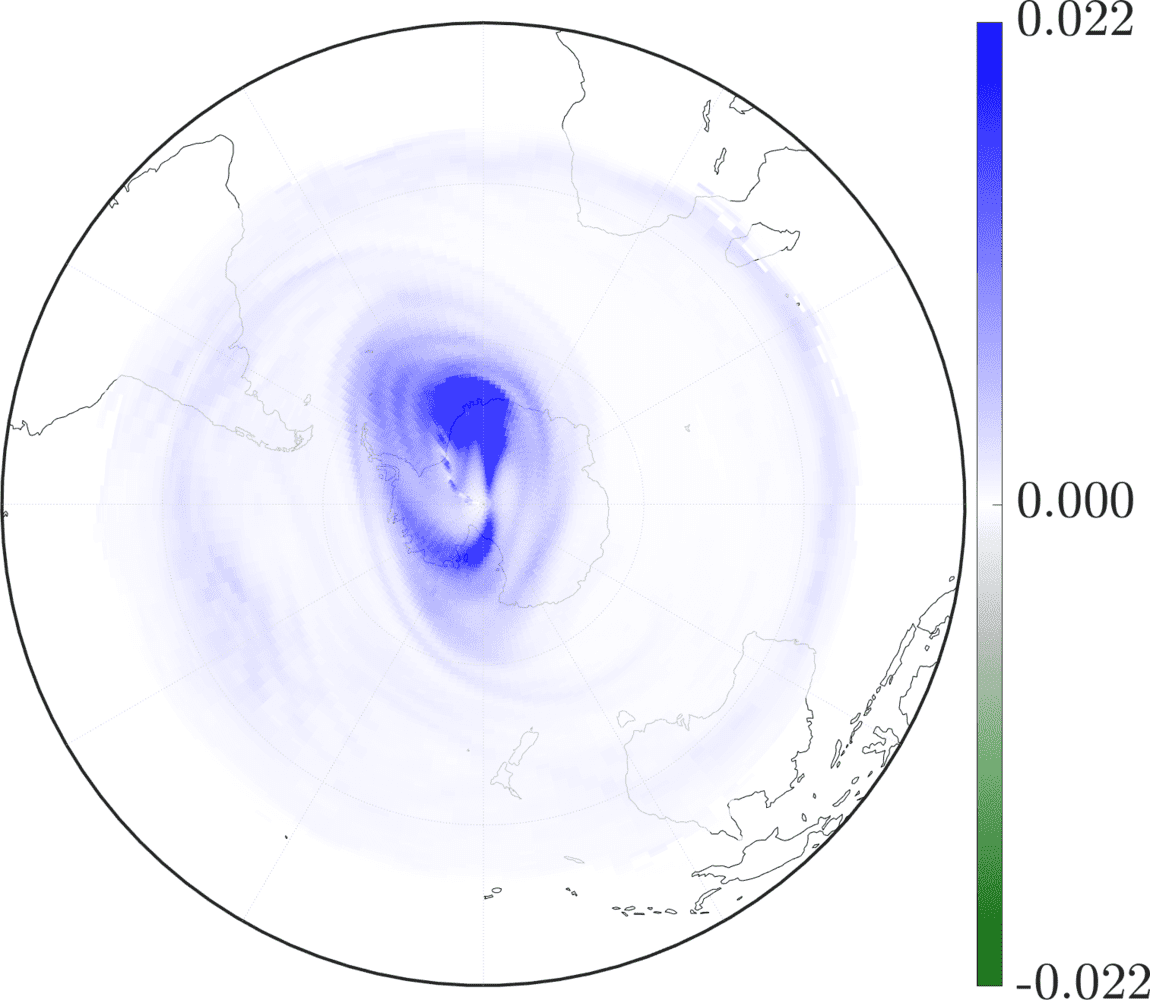}
\subcaption{$1800$ on $19$ Sep.}
\end{minipage}
\begin{minipage}[b]{0.245\textwidth}
\centering
\includegraphics[width=\textwidth,center]{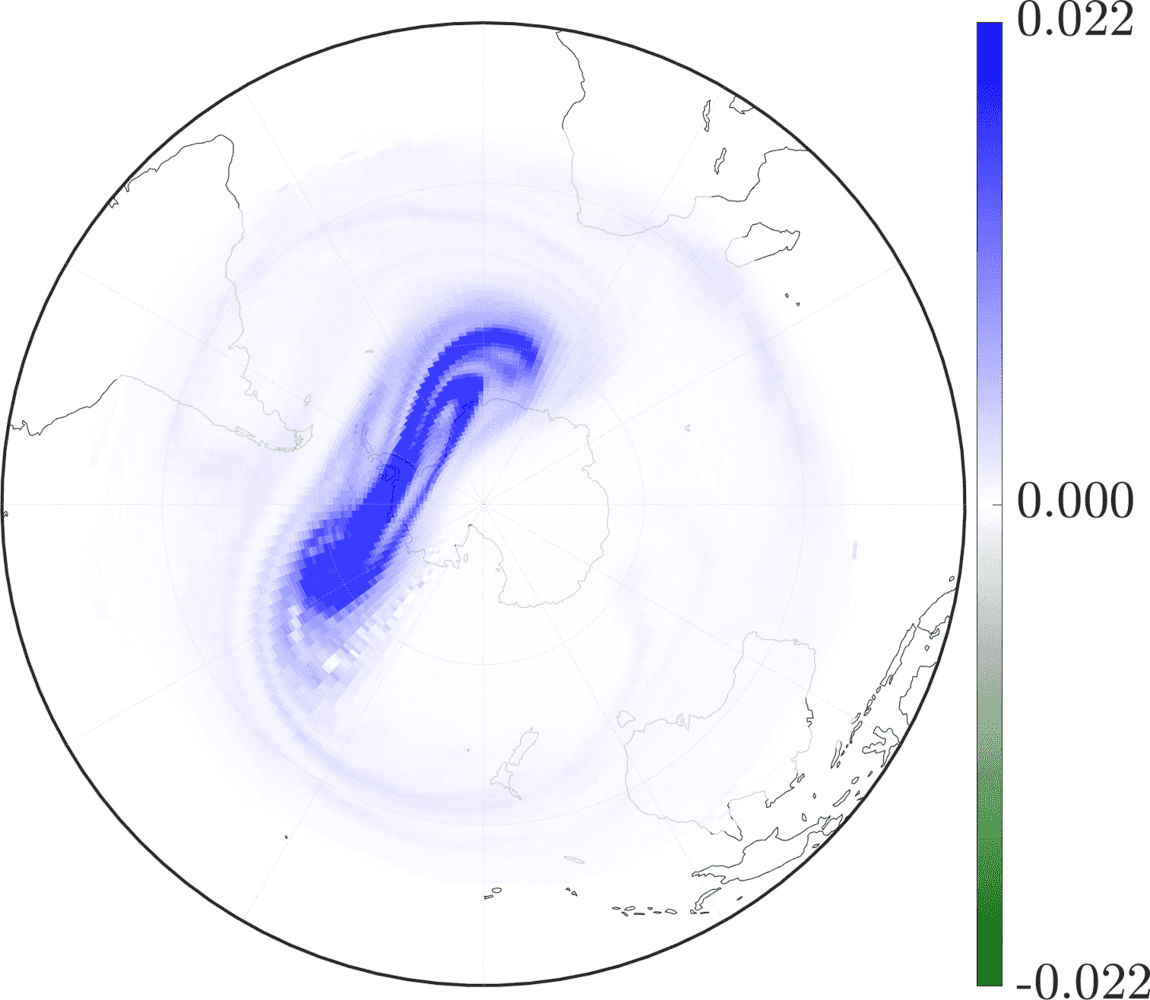}
\subcaption{$1200$ on $22$ Sep.}
\end{minipage}
\begin{minipage}[b]{0.245\textwidth}			
\centering
\includegraphics[width=\textwidth,center]{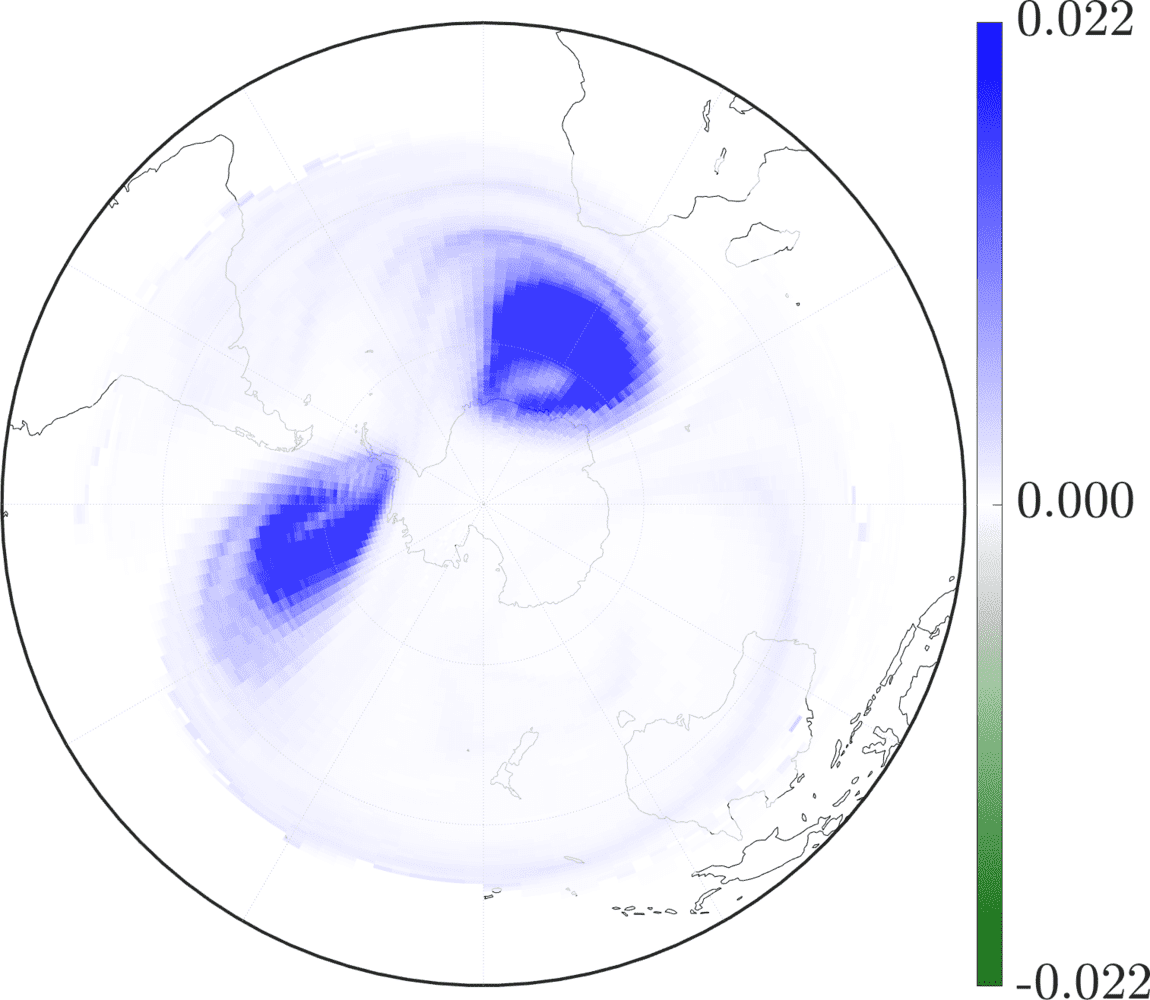}
\subcaption{$1200$ on $25$ Sep.}
\end{minipage}
\begin{minipage}[b]{0.245\textwidth}
\centering
\includegraphics[width=\textwidth,center]{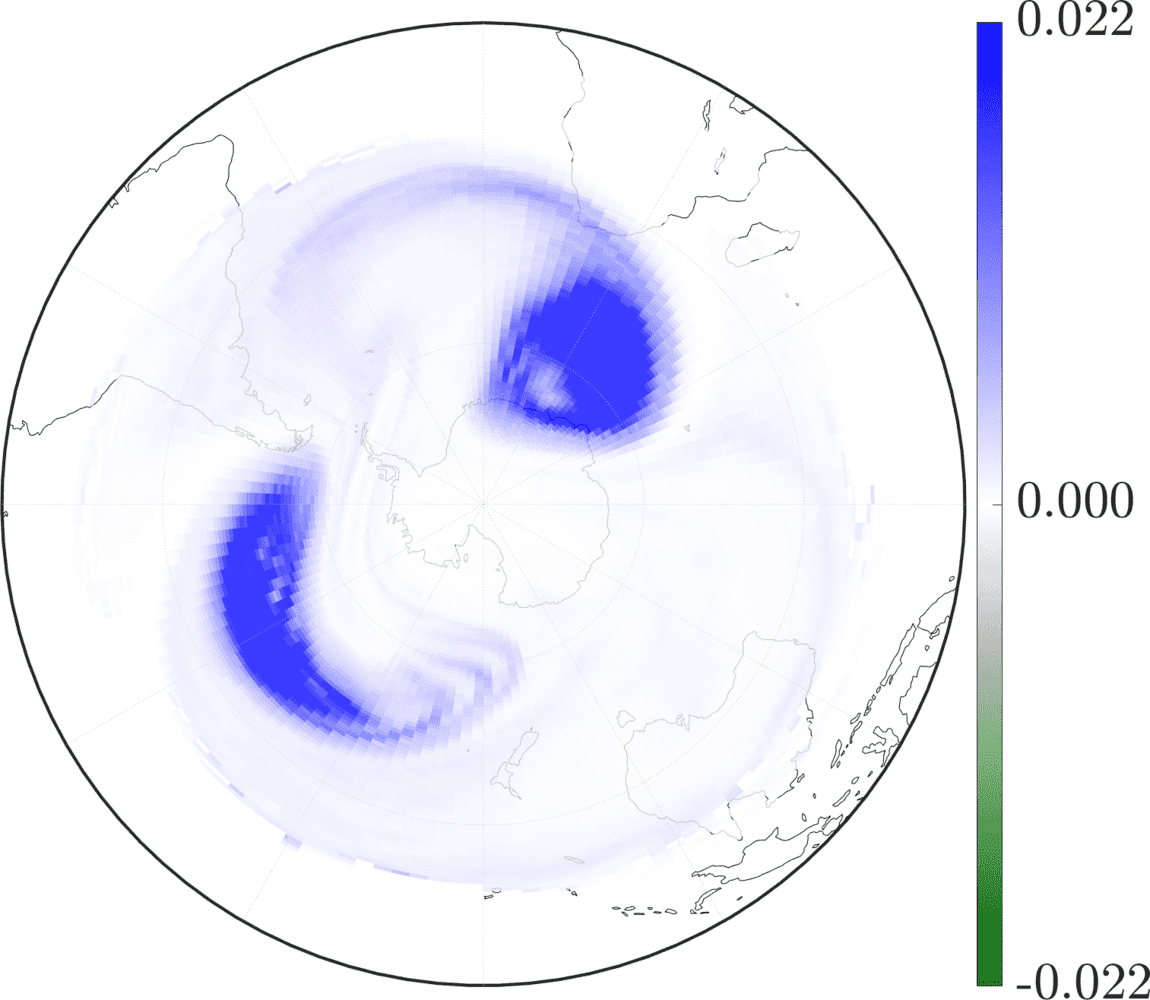}
\subcaption{$1200$ on $27$ Sep.}
\end{minipage}
\caption{Evolved leading singular vectors for time windows centred at $0000$ on $24$ Sep. for $m=12,800$ initially seeded bins whose centres are south of $20^{\circ}$S. This is illustrated on the full southern hemisphere.} 
\label{fig:SPV_res14_n56_lim20}
\end{figure}

In this case we begin with $m=12,800$ bins seeded at $t_{0}=$ $0000$ $17$ September. We close the time window at $0000$ on $1$ October with a collection of $m'=13,604$ bins. Whilst one notes that the wider initial seeding leads to slightly larger matrices, we still have $2^{13}<m,m'<2^{14}$. Thus a smoother distribution of mass throughout the evolved system is achieved without the increase in resolution we would need for the case of $m=m'$. This is evident in the results presented in Figure~\ref{fig:SPV_res14_n56_lim20}. 

\begin{figure}[H]
\begin{minipage}[b]{0.245\textwidth}
\centering
\includegraphics[width=\textwidth,center]{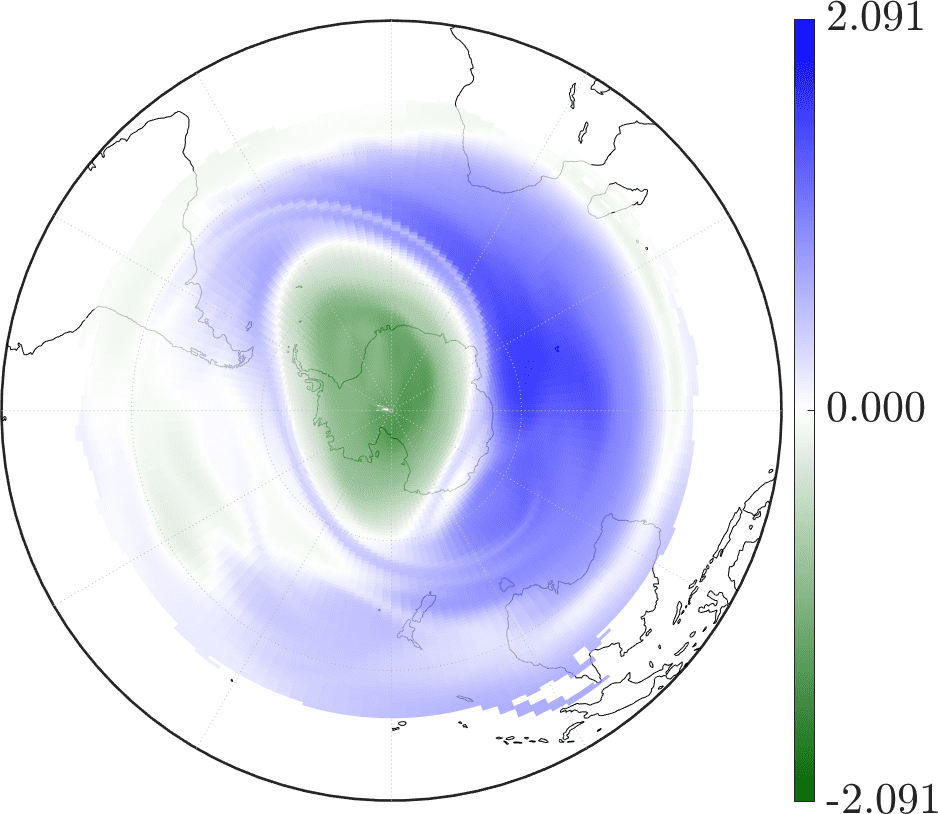}
\subcaption{$1800$ on $19$ Sep.} 
\label{fig:SPV_res14_n56_lim20_a}
\end{minipage}
\begin{minipage}[b]{0.245\textwidth}
\centering
\includegraphics[width=\textwidth,center]{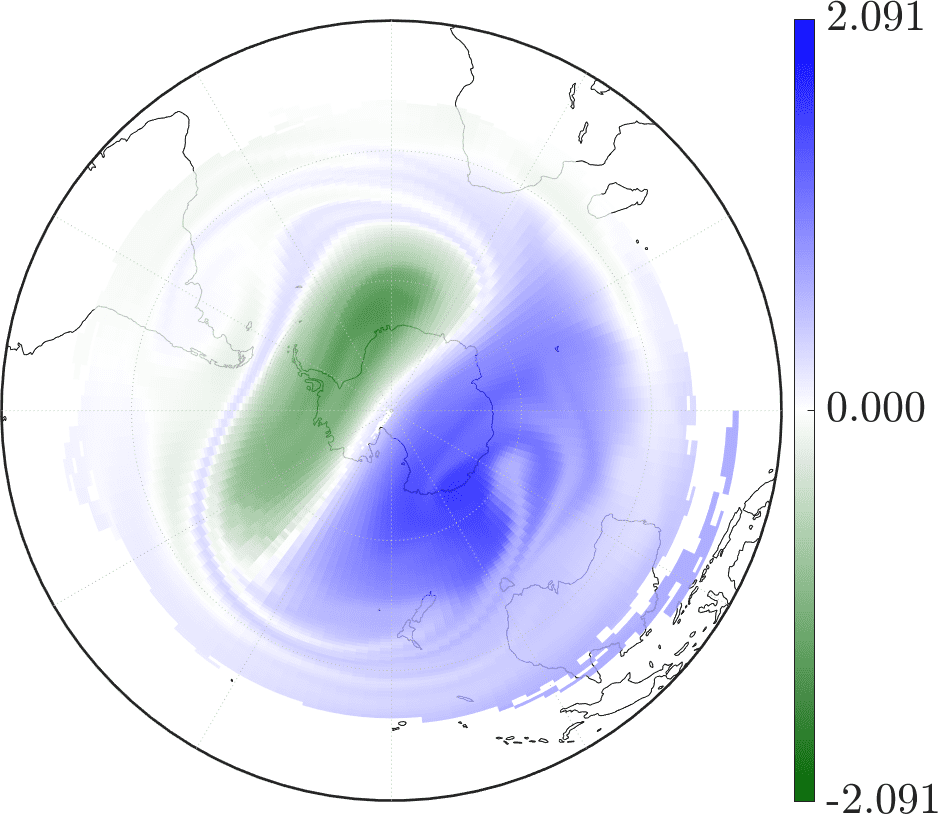}
\subcaption{$1200$ on $22$ Sep.} 
\label{fig:SPV_res14_n56_lim20_b}
\end{minipage}
\begin{minipage}[b]{0.245\textwidth}			
\centering
\includegraphics[width=\textwidth,center]{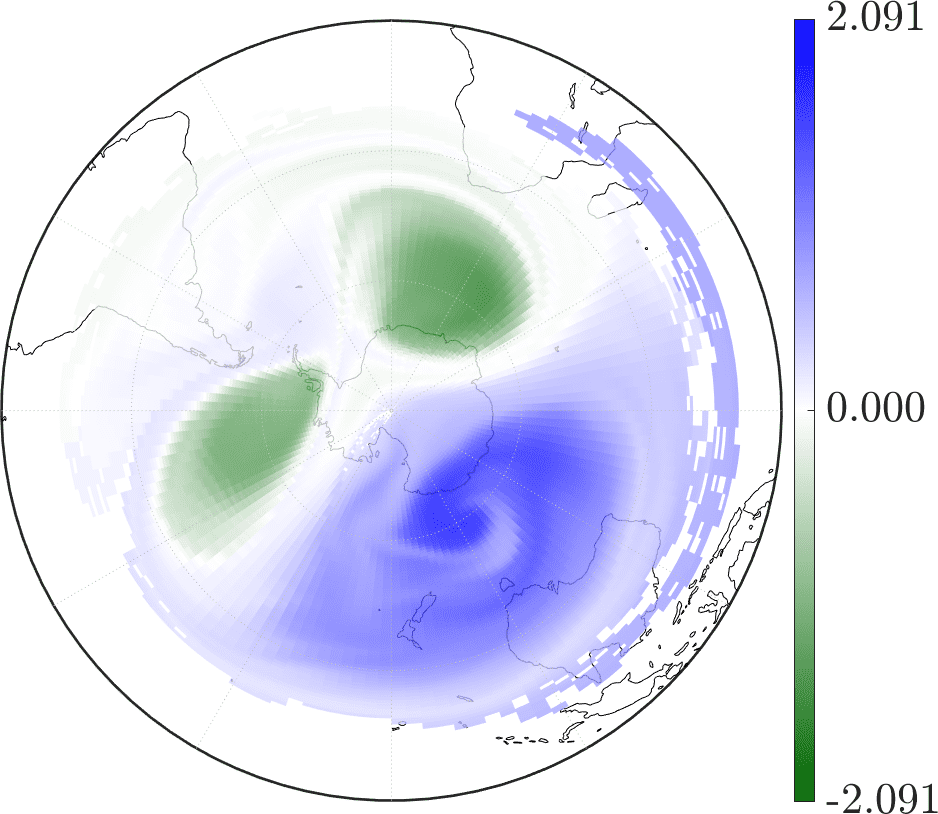}
\subcaption{$1200$ on $25$ Sep.} 
\label{fig:SPV_res14_n56_lim20_c}
\end{minipage}
\begin{minipage}[b]{0.245\textwidth}
\centering
\includegraphics[width=\textwidth,center]{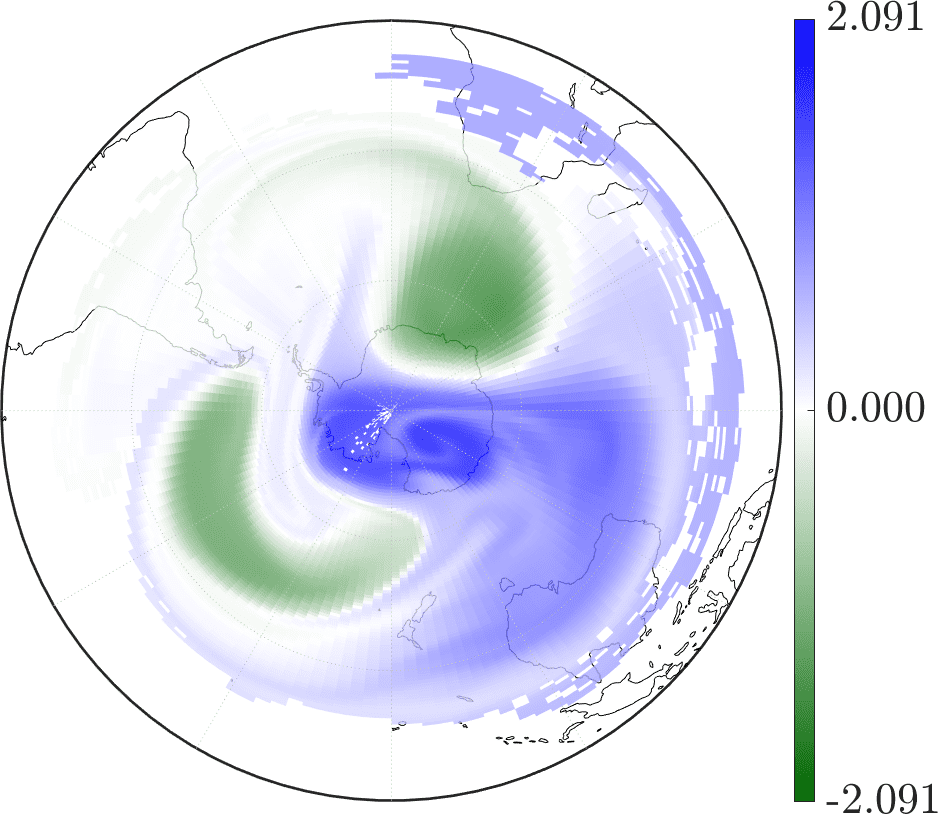}
\subcaption{$1200$ on $27$ Sep.} 
\label{fig:SPV_res14_n56_lim20_d}
\end{minipage}
\caption{Evolved subdominant mode normalised as in \cite{FSM_2010} for time windows centred at $0000$ on $24$ Sep. for $m=12,800$ initially seeded bins whose centres are south of $20^{\circ}$S. This is illustrated on the full southern hemisphere.}
\label{fig:SPV_res14_n56_lim20_FSM}
\end{figure}

\vspace{-0.25cm}
Let us now consider the subdominant singular vector. Normalising this, as per \cite{FSM_2010}, results in Figure~\ref{fig:SPV_res14_n56_lim20_FSM}. The normalised first singular vector now returns a uniform density by design but the subdominant one identifies structures closer to Ertel's PV on the $850$ K isentropic surface as shown in \cite{Charltonetal_SPV}. Whilst in this case the weaker of the two anticyclones is not distinctly identified, the SPV itself is clearly separated from surrounding areas by streams of filaments. This can be seen in Figure~\ref{fig:SPV_res14_n56_lim20_b} which shows the vortex elongated and preconditioned to separate.

\vspace{-0.25cm}
Figure~\ref{fig:SPV_res14_n56_lim20_c} illustrates two distinct daughter vortices.
There is also a separation of space into two clear components (green and blue in the electronic version). The more strongly negative values are associated with the SPV and daughter vortices whilst positive values indicate the stronger anticyclone south of Australia. Figure~\ref{fig:SPV_res14_n56_lim20_d} illustrates how this anticyclone then moves between the daughter vortices to merge with the weaker anticyclone near the tip of South America. 

\section{Conclusion}
Combining the transfer operator technology with existing numerical and data analysis techniques for the purpose of identifying finite-time coherent structures is an active area of research \cite{FroylandRossSakeralliou,KPG_SR_SP_AV_2017,ndour2020predicting,BanischKoltai}. One expects detailed information about the global dynamics of a system can be extracted directly from numerical models, using MET tools and ideas. Indeed, we have found that fundamental changes to coherent structures can be detected using the tools of MET and transfer operators. Our algorithms detected a number of merging and separation events, both in periodically and quasi-periodically driven idealised models and in the real world example of the Southern Polar Vortex.

Events such as merging and splitting not only affect structural boundaries, they are also associated with the expected lifespan of the associated structures. We found that phases of increasing or decreasing singular values of neighbouring matrix cocycles were indicative of fundamental changes in how the boundaries of the associated structures were defined over time. However, changes in dominance, associated with crossings, were related to the complementary relationship between modes associated with different aspects of the same structure experiencing varied coherency through time.

In the double well potential model, it was clear that trajectories of increasing singular values were related to an increasingly coherent structure working to consolidate external boundaries. This held whether the identified structure was characterised as a single core component or as consisting of two clearly distinct cores of opposite sign. The converse was true for trajectories of decreasing dominance. Trajectories of decreasing singular values were associated with structures whose boundaries would disintegrate. 

The efficacy of our pairing algorithms was assessed using a measure of equivariance mismatch. Interestingly, these values were lowest when time windows were longer. This occurred despite the fact that more matrices were being multiplied and numerical error was likely to increase. This is in alignment with the outcomes of the multiplicative ergodic theorem. As time increases, we see a more distinct separation of structures based on their expected survival rates. The less coherent, more short lived structures, will mix more freely with their surroundings and lose dynamical significance. 

Whilst it becomes progressively more difficult to track modes in increasingly complex examples, our algorithms could still identify the splitting of the Southern Polar Vortex. In all cases the associated singular vectors were useful for identifying the spatial region where fundamental changes occurred. However, further research is required to assess the suitability of these methods to wider applications. 

Likewise, further research is needed to clarify exactly how the onset of fundamental changes is signalled in geophysical models. One might also be interested to explore how to optimise time window lengths or bin size for a variety of models and dynamical behaviour. This future work should not treat these methods as a mere black box approach. Indeed, we anticipate effective work in this direction to incorporate specific disciplinary knowledge alongside ideas from ergodic theory. 

\section*{Acknowledgments} The authors would like to acknowledge Sanjeeva Balasuriya (Adelaide), Gary Froyland (UNSW) and Tony Roberts (UQ) for useful conversations and suggestions, and two anonymous referees for providing thoughtful and constructive feedback.


\begin{thebibliography}{10}

\bibitem{AllshouseMichaelR.2015Lbmf}
\newblock M.~R. Allshouse and T.~Peacock,
\newblock Lagrangian based methods for coherent structure detection,
\newblock \emph{Chaos}, \textbf{25} (2015), 097617--1, 097617--13,
\newblock \urlprefix\url{https://doi.org/10.1063/1.4922968}.

\bibitem{BalasuriyaSanjeeva2018GLcs}
\newblock S.~Balasuriya, N.~T. Ouellette and I.~I. Rypina,
\newblock Generalized {L}agrangian coherent structures,
\newblock \emph{Physica D: Nonlinear Phenomena}, \textbf{372} (2018), 31--51.

\bibitem{BanischKoltai}
\newblock R.~Banisch and P.~Koltai,
\newblock {Understanding the geometry of transport: Diffusion maps for
  Lagrangian trajectory data unravel coherent sets},
\newblock \emph{Chaos: An Interdisciplinary Journal of Nonlinear Science},
  \textbf{27} (2017), 035804,
\newblock \urlprefix\url{https://doi.org/10.1063/1.4971788}.

\bibitem{BudisicMarko2012AK}
\newblock M.~Budi{\v{s}}i{\'{c}}, R.~Mohr and I.~Mezi{\'{c}},
\newblock Applied {K}oopmanism,
\newblock \emph{Chaos}, \textbf{22} (2012), 047510, 33,
\newblock \urlprefix\url{https://doi.org/10.1063/1.4772195}.

\bibitem{Charltonetal_SPV}
\newblock A.~J. Charlton, A.~O’Neill, W.~A. Lahoz and P.~Berrisford,
\newblock The splitting of the stratospheric polar vortex in the {S}outhern
  {H}emisphere, {S}eptember 2002: Dynamical evolution,
\newblock \emph{Journal of the Atmospheric Sciences}, \textbf{62} (2005),
  590--602,
\newblock \urlprefix\url{https://doi.org/10.1175/JAS-3318.1}.

\bibitem{DeeD2011TErc}
\newblock D.~Dee, S.~Uppala, A.~Simmons, P.~Berrisford, P.~Poli, S.~Kobayashi,
  U.~Andrae, M.~Balmaseda, G.~Balsamo, P.~Bauer, P.~Bechtold, A.~Beljaars,
  L.~van~de Berg, J.~Bidlot, N.~Bormann, C.~Delsol, R.~Dragani, M.~Fuentes,
  A.~Geer and L.~Haimberger,
\newblock The {ERA}-{I}nterim reanalysis: configuration and performance of the
  data assimilation system,
\newblock \emph{Quarterly Journal Of The Royal Meteorological Society},
  \textbf{137} (2011), 553--597.

\bibitem{DellnitzMichael2009OtAo}
\newblock M.~Dellnitz, G.~Froyland, C.~Horenkamp and K.~Padberg,
\newblock On the approximation of transport phenomena---a dynamical systems
  approach,
\newblock \emph{GAMM-Mitt.}, \textbf{32} (2009), 47--60,
\newblock \urlprefix\url{https://doi.org/10.1002/gamm.200910004}.

\bibitem{Dellnitz2001}
\newblock M.~Dellnitz, G.~Froyland and O.~Junge,
\newblock The algorithms behind {GAIO}-set oriented numerical methods for
  dynamical systems,
\newblock in \emph{Ergodic theory, analysis, and efficient simulation of
  dynamical systems},
\newblock Springer, Berlin, 2001,
\newblock 145--174, 805--807.

\bibitem{DellnitzMichael1999OtAo}
\newblock M.~Dellnitz and O.~Junge,
\newblock On the approximation of complicated dynamical behavior,
\newblock \emph{SIAM J. Numer. Anal.}, \textbf{36} (1999), 491--515,
\newblock \urlprefix\url{https://doi.org/10.1137/S0036142996313002}.

\bibitem{DeuflhardDellnitzJungeetal.1999}
\newblock P.~Deuflhard, M.~Dellnitz, O.~Junge and C.~Sch{\"u}tte,
\newblock Computation of essential molecular dynamics by subdivision
  techniques,
\newblock in \emph{Computational Molecular Dynamics: Challenges, Methods,
  Ideas: Proceedings of the 2nd International Symposium on Algorithms for
  Macromolecular Modelling, Berlin, May 21–24, 1997} (eds. P.~Deuflhard,
  J.~Hermans, B.~Leimkuhler, A.~E. Mark, S.~Reich and R.~D. Skeel), vol.~4 of
  Lecture Notes in Computational Science and Engineering,,
\newblock Springer Berlin Heidelberg, Berlin, Heidelberg, 1999,
\newblock 98 -- 115.

\bibitem{FroylandGary2013CcLv}
\newblock G.~Froyland, T.~H{\"{u}}ls, G.~P. Morriss and T.~M. Watson,
\newblock Computing covariant {L}yapunov vectors, {O}seledets vectors, and
  dichotomy projectors: a comparative numerical study,
\newblock \emph{Phys. D}, \textbf{247} (2013), 18--39,
\newblock \urlprefix\url{https://doi.org/10.1016/j.physd.2012.12.005}.

\bibitem{FroylandGary2010Csai}
\newblock G.~Froyland, S.~Lloyd and A.~Quas,
\newblock Coherent structures and isolated spectrum for {P}erron-{F}robenius
  cocycles,
\newblock \emph{Ergodic Theory Dynam. Systems}, \textbf{30} (2010), 729--756,
\newblock \urlprefix\url{https://doi.org/10.1017/S0143385709000339}.

\bibitem{FLQ2}
\newblock G.~Froyland, S.~Lloyd and A.~Quas,
\newblock A semi-invertible {O}seledets theorem with applications to transfer
  operator cocycles,
\newblock \emph{Discrete Contin. Dyn. Syst.}, \textbf{33} (2013), 3835--3860,
\newblock \urlprefix\url{http://dx.doi.org/10.3934/dcds.2013.33.3835}.

\bibitem{FLS_2010}
\newblock G.~Froyland, S.~Lloyd and N.~Santitissadeekorn,
\newblock Coherent sets for nonautonomous dynamical systems,
\newblock \emph{Phys. D}, \textbf{239} (2010), 1527--1541,
\newblock \urlprefix\url{https://doi.org/10.1016/j.physd.2010.03.009}.

\bibitem{FP_2009}
\newblock G.~Froyland and K.~Padberg,
\newblock Almost-invariant sets and invariant manifolds---connecting
  probabilistic and geometric descriptions of coherent structures in flows,
\newblock \emph{Phys. D}, \textbf{238} (2009), 1507--1523,
\newblock \urlprefix\url{https://doi.org/10.1016/j.physd.2009.03.002}.

\bibitem{FroylandGary2007Doco}
\newblock G.~Froyland, K.~Padberg, M.~H. England and A.~M. Treguier,
\newblock Detection of coherent oceanic structures via transfer operators.,
\newblock \emph{Physical review letters}, \textbf{98} (2007),
  224503--1,224503--4,
\newblock \urlprefix\url{http://search.proquest.com/docview/68131950/}.

\bibitem{FroyPG2014}
\newblock G.~Froyland and K.~Padberg-Gehle,
\newblock Almost-invariant and finite-time coherent sets: directionality,
  duration, and diffusion,
\newblock in \emph{Ergodic theory, open dynamics, and coherent structures},
  vol.~70 of Springer Proc. Math. Stat.,
\newblock Springer, New York, 2014,
\newblock 171--216,
\newblock \urlprefix\url{https://doi.org/10.1007/978-1-4939-0419-8_9}.

\bibitem{FroylandRossSakeralliou}
\newblock G.~Froyland, C.~P. Rock and K.~Sakellariou,
\newblock Sparse eigenbasis approximation: Multiple feature extraction across
  spatiotemporal scales with application to coherent set identification,
\newblock \emph{Communications in Nonlinear Science and Numerical Simulation},
  \textbf{77} (2019), 81 -- 107,
\newblock
  \urlprefix\url{http://www.sciencedirect.com/science/article/pii/S1007570419301236}.

\bibitem{FSM_2010}
\newblock G.~Froyland, N.~Santitissadeekorn and A.~Monahan,
\newblock Transport in time-dependent dynamical systems: finite-time coherent
  sets,
\newblock \emph{Chaos}, \textbf{20} (2010), 043116--1, 043116--10,
\newblock \urlprefix\url{https://doi.org/10.1063/1.3502450}.

\bibitem{GinelliEtAl}
\newblock F.~Ginelli, P.~Poggi, A.~Turchi, H.~Chat\'e, R.~Livi and A.~Politi,
\newblock Characterizing dynamics with covariant lyapunov vectors,
\newblock \emph{Phys. Rev. Lett.}, \textbf{99} (2007), 130601,
\newblock
  \urlprefix\url{https://link.aps.org/doi/10.1103/PhysRevLett.99.130601}.

\bibitem{GonzlezTokman2018MultiplicativeET}
\newblock C.~Gonz{\'a}lez-Tokman,
\newblock Multiplicative ergodic theorems for transfer operators: towards the
  identification and analysis of coherent structures in non-autonomous
  dynamical systems,
\newblock in \emph{Contributions of {M}exican mathematicians abroad in pure and
  applied mathematics}, vol. 709 of Contemp. Math.,
\newblock Amer. Math. Soc., Providence, RI, 2018,
\newblock 31--52,
\newblock \urlprefix\url{https://doi.org/10.1090/conm/709/14290}.

\bibitem{GTQuas1}
\newblock C.~Gonz{\'a}lez-Tokman and A.~Quas,
\newblock A semi-invertible operator {O}seledets theorem,
\newblock \emph{Ergodic Theory Dynam. Systems}, \textbf{34} (2014), 1230--1272,
\newblock \urlprefix\url{https://doi.org/10.1017/etds.2012.189}.

\bibitem{GTQuas2}
\newblock C.~Gonz{\'a}lez-Tokman and A.~Quas,
\newblock A concise proof of the multiplicative ergodic theorem on {B}anach
  spaces,
\newblock \emph{J. Mod. Dyn.}, \textbf{9} (2015), 237--255,
\newblock \urlprefix\url{https://doi.org/10.3934/jmd.2015.9.237}.

\bibitem{HallerGeorge2015LCS}
\newblock G.~Haller,
\newblock Lagrangian coherent structures,
\newblock in \emph{Annual review of fluid mechanics. {V}ol. 47}, vol.~47 of
  Annu. Rev. Fluid Mech.,
\newblock Annual Reviews, Palo Alto, CA, 2015,
\newblock 137--162.

\bibitem{HallerKarraschKogelbauer}
\newblock G.~Haller, D.~Karrasch and F.~Kogelbauer,
\newblock Material barriers to diffusive and stochastic transport,
\newblock \emph{Proceedings of the National Academy of Sciences}, \textbf{115}
  (2018), 9074--9079,
\newblock \urlprefix\url{https://www.pnas.org/content/115/37/9074}.

\bibitem{BernardBinson2002}
\newblock B.~Joseph and B.~Legras,
\newblock Relation between kinematic boundaries, stirring, and barriers for the
  {A}ntarctic polar vortex,
\newblock \emph{Journal of the Atmospheric Sciences}, \textbf{59} (2002),
  1198--1212,
\newblock
  \urlprefix\url{https://doi.org/10.1175/1520-0469(2002)059<1198:RBKBSA>2.0.CO;2}.

\bibitem{KlusStefan2015Otna}
\newblock S.~Klus, P.~Koltai and C.~Sch\"{u}tte,
\newblock On the numerical approximation of the {P}erron-{F}robenius and
  {K}oopman operator,
\newblock \emph{J. Comput. Dyn.}, \textbf{3} (2016), 51--79,
\newblock \urlprefix\url{https://doi.org/10.3934/jcd.2016003}.

\bibitem{KoltaiRenger}
\newblock P.~Koltai and D.~R.~M. Renger,
\newblock From large deviations to semidistances of transport and mixing:
  Coherence analysis for finite {L}agrangian data,
\newblock \emph{Journal of Nonlinear Science}, \textbf{28} (2018), 1915--1957,
\newblock \urlprefix\url{https://doi.org/10.1007/s00332-018-9471-0}.

\bibitem{LasotaAndrzej1994CFaN}
\newblock A.~Lasota and M.~C. Mackey,
\newblock \emph{Chaos, fractals, and noise : stochastic aspects of dynamics /
  {A}ndrzej {L}asota, {M}ichael {C}. {M}ackey.},
\newblock 2nd edition,
\newblock Applied mathematical sciences (Springer-Verlag New York Inc.); v. 97,
  Springer-Verlag, New York, 1994.

\bibitem{LekienRoss}
\newblock F.~Lekien and S.~D. Ross,
\newblock The computation of finite-time {L}yapunov exponents on unstructured
  meshes and for non-{E}uclidean manifolds,
\newblock \emph{Chaos: An Interdisciplinary Journal of Nonlinear Science},
  \textbf{20} (2010), 017505,
\newblock \urlprefix\url{https://doi.org/10.1063/1.3278516}.

\bibitem{MosovskyBA2011TiTD}
\newblock B.~A. Mosovsky and J.~D. Meiss,
\newblock Transport in transitory dynamical systems,
\newblock \emph{SIAM Journal on Applied Dynamical Systems}, \textbf{10} (2011),
  35--65.

\bibitem{ndour2020predicting}
\newblock M.~Ndour and K.~Padberg-Gehle,
\newblock Predicting bifurcations of almost-invariant patterns: a set-oriented
  approach, 2020,
\newblock ArXiv:2001.01099 [math.DS].

\bibitem{NewmanPaul2005TUSH}
\newblock P.~Newman and E.~Nash,
\newblock The unusual {S}outhern {H}emisphere stratosphere winter of 2002,
\newblock \emph{Journal of the Atmospheric Sciences}, \textbf{62} (2005),
  614--628,
\newblock \urlprefix\url{http://search.proquest.com/docview/20651113/}.

\bibitem{Noethen}
\newblock F.~Noethen,
\newblock A projector-based convergence proof of the {G}inelli algorithm for
  covariant {L}yapunov vectors,
\newblock \emph{Phys. D}, \textbf{396} (2019), 18--34,
\newblock \urlprefix\url{https://doi.org/10.1016/j.physd.2019.02.012}.

\bibitem{oNeill2017vortex}
\newblock A.~O'Neill, C.~L. Oatley, A.~J. Charlton‐Perez, D.~M. Mitchell and
  T.~Jung,
\newblock Vortex splitting on a planetary scale in the stratosphere by
  cyclogenesis on a subplanetary scale in the troposphere,
\newblock \emph{Quarterly Journal of the Royal Meteorological Society},
  \textbf{143} (2017), 691--705,
\newblock \urlprefix\url{https://doi.org/10.1002/qj.2957}.

\bibitem{OrsoliniYvanJ.2005Aoso}
\newblock Y.~Orsolini, R.~C.E, G.~Manney and A.~D.R,
\newblock An observational study of the final breakdown of the {S}outhern
  {H}emisphere stratospheric vortex in 2002,
\newblock \emph{Journal of the Atmospheric Sciences}, \textbf{62} (2005),
  735--747.

\bibitem{MR0240280}
\newblock V.~I. Oseledec,
\newblock A multiplicative ergodic theorem. {C}haracteristic {L}japunov,
  exponents of dynamical systems,
\newblock \emph{Trudy Moskov. Mat. Ob{\v s}{\v c}.}, \textbf{19} (1968),
  179--210.

\bibitem{KPG_SR_SP_AV_2017}
\newblock K.~Padberg-Gehle, S.~Reuther, S.~Praetorius and A.~Voigt,
\newblock Transfer operator-based extraction of coherent features on surfaces,
\newblock in \emph{Topological Methods in Data Analysis and Visualization IV}
  (eds. H.~Carr, C.~Garth and T.~Weinkauf),
\newblock Springer International Publishing, Cham, 2017,
\newblock 283--297.

\bibitem{RaghunathanM.1979ApoO}
\newblock M.~S. Raghunathan,
\newblock A proof of {O}seledec's multiplicative ergodic theorem,
\newblock \emph{Israel J. Math.}, \textbf{32} (1979), 356--362,
\newblock \urlprefix\url{https://doi.org/10.1007/BF02760464}.

\bibitem{ShaddenShawnC2005Dapo}
\newblock S.~C. Shadden, F.~Lekien and J.~E. Marsden,
\newblock Definition and properties of {L}agrangian coherent structures from
  finite-time {L}yapunov exponents in two-dimensional aperiodic flows,
\newblock \emph{Physica D: Nonlinear Phenomena}, \textbf{212} (2005), 271--304.

\bibitem{UlamStanislawM.1960Acom}
\newblock S.~M. Ulam,
\newblock \emph{A collection of mathematical problems},
\newblock Interscience Tracts in Pure and Applied Mathematics, {N}o. 8,
  Interscience Publishers, New York-London, 1960.

\bibitem{WilliamsKevrekidisRowley2015}
\newblock M.~O. Williams, I.~G. Kevrekidis and C.~W. Rowley,
\newblock A data-driven approximation of the {K}oopman operator: extending
  dynamic mode decomposition,
\newblock \emph{J. Nonlinear Sci.}, \textbf{25} (2015), 1307--1346,
\newblock \urlprefix\url{https://doi.org/10.1007/s00332-015-9258-5}.

\end{thebibliography}
\providecommand{\href}[2]{#2}
\providecommand{\arxiv}[1]{\href{http://arxiv.org/abs/#1}{arXiv:#1}}
\providecommand{\url}[1]{\texttt{#1}}
\providecommand{\urlprefix}{URL }

\medskip
Received xxxx 20xx; revised xxxx 20xx.
\medskip

\end{document}